\documentstyle[10pt]{amsart}

\theoremstyle{definition}
\newtheorem{thm}{Theorem}[section]
\newtheorem{lem}[thm]{Lemma}
\newtheorem{prp}[thm]{Proposition}
\newtheorem{dfn}[thm]{Definition}
\newtheorem{cor}[thm]{Corollary}

\newtheorem{rmk}[thm]{Remark}
\newtheorem{ntn}[thm]{Notation}
\newtheorem{exa}[thm]{Example}

\newenvironment{pff}{{\em Proof:}}{\QED}

\newcommand{\beq}{\begin{equation}}
\newcommand{\eeq}{\end{equation}}
\newcommand{\beqr}{\begin{eqnarray*}}
\newcommand{\eeqr}{\end{eqnarray*}}

\newcommand{\dirlim}{\displaystyle \lim_{\longrightarrow}}
\newcommand{\limi}[1]{\lim_{{#1} \to \infty}}
\newcommand{\limo}[1]{\lim_{{#1} \to \om}}

\newcommand{\diag}{\operatorname{diag}}
\newcommand{\ev}{\operatorname{ev}}

\newcommand{\Aut}{\operatorname{Aut}}
\newcommand{\id}{\operatorname{id}}

\newcommand{\re}{\operatorname{Re}}
\newcommand{\cb}{\operatorname{cb}}
\newcommand{\spn}{\operatorname{span}}
\newcommand{\sgn}{\operatorname{sgn}}
\newcommand{\Lip}{\operatorname{Lip}}
\newcommand{\spec}{\operatorname{sp}}
\newcommand{\Real}{\operatorname{Re}}

\newcommand{\xfrac}[2]{{\textstyle{\frac{#1}{#2}}}}

\newcommand{\half}{\frac{1}{2}}
\newcommand{\QED}{\rule{1.5mm}{3mm}}

\newcommand{\ds}{d_{\operatorname{S}}}
\newcommand{\C}{{\Bbb C}}
\newcommand{\R}{{\Bbb R}}
\newcommand{\Z}{{\Bbb Z}}
\newcommand{\N}{{\Bbb N}}
\newcommand{\unit}{^{\dagger}}

\newcommand{\af}{\alpha}
\newcommand{\bt}{\beta}
\newcommand{\gm}{\gamma}
\newcommand{\dt}{\delta}
\newcommand{\ep}{\varepsilon}
\newcommand{\zt}{\zeta}
\newcommand{\et}{\eta}

\newcommand{\io}{\iota}
\newcommand{\te}{\theta}
\newcommand{\ld}{\lambda}
\newcommand{\sm}{\sigma}

\newcommand{\ph}{\varphi}
\newcommand{\ps}{\psi}
\newcommand{\rh}{\rho}
\newcommand{\om}{\omega}
\newcommand{\ta}{\tau}

\newcommand{\Gm}{\Gamma}

\newcommand{\ca}{C*-algebra}
\newcommand{\mops}{mutually orthogonal projections}
\newcommand{\hm}{homomorphism}
\newcommand{\pisca}{purely infinite simple \ca}
\newcommand{\sep}{separable}
\newcommand{\ct}{continuous}
\newcommand{\mvn}{Murray-von Neumann equivalent}

\newcommand{\cpt}{compact Hausdorff}
\newcommand{\wolog}{without loss of generality}
\newcommand{\Wolog}{Without loss of generality}

\newcommand{\ifo}{if and only if}
\newcommand{\pj}{projection}
\newcommand{\fd}{finite dimensional}
\newcommand{\andeqn}{\,\,\,\,\,\, {\mathrm{and}} \,\,\,\,\,\,}

\newcommand{\aue}{approximate unitary equivalence}
\newcommand{\ayue}{approximately unitarily equivalent}
\newcommand{\kfalg}{separable nuclear unital \pisca}
\newcommand{\snus}{separable, nuclear, unital, and simple}
\newcommand{\kfive}
        {separable, nuclear, unital, purely infinite, and simple}

\newcommand{\ucp}{unital completely positive}
\newcommand{\uacp}{unital and completely positive}

\newcommand{\OI}{{\cal O}_{\infty}}
\newcommand{\OIA}[1]{\OI \otimes {#1}}
\newcommand{\OA}[1]{{\cal O}_{#1}}
\newcommand{\Mt}[1]{M_{#1} \otimes}
\newcommand{\Kt}{K \otimes}

\newcommand{\OT}[1]{\OA{2} \otimes {#1}}
\newcommand{\tmin}{\otimes_{\operatorname{min}}}
\newcommand{\tmax}{\otimes_{\operatorname{max}}}
\newcommand{\talg}{\otimes_{\operatorname{alg}}}

\newcommand{\LI}[1]{l^{\infty} ({#1})}
\newcommand{\SCC}{\bt \N - \N}

\title[Embedding of exact C*-algebras]{Embedding of exact
C*-algebras and continuous fields in the Cuntz algebra $\OA{2}$}

\date{\today}

\author[Kirchberg]{Eberhard Kirchberg}
\address[]{Institut f\"{u}r Mathematik \\
Humboldt-Universit\"{a}t zu Berlin \\
Unter den Linden 6 \\
D-10099 Berlin  \\
Germany}
\email[]{kirchbrg@@mathematik.hu-berlin.de}

\author[Phillips]{N. Christopher Phillips}
\address[]{Department of Mathematics \\
University of Oregon \\
Eugene OR 97403-1222  \\
USA}
\email[]{phillips@@math.uoregon.edu}
\thanks{Research partially supported by NSF grant DMS-94 00904}



\begin{document}

\maketitle
\setcounter{section}{-1}

\begin{center}
{\em Dedicated to Edward Effros on his 60th birthday.}
\end{center}

\begin{abstract}

We prove that any separable exact \ca\  is isomorphic to a
subalgebra of the Cuntz algebra $\OA{2}.$  We further prove that
if $A$ is a simple separable unital
nuclear \ca, then $\OT{A} \cong \OA{2},$ and if, in addition,
$A$ is purely infinite, then $\OIA{A} \cong A.$

The embedding of exact \ca s in $\OA{2}$ is continuous in the
following sense. If $A$ is a
continuous field of \ca s over a compact manifold or finite
CW complex $X$ with fiber $A (x)$ over $x \in X,$
such that the algebra of continuous
sections of $A$ is separable and exact, then there is a family
of injective \hm s $\ph_x  : A (x) \to \OA{2}$ such that
for every continuous section $a$ of $A$ the function
$x \mapsto \ph_x (a (x))$ is continuous. Moreover, one can say
something about the modulus of continuity of the functions
$x \mapsto \ph_x (a (x))$ in terms of the structure of the
continuous field. In particular, we show that the continuous field
$\te \mapsto A_{\te}$ of rotation algebras posesses unital
embeddings $\ph_{\te}$ in $\OA{2}$ such that  the standard generators
$u (\te)$ and $v (\te)$ satisfy
\[
\max ( \| \ph_{\te_1} (u (\te_1)) - \ph_{\te_2} (u (\te_2)) \|, \,
       \| \ph_{\te_1} (v (\te_1)) - \ph_{\te_2} (v (\te_2)) \| )
  < C | \te_1 - \te_2 |^{1/2}
\]
for some constant $C.$
\end{abstract}

\section{Introduction}

It has recently become clear that the exact \ca s form an important
class of \ca s more general than the nuclear \ca s. (A \ca\  $A$
is called {\em exact} if the functor $A \tmin -$ preserves short
exact sequences.) For example, every C*-subalgebra of a nuclear
\ca\  is exact. The class of exact \ca s has a number of good
functorial properties (see Section 7 of \cite{Kr2}); in particular,
unlike the class of nuclear \ca s, it is closed under passage to
subalgebras. (Unfortunately, though, it is not closed under arbitrary
extensions, only under ``locally liftable'' ones. See \cite{Kr0.5}
and Section 7 of \cite{Kr2}.) The reduced \ca s
of some discrete groups (including free groups), and perhaps all
discrete groups, are exact. (See Remark 7.8 of \cite{Kr0.5}.)
Separable exact \ca s can be characterized as exactly those \ca s
which occur as subquotients of the (nuclear) CAR (or
$2^{\infty}$ UHF) algebra (\cite{Kr1}).

In this paper, we show that every separable exact \ca\  is isomorphic
to a subalgebra of the Cuntz algebra $\OA{2}.$ Thus, a separable
\ca\  is exact \ifo\  it is isomorphic to a subalgebra of
of a nuclear \ca, \ifo\  it is isomorphic to a subalgebra of
the particular nuclear \ca\  $\OA{2}.$

The methods used to prove the embedding in $\OA{2}$ show that
separable exact \ca s which are ``close'' in a certain sense
have nearby embeddings in $\OA{2}.$ We prove that if $X$ is
a compact metric space which is sufficiently nice (certainly
including all compact manifolds and all finite CW complexes),
and if $A$ is a continuous field over $X$ in the sense of Dixmier
(Chapter 10 of \cite{Dx}) such that the algebra of continuous
sections of $A$ is separable and exact, than $A$ has a continuous
representation in $\OA{2}.$ That is, there is a collection of
injective \hm s $\ph_x$ from the fibers $A (x)$ of $A$ to $\OA{2}$
such that, for every continuous section $a$ of $A,$ the function
$x \mapsto \ph_x (a (x))$ is continuous. Moreover, one can say something
about the ``smoothness'' of these functions. For example, we show that
there
are injective \hm s $\ph_{\te}$ from the rotation algebras $A_{\te}$
(rational and irrational) to $\OA{2}$ such that, if $u (\te)$ and
$v (\te)$ denote the standard generators of $A_{\te},$ then
there is a constant $C$ such that
\[
\max ( \| \ph_{\te_1} (u (\te_1)) - \ph_{\te_2} (u (\te_2)) \|, \,
       \| \ph_{\te_1} (v (\te_1)) - \ph_{\te_2} (v (\te_2)) \| )
  < C | \te_1 - \te_2 |^{1/2}
\]
for all $\te_1,$ $\te_2 \in \R.$ Haagerup and R\o rdam \cite{HR}
obtained representations on a Hilbert space with these properties,
and showed that even there the exponent $\half$ can't be improved.
Our work provides an independent proof of the Haagerup and R\o rdam
representation--a proof not using unbounded operators.

Blanchard \cite{Bl} has a somewhat different approach to the
representation of continuous fields in $\OA{2}.$ He is able to
allow more general base spaces $X,$ but obtains no information
on smoothness. In particular, his methods do not provide the
$\Lip^{1/2}$ representation (or a $\Lip^{\af}$ representation for
any $\af$) of the field of rotation algebras in $\OA{2}.$

In \cite{El1}, George Elliott initiated a program for the classification
of separable nuclear simple \ca s of real rank zero in terms of
K-theoretic invariants. The purely infinite
case was first tackled in \cite{Rr1} (building on the work of
\cite{BKRS}), and has since been investigated in a number of papers.
(See \cite{BEEK}, \cite{ElR}, \cite{Ln2}, \cite{Ln3}, \cite{LP0},
\cite{LP}, \cite{Rr1}, \cite{Rr2}, \cite{Rr3}, and
\cite{Rr4}.) The classification program predicts that if $A$ is
\snus, then $\OT{A}$ should be isomorphic to $\OA{2},$ and that if, in
addition, $A$ is purely infinite, then $\OIA{A}$
should be isomorphic to $A.$ (The K\"{u}nneth formula for
\ca s \cite{Sch} implies that $K_* (\OT{A}) \cong K_* (\OA{2})$
and $K_* (\OIA{A}) \cong K_* (A).$) We use the result
on embeddings in $\OA{2},$ together with some of the lemmas in its
proof and some additional results, to prove that these isomorphisms
do in fact hold: $\OT{A} \cong \OA{2}$ for separable nuclear unital
simple $A,$ and $\OIA{A} \cong A$ for \kfive\  $A.$

These results are the starting point for a proof by the second
author of a general classification theorem for \kfalg s satisfying
the universal coefficient theorem \cite{Ph2}. The first author
also has an independent proof \cite{Kr3} of this classification theorem,
which does not use the isomorphisms of tensor products above
directly, but rather obtains them as a consequence of the general
classification result.

This paper is organized as follows. The rest of the introduction
contains some general terminology and notation, and some more or less
well known results that we use often enough that it is convenient
to have them restated here. In Section 1 we prove several technical
results on approximation and perturbation of \ucp\  maps.
In the second section, we use these results to prove that separable
exact \ca s can be embedded in $\OA{2}.$ Section 3 contains the
proofs of $\OT{A} \cong \OA{2}$ and $\OT{A} \cong \OA{2}.$ Those
interested only in the classification program need read no further.
In the fourth section, we present some preliminaries on continuous
fields, including results based on earlier papers and results based
on Section 1 of this paper. In the fifth section we give a version
of the argument Haagerup and R\o rdam use in \cite{HR} to obtain
continuous representations of continuous fields when one merely knows
that any two nearby fibers have nearby embeddings. Section 6 is devoted
to the detailed examination of the field of rotation algebras.

The first author would like to thank Mikael R\o rdam and George Elliott
for valuable discussions.
The second author would like to thank Ed Effros and Gilles Pisier for
valuable discussions, Uffe Haagerup for a stimulating question, and
K\o benhavns Universitet for its hospitality during the fall semester
of 1995, when some of this paper was written.

Here is some general notation we use throughout. $M_n$ is the
algebra of $n \times n$ matrices, with matrix units $\{e_{ij}\}.$
If $H$ is a separable infinite dimensional Hilbert space, then
$K = K (H)$ denotes the algebra of compact operators on $H,$ and $L (H)$
denotes the algebra of bounded operators on $H.$ The unitization of a
\ca\  $A$ is denoted $A\unit;$ this means we add a new unit
even if $A$ already has one. We let $\tilde{A}$ denote $A$ if $A$ is
unital and $A\unit$ if not. The unitary group of a unital \ca\  $A$
is denoted $U (A),$ and $U_0 (A)$ is the connected component of
$U (A)$ containing the identity. When we refer to a unital subalgebra
$B$ of a \ca\  $A,$ we implicitly mean that the $B$ is supposed to
contain the identity of $A.$

We present here some (well) known results which are used frequently
enough that is is convenient to restate them.

The first part of the inequality in the following estimate is the best
possible, as can be seen by taking $p = 1$ and $x$ to be a positive
real scalar less than $1.$ Note that $1 - ( 1 - \dt)^{1/2}$ is
approximately $\dt / 2$ for small $\dt.$ The second part of the
inequality gives the best linear estimate over the relevant range, as
can be seen by letting $\| x^* x - p\| \to 1$ (for example, letting
$x \to 0$).

\begin{lem} 
Let $A$ be a \ca, let $p \in A$ be a \pj, and let $x \in A$
satisfy $xp = x$ and $\| x^* x - p\| < 1.$ Then the formula
$v = x (x^* x)^{-1/2}$ (functional calculus evaluated in $p A p$)
defines a partial isometry in $A$ with $v^* v = p.$ Moreover,
\[
\|v - x\| \leq 1 - ( 1 - \| x^* x - p\|)^{1/2} \leq \| x^* x - p\|.
\]
\end{lem}

\begin{pff}
It is well known that if $\| x^* x - p\| < 1$ then $v$ is a
partial isometry satisfying $v^* v = p.$ For the rest, let
$\dt = \| x^* x - p\|.$ A calculation, using the fact that
$(x^* x)^{1/2} \in p A p,$ shows that
\[
(v - x)^* (v - x) = [(x^* x)^{1/2} - p]^2
\]
and
\[
\|v - x\| = \|(x^* x)^{1/2} - p\| \leq
           \sup \{ |t^{1/2} - 1| : |t - 1| \leq \dt \} =
           1 - ( 1 - \dt)^{1/2}.
\]
This gives the first inequality. For the second, note that
$(1 - \dt)^{1/2} \geq 1 - \dt,$ whence $(1 - \dt)^{1/2} - 1 \geq - \dt.$
\end{pff}

\begin{thm} 
(Theorem 3.1 of \cite{CE1}; also see \cite{Kr0}.)
A separable unital \ca\  $A$ is nuclear if and only if there are
sequences of \ucp\  maps $S_k : A \to M_{n(k)}$ and
$S_k : M_{n(k)} \to A,$ for suitable $n (k),$ such that
$\limi{k} T_k \circ S_k (a) = a$ for all $a \in A.$
\end{thm}

The maps in \cite{CE1} are not required to be unital, but this is easily
fixed. See the note on this point in the proof of Proposition 4.3 of
\cite{EH}.

\begin{thm} 
(Choi-Effros Lifting Theorem; see the corollary to Theorem 7 of
\cite{Ar}.)
Let $A$ be a separable nuclear unital \ca, let $B$ be a
unital \ca, and let $J$ be an ideal of $B,$ with quotient
map $\pi : B \to B/J.$ Then for every
\ucp\  map $S : A \to B / J,$ there is a \ucp\  map
$T : A \to B$ which lifts $S,$ that is, such that $\pi \circ T = S.$
\end{thm}

\begin{prp} 
Let $A$ and $B$ be unital \ca s. Let $E \subset A$ be
an operator system (in the sense of Choi and Effros; see \cite{CE2}),
and let $S : E \to B$ be a nuclear \ucp\  map. Then for every
\fd\  operator system
$E_0 \subset E$ and every $\ep > 0,$ there is a \ucp\  map
$T : A \to B$ such that $\| T |_{E_0} - S |_{E_0} \| < \ep.$
\end{prp}

\begin{pff}
This is similar to Proposition 4.3 of \cite{EH}.
Since $S$ is nuclear, there are $n$ and  \ucp\  maps
$P_0 : E \to M_n$ and $Q : M_n \to B$ such that
$\| Q \circ P_0 |_{E_0} - S |_{E_0} \| < \ep.$
The Arveson extension theorem (Theorem 6.5 of \cite{Pl})
provides a \ucp\  map $P : A \to M_n$ such that
$P |_{E_0}  =  P_0 |_{E_0}.$ Set $T = Q \circ P.$
\end{pff}

\begin{dfn} 
Let $A$ and $B$ be \ca s, with  $A$ separable and
$B$ unital. Two \hm s $\ph, \, \ps : A \to B$ are {\em \ayue}\  if
there is a sequence $(u_n)$ of unitaries in $B$ such that
$\limi{n} \| u_n \ph (a) u_n^* - \ps (a) \| = 0$ for all $a \in A.$
\end{dfn}

We will frequently use the following special case of Elliott's
approximate intertwining argument, of which the original form is
Theorem 2.1 of \cite{El1}.

\begin{lem}      
Let $A$ and $B$ be separable unital \ca s, and let
$\ph: A \to B$ and $\ps : B \to A$ be \hm s such that
$\ps \circ \ph$ is \ayue\  to $\id_A$ and $\ph \circ \ps$ is \ayue\  to
$\id_B.$ Then $A \cong B.$
\end{lem}

\begin{pff}
This is contained in the proof of Theorem 6.2 (1) of \cite{Rr2}.
(Or see Proposition A of \cite{Rr3}.)
\end{pff}

\begin{prp}      
Let $D$ be a unital \pisca. Then any two unital \hm s from
$\OA{2}$ to $D$ are \ayue.
\end{prp}

\begin{pff}
This is a special case of Theorem 3.6 of \cite{Rr1}. The required
two conditions on $D$ (that $U (D) / U_0 (D) \to K_1 (D)$ be an
isomorphism and that $D$ have finite exponential length in the sense
of \cite{Rn}) follow from Theorem 1.9 of \cite{Cu1} and
from \cite{Ph} or \cite{Ln} respectively.
\end{pff}

\begin{thm}      
(\cite{Rr3}) $\OT{\OA{2}} \cong \OA{2}.$
\end{thm}

\section{Completely positive maps and perturbation}

In this section we prove several approximation and perturbation
results which are frequently required  later in the paper. At the
end of the section, we combine these results (without using their
full strength) to prove a lemma on \aue\  of \hm s to tensor products
with $\OA{2}.$ Combined with the isomorphism
$\OT{\OA{2}} \cong \OA{2},$ it implies that
two unital injective \hm s from a
separable exact \ca\  to $\OA{2}$ are \ayue.

The first of the three
main results of this section is Proposition 1.7,
which shows that a nuclear \ucp\  map from a unital \pisca\  to
itself can be approximated on finite sets by maps of the form
$a \mapsto s^* a s$ for isometries $s.$ The proof is somewhat different
from that of a similar result in \cite{Kr3}, relying much more
heavily on properties of \pisca s. The next important result is
Lemma 1.10,
in which we show, under suitable exactness
and nuclearity assumptions, that if $S,$ $T : A \to B$ are
\uacp\  and $S$ is sufficiently close to being completely
isometric on a \fd\  operator system $E,$ then $T \circ S^{-1}$
can be approximated on $S (E)$ by a \ucp\  map. The last main result
is Lemma 1.12,
in which approximate ``similarity'' via isometries
in $A$ is shown to imply \aue\  in $\OT{A}.$
We use it to show that two injective unital \hm s from
a separable unital exact \ca\  to $\OA{2}$ are \ayue.

We begin with some preliminary lemmas on \pisca s.

\begin{lem}  
Let $A$ be a \ca, let $a,$ $h \in A$ be selfadjoint elements with
$0 \leq a \leq 1$ and $0 \leq h \leq 1,$ and let $q \in A$ be
a projection. Then $\| q a - q \| \leq 12 \| q h a h - q \|^{1/3}.$
\end{lem}

\begin{pff}
Represent $A$ faithfully on a Hilbert space $H.$
Note that $\| q a - q \| = \| a q - q \|$ and
$\| q h a h - q \| = \| h a h q - q \|.$
It suffices to show that if $\xi \in q H,$ then
$\| a \xi - \xi \| \leq 12 \| h a h \xi - \xi\|^{1/3}.$

We first claim that if $b \in L (H)$ satisfies $0 \leq b \leq 1$ and
$\et \in H$ satisfies $\| \et \| = 1,$ then
$\|b \et - \et \| \leq 4 ( 1 - \| b \et \|)^{1/3}.$ To see this, set
$\dt = 1 - \|b \et \|,$ let  $\rh = \dt^{1/3},$ and let
$p \in L (H)$ be the spectral projection
for $b$ corresponding to $[ 1 - \rh, 1].$ Then
\beqr
(1 - \dt)^2 & = & \|b \et \|^2 = \|b p \et \|^2 + \|b (1 - p) \et \|^2
          \leq \| p \et \|^2 + (1 - \rh)^2 \| (1 - p) \et \|^2    \\
 & = & 1 - \| (1 - p) \et \|^2 + (1 - \rh)^2 \| (1 - p) \et \|^2
           = 1 - \rh (2 - \rh) \| (1 - p) \et \|^2.
\eeqr
It follows that
\[
\| (1 - p) \et \| \leq
 \sqrt{ \left( \frac{2 - \dt}{2 - \rh} \right)
                              \left( \frac{\dt}{\rh} \right)}
  \leq \sqrt{2} \dt^{1/3}.
\]
So
\[
\|b \et - \et \| \leq \| b \| \| (1 - p) \et \| + \| b p \et - p \et \|
                         + \| (1 - p) \et \|
     \leq 2 \sqrt{2} \dt^{1/3} + \rh < 4 \dt^{1/3}.
\]
This proves the claim.

Now let $\xi \in q H$ satisfy $\| \xi \| = 1.$ Since
$\|a\|,$ $\|h\| \leq 1,$ we have
$\|h \xi \| \geq 1 - \| h a h \xi - \xi\|.$
Applying the claim to $h$ and $\xi,$ we get
$\| h \xi - \xi\| \leq 4 \| h a h \xi - \xi\|^{1/3}.$
Also, with $\et = \frac{1}{\|h \xi \|} h \xi,$ we have
\[
\| a \et \| \geq \frac{1}{\|h \xi \|} \| h a h \xi \| \geq
        \frac{1}{\|h \xi \|} (1 - \| h a h \xi - \xi\|)
         \geq 1 - \| h a h \xi - \xi\|,
\]
so
\[
\| a h \xi - h \xi \| = \|h \xi\| \| a \et - \et \|
     \leq \|h \xi\| \cdot 4 \| h a h \xi - \xi\|^{1/3}
    \leq 4 \| h a h \xi - \xi\|^{1/3}.
\]
Therefore
\[
\| a \xi - \xi \| \leq \| a \| \| \xi - h \xi \| +
      \| a h \xi - h \xi \| + \| h \xi - \xi\|
  \leq 12 \| h a h \xi - \xi\|^{1/3}.
\]
\end{pff}

\begin{lem}  
Let $A$ be a \pisca, and let $a_1, \dots, a_n \in A$ be positive
elements with $\| a_j \| = 1$ for all $j.$
Then for every $\ep > 0$ there are nonzero
\mops\  $p_1, \dots, p_n \in A$
such that $\| p_j a_j - p_j \| < \ep$ for all $j.$
\end{lem}

\begin{pff}
Choose $\dt > 0$ with $12 (2 \dt)^{1/3} < \ep.$
Choose an irreducible representation $\pi$ of $A$ on a Hilbert
space $H.$ By induction, we construct a sequence
$\xi_1, \dots, \xi_n$  of orthogonal unit vectors in $H$ with
$\|\pi (a_j) \xi_j - \xi_j \| < \dt$ for $1 \leq j \leq n.$
Choose $\xi_1$ to be any unit vector in the spectral subspace
for $\pi (a_1)$ corresponding to $[1 - \dt, 1].$  Given
$\xi_1, \dots, \xi_j,$ let $H_0$ be the spectral subspace
for $\pi (a_{j + 1})$ corresponding to $[1 - \dt, 1].$ Since
$\pi (A)$ contains no compact operators, this space must be
infinite dimensional. Therefore it must nontrivially intersect the
finite codimension subspace $\spn (\xi_1, \dots, \xi_j)^{\perp},$
and we take $\xi_{j + 1}$ to be any unit vector in the intersection.

Let $p_j \in L (H)$ be the projection onto $\C \xi_j,$ and
let $p = p_1 + \cdots + p_n$ be the projection onto
$\spn (\xi_1, \dots, \xi_n).$ Let
\[
L = \{ a \in A : \pi (a) p = 0\} \andeqn
N = \{ a \in A : \pi (a) p = p \pi (a) \}.
\]
Then $L$ is a left ideal of $A,$ $N$ is a C*-subalgebra of $A,$ and
$L \cap L^*$ is an ideal in $N.$ Define a \ucp\  map
$T : A \to L (p H) \cong M_n$ by $T (a) = p \pi (a) p.$
Then $T |_N$ is a \hm\  with kernel $L \cap L^*.$ The Kadison
Transitivity Theorem implies it is surjective. Indeed, let
$u \in L (p H)$ be unitary.
Since $\pi$ is injective, Theorem 5.4.5 of \cite{KR} provides a unitary
$v \in \tilde{A}$ such that
$p \pi (v) p = u.$ Since $\pi (v)$ and $p \pi (v) p$ are both unitary,
$p$ must commute with $\pi (v).$ (This is well known, but see a closely
related result in Lemma 1.11 below.) So $v \in \tilde{N}$ and
$T (v) = u.$ This shows that the image of $\tilde{N}$ contains all
unitaries in $L (p H),$ and so is
all of $L (p H).$ The image of $N$ is an ideal of
codimension at most $1.$ We
may clearly assume $n \geq 2;$ then $L (p H)$ has no proper ideals
of codimension at most $1,$ so $T |_N$ must be surjective.

By Theorem 4.6 of \cite{Lr} (essentially Proposition 2.6 of \cite{AP1}),
there are $h_1, \dots, h_n \in N$ satisfying $T (h_j) = p_j,$
$0 \leq h_j \leq 1,$ and $h_j h_k = 0$ for $j \neq k.$
Then
\[
\| h_j a_j h_j \| \geq \| p \pi (h_j) \pi (a_j) \pi (h_j) p \|
   = \| p_j \pi (a_j) p_j \| = \langle \pi (a_j) \xi_j, \xi_j \rangle
   \geq 1 - \dt.
\]
The proof of Lemma 1.7 of \cite{Cu1} provides a projection
$q_j \in \overline{h_j A h_j}$ such that
$\| q_j h_j a_j h_j - q_j \| < 2 \dt.$ Since $h_j h_k = 0$ for
$j \neq k,$ we also have $q_j q_k = 0$ for $j \neq k.$
Moreover, $\| q_j a_j - q_j \| \leq 12 (2 \dt)^{1/3} < \ep$ by the
previous lemma.
\end{pff}

\begin{lem}  
Let $A$ be a unital \pisca, and let $F \subset A$ be a finite subset
consisting of positive elements $a$ satisfying
$\af \leq a \leq \bt$ for fixed positive real numbers $\af$ and $\bt.$
Then for all $\rh > 0$ there is $\dt > 0$ such that the following holds:
If there are positive elements $c_1, \dots, c_n \in A$ with
$\|c_j\| = 1$ such that $\|c_j a c_j - \ld_j (a) c_j^2 \| < \dt$ for
$1 \leq j \leq n,$ $a \in F,$ and some numbers
$\ld_j (a) \in [\af, \bt],$ then there are nonzero
\mops\  $p_1, \dots, p_n \in A$ such that
$\|p_j a p_k - \dt_{jk} \ld_j (a) p_j \|  < \rh$
for $1 \leq j \leq n$ and $a \in F,$
where $\dt_{jk}$ is the Kronecker delta.
\end{lem}

\begin{pff}
\Wolog\  $\bt = 1,$ $1 \in F,$ and $\ld_j (1) = 1$ for all $j.$
(If $\bt \neq 1,$ we can rescale by multiplying by $1 / \bt.$
If $1$ is already in $F$ but $\ld_j (1) \neq 1,$
it does no harm to redefine $\ld_j (1).$)
The proof is now by induction on the number of elements of $F,$
but to make the argument work we require the additional conclusion
$\|p_j c_j - p_j\| < \rh$ for $1 \leq j \leq n.$ If
$F$ has only one element, then it is $1,$ and the existence of the
required projections $p_j$ is just Lemma 1.2.

Assume now that $F$ has more than one element, and that the lemma
is known to hold for all smaller such sets.
Define $\mu = \af \rh / 15.$ (Then in particular
$0 < \mu < \rh/5.$) Choose $b \in F$ with
$b \neq 1,$ and let $F_0 = F \setminus \{b\}.$
Use the induction hypothesis to choose $\dt$ with $0 < \dt < \mu$
such that the
conclusion of the lemma holds for $F_0$ in place of $F$ and
$\mu$ in place of $\rh.$ Let $q_1, \dots, q_n$ be the nonzero
\mops\ provided by the conclusion.

Set $x_j = \ld_j (b)^{-1/2} q_j b^{1/2}.$ Then
\beqr
\| x_j x_j^* - q_j \| & \leq &
          \frac{1}{\af} \|q_j b q_j - \ld_j (b) q_j \|   \\
 & \leq & \frac{1}{\af} \left( 4 \|q_j c_j - q_j\| +
        \| q_j \| \|c_j b c_j - \ld_j (b) c_j^2 \| \| q_j \| \right)
  <  \frac{1}{\af} (4 \mu + \dt) \leq \frac{5 \mu}{\af} =
      \frac{\rh}{3}.
\eeqr
By Lemma 0.1
there exist partial isometries $v_j$ such that
$v_j v_j^* = q_j$ and $\|v_j - x_j \| < \rh /3.$

Lemma 1.2
provides nonzero \mops\  $e_1, \dots, e_n$ such that
$\|e_j v_j^* v_j - e_j \| < \rh /24$ for $1 \leq j \leq n.$ It
follows that $\|v_j^* v_j e_j v_j^* v_j - e_j \| < \rh /12,$ and
applying functional calculus to $v_j^* v_j e_j v_j^* v_j$ yields
a \pj\ $f_j \leq v_j^* v_j$ such that $\|f_j - e_j \| < \rh /6.$
Since the $e_j$ are mutually orthogonal, it follows that
$\|f_j f_k\| < \rh /3$ for $j \neq k.$
Define $p_j = v_j f_j v_j^* \leq q_j.$

Since $p_j \leq q_j,$ the $p_j$ are mutually orthogonal, and the
corresponding properties of the $q_j$ imply that
$\|p_j a p_k - \dt_{jk} \ld_j (a) p_j \| < \mu \leq \rh$
for $1 \leq j, \, k \leq n$ and $a \in F_0$ and that
$\|p_j c_j - p_j\| < \mu \leq \rh$ for $1 \leq j \leq n.$
Furthermore, in the estimate of $\| x_j x_j^* - q_j \|$ above we saw
that $\|q_j b q_j - \ld_j (b) q_j \| < 5 \mu,$ which is at most
$\rh.$ The same estimate therefore holds with $p_j$ in place of $q_j.$

It remains to estimate $\|p_j b  p_k\|$ for $j \neq k.$
Using $\ld_j (b)^{1/2}, \, \ld_k (b)^{1/2} \leq 1$ in the last step,
we have, for $j \neq k,$
\beqr
\|p_j b  p_k\| & = & \|p_j q_j b q_k p_k\| =
        \ld_j (b)^{1/2} \ld_k (b)^{1/2} \|p_j x_j x_k^* p_k\|
   \leq  \ld_j (b)^{-1/2} \ld_k (b)^{-1/2}
                          \| f_j v_j^* x_j x_k^* v_k f_k \|          \\
  & \leq & \ld_j (b)^{1/2} \ld_k (b)^{1/2}
                     \left( \| f_j v_j^* v_j v_k^* v_k f_k \| +
              \| x_k \| \|x_j - v_j \| + \|x_k - v_k \| \right)  \\
  & < & \ld_j (b)^{1/2} \ld_k (b)^{1/2}
                \left( \|f_j f_k\| + \ld_k (b)^{-1/2} \rh / 3 +
                            \rh /3 \right) \leq \rh.
\eeqr
\end{pff}

The main technical parts of the proof of the next
lemma are the excision of pure states
(see \cite{AAP}) and the previous lemma.
(Note that Proposition 2.3 of \cite{AAP} and Lemma 1.4
imply that every state on a unital \pisca\  is a
weak* limit of pure states.)

\begin{lem}  
Let $A$ be a unital \pisca, and let $\om$ be a state on $A.$ Then for
every $\ep > 0$ and every finite subset $F \subset A$ there exists
a nonzero \pj\  $p \in A$ such that
$\| p a p - \om (a) p \| < \ep$ for all $a \in F.$
\end{lem}

\begin{pff}
Without loss of generality, we may assume $1 \in F$ and that
$F$ consists of positive elements of norm at most 1. It further
suffices to prove the lemma using the set
$\{\frac{1}{2} (1 + a) : a \in F\}$ instead of $F;$
thus we may assume \wolog\  that $a \geq {\frac{1}{2}}$
for all $a \in F.$ In particular, if $\mu$ is any state on $A,$ then
$\mu (a) \geq {\frac{1}{2}}$ for all $a \in F.$

Since the set of all states is the weak* closed convex hull of the
set of pure states, there are $\af_1, \dots, \af_n \in [0,1]$ and
pure states $\om_1, \dots, \om_n$ of $A$ such that
$\sum_{j = 1}^n \af_j = 1$ and
$| \om(a) - \sum_{j = 1}^n \af_j \om_j (a) | < \ep/2$
for all $a \in F.$ Choose $\dt> 0$ as in Lemma 1.3
for $F$ and the number $\rh = \ep / (2 n^2).$ Excision of pure states
(see Proposition 2.2 of \cite{AAP})
implies that there are positive elements $c_1, \dots, c_n \in A$
of norm $1$ such that
$\|c_j a c_j - \om_j (a) c_j^2 \| < \dt$ for $1 \leq j \leq n$
and $a \in F.$ Apply Lemma 1.3
with $\ld_j = \om_j |_F$ to obtain
nonzero \mops\  $q_1, \dots, q_n \in A$
such that $\|q_j a q_k - \dt_{jk} \om_j (a) q_j \| < \ep / (2 n^2).$

Since $a$ is purely infinite and simple, there exist isometries
$s_1, \dots, s_n \in A$ such that $p_j = s_j s_j^* \leq q_j.$
Define $s = \sum_{j = 1}^n \af_j^{1/2} s_j.$ Since the $p_j$ are
orthogonal and $\sum_{j = 1}^n \af_j = 1,$ one immediately checks that
$s$ is an isometry. Define $p = s s^*.$ Then, for $a \in F$ we have
\beqr
\| p a p - \om(a) p\| & = & \| s^* a s - \om(a) \cdot 1 \|
    <  \frac{\ep}{2} +
    \sum_{j, \, k = 1}^n \af_j^{1/2} \af_k^{1/2}
                      \| s_j^* a s_k - \dt_{jk} \om_j (a) \cdot 1 \|  \\
 & = & \frac{\ep}{2} +
    \sum_{j, \, k = 1}^n \af_j^{1/2} \af_k^{1/2}
                      \| p_j a p_k - \dt_{jk} \om_j (a) p_j \|
    <  \frac{\ep}{2} + n^2 \frac{\ep}{2 n^2} = \ep,
\eeqr
since $\af_j \leq 1$ and $p_j \leq q_j.$
\end{pff}

\begin{lem}  
Let $A$ be a unital \pisca, let $T : A \to M_n$ be a \ucp\  map,
and let $\ph : M_n \to A$ be a (not necessarily unital) \hm. Then for
every $\ep > 0$ and every finite subset $F \subset A$ there exists
a partial isometry $s \in A$ such that $s^* s = \ph (1)$ and
$\| s^* a s - \ph (T (a)) \| < \ep$ for all $a \in F.$
\end{lem}

\begin{pff}
\Wolog\  we may assume $1 \in F$ and that all elements of $F$ have norm
at most $1.$ Choose $\rh > 0$ such that $\rh \leq \min (1, \ep / 4),$
and also so small that if $q$ and $q'$ are \pj s such that
$\|q q'\| < 4 \rh$ then
there are orthogonal \pj s $r$ and $r'$ unitarily equivalent
to $q$ and $q'$ respectively. Define $\dt = \rh/ (3 n^3);$
then $0 < \dt \leq \ep/n^3.$

Let $\{ \xi_1, \dots, \xi_n \}$ be the standard orthonormal basis of
$\C^n,$ and let $e_{kl},$ for $1 \leq k, \, l \leq n,$ be the
standard matrix units satisfying $e_{kl} \xi_j = \dt_{jl} \xi_k.$
Following Theorem 5.1 of \cite{Pl} and the preceding discussion,
define a state on $M_n \otimes A$ by
\[
\om \left(\sum_{k, \, l = 1}^n e_{kl} \otimes a_{kl} \right) =
   \frac{1}{n} \sum_{k, \, l = 1}^n \langle T (a_{kl}) \xi_l, \xi_k
                                                            \rangle.
\]
(Note that $\om (1) = 1$ because $T (1) = 1.$) As there, we then have
\[
T (a) = n \sum_{k, \, l = 1}^n \om (e_{kl} \otimes a) e_{kl}
\]
for $a \in A.$
By Lemma 1.4,
there is a nonzero \pj\  $p_0 \in M_n \otimes A$
such that
$\| p_0 (e_{ij} \otimes a) p_0 - \om (e_{ij} \otimes a) p_0 \| < \dt$
for all $a \in F$ and $1 \leq i, \, j \leq n.$

Since
$M_n \otimes A$ is purely infinite and simple, there is a nonzero
\pj\  $p \in M_n \otimes A,$ with $p \leq p_0,$
such that there are partial isometries
$s_1, \dots, s_n \in M_n \otimes A$ satisfying
$s_j s_j^* = p$ and $s_j^* s_j = e_{11} \otimes \ph (e_{jj}).$
We then have
\[
\| s_i^* (e_{kl} \otimes a) s_j -
           \om (e_{kl} \otimes a) (e_{11} \otimes \ph (e_{ij})) \| < \dt
\]
for all $a \in F$ and $1 \leq i, \, j \, k, \, l \leq n.$

Define $c = \sum_{k = 1}^n (e_{1k} \otimes 1) s_k \in M_n \otimes A.$
One checks that for $a \in F$ we have
\[
c^* (e_{11} \otimes a) c =
         \sum_{k, \, l = 1}^n s_k^* (e_{kl} \otimes a) s_l,
\]
from which it follows that
\beqr
\lefteqn{\| n c^* (e_{11} \otimes a) c - e_{11} \otimes \ph ( T (a)) \|}
                                     \\
& \leq & n \sum_{k, \, l = 1}^n
   \| s_k^* (e_{kl} \otimes a) s_l
               - \om (e_{kl} \otimes a) (e_{11} \otimes \ph (e_{kl})) \|
         < n^3 \dt.
\eeqr
Putting in particular $a = 1,$ we obtain
$\| n c^* (e_{11} \otimes 1) c - e_{11} \otimes \ph (1) \|< n^3 \dt.$
Set
\[
d = \sqrt{n} (e_{11} \otimes 1) c (e_{11} \otimes \ph (1))
             \in (e_{11} \otimes 1) (M_n \otimes A) (e_{11} \otimes 1).
\]
Then $\|d^* d - e_{11} \otimes \ph (1)\|  < n^3 \dt,$ so Lemma 0.1
implies that $d (d^* d)^{-1/2}$ (with functional calculus taken in
$(e_{11} \otimes \ph (1)) (M_n \otimes A) (e_{11} \otimes \ph (1))$)
is a partial isometry in
$(e_{11} \otimes 1) (M_n \otimes A) (e_{11} \otimes 1)$ satisfying
\[
[d (d^* d)^{-1/2}]^* [d (d^* d)^{-1/2}] = e_{11} \otimes \ph (1)
\andeqn
\|d (d^* d)^{-1/2} - d \| < n^3 \dt \leq \rh.
\]

Let $s \in A$ be the partial
isometry such that $d (d^* d)^{-1/2} = e_{11} \otimes s.$
We check that  $s$ satisfies the required estimates. First note that
$\| n c^* (e_{11} \otimes a) c - e_{11} \otimes \ph ( T (a)) \|
                                                          < n^3 \dt$
implies that
\[
\| d^* (e_{11} \otimes a) d - e_{11} \otimes \ph ( T (a)) \| < n^3 \dt.
\]
Moreover, $\| e_{11} \otimes s - d \| < \rh$ implies that
$\|d\| < 1 +\rh < 2.$ Therefore
\[
\| s^* a s - \ph ( T (a)) \| =
  \| (e_{11} \otimes s)^* (e_{11} \otimes a) (e_{11} \otimes s)
                             - e_{11} \otimes \ph ( T (a)) \|
        < 2 \rh + \rh + n^3 \dt \leq \ep,
\]
since $\| a \| \leq 1$ for all $a \in F.$
\end{pff}

The following lemma and its proof were inspired by Kasparov's
Stinespring theorem for Hilbert modules (\cite{Ks1}).

\begin{lem}     
Let $A$ be a unital \ca, and let $T : M_n \to A$ be \uacp.
Denote by $\{e_{ij}\}$ the standard system of matrix units in $M_n.$
Then there exists a partial isometry $t \in \Mt{n} \Mt{n} A$
such that
\[
t^* t = e_{11} \otimes e_{11} \otimes 1 \andeqn
  t^* (b \otimes 1 \otimes 1) t = e_{11} \otimes e_{11} \otimes T (b)
\]
for $b \in M_n.$
\end{lem}

\begin{pff}
Set $x = \sum_{i, \, j = 1}^n e_{ij} \otimes e_{ij} \in \Mt{n} M_n.$
Note that $n^{-1} x$ is a \pj, so $x$ is positive. Therefore so is
\[
y = (\id_{M_n} \otimes T) (x) =
       \sum_{i, \, j = 1}^n e_{ij} \otimes T (e_{ij}) \in \Mt{n} A.
\]
Write $y^{1/2} = \sum_{i, \, j = 1}^n e_{ij} \otimes a_{ij}$ for
suitable $a_{ij} \in \Mt{n} A.$ Since $y^{1/2}$ is selfadjoint and has
square $y,$ we obtain
$a_{ik}^* = a_{ki}$ and  $\sum_{j = 1}^n a_{ij} a_{jk} = T (e_{ik})$
for $1 \leq i, \, k \leq n.$

Define $t = \sum_{i, \, j = 1}^n e_{i1} \otimes e_{j1} \otimes a_{ji}.$
A computation shows that
\[
t^* (e_{ik} \otimes 1 \otimes 1) t =
          e_{11} \otimes e_{11} \otimes \sum_{j = 1}^n a_{ij} a_{jk}
   = e_{11} \otimes e_{11} \otimes  T (e_{ik}).
\]
It follows that
$t^* (b \otimes 1 \otimes 1) t = e_{11} \otimes e_{11} \otimes T (b)$
for all $b \in M_n.$ In particular,
$t^* t = e_{11} \otimes e_{11} \otimes T (1)
                          = e_{11} \otimes e_{11} \otimes 1.$
The last equation implies that $t$ is a partial isometry.
\end{pff}

\begin{prp}     
Let $A$ be a unital \pisca, and let $V : A \to A$ be a nuclear
\ucp\  map.
Then for
every $\ep > 0$ and every finite subset $F \subset A$ there exists
a nonunitary isometry $s \in A$ such that
$\| s^* a s - V (a) \| < \ep$ for all $a \in F.$
\end{prp}

\begin{pff}
By the definition of nuclearity, $V$ is a pointwise norm limit
of maps of the form $T \circ S,$ with $n \in \N$ and $S : A \to M_n,$
$T : M_n \to A$ \uacp. Therefore it suffices to prove the
proposition for maps of the form $T \circ S.$

Let $\{e_{ij}\}$ be the standard system of matrix units in $M_n.$
Since $A$ is purely infinite and simple, so is $\Mt{n} \Mt{n} A.$
Therefore there is $t_1 \in \Mt{n} \Mt{n} A$ such that
\[
t_1^* t_1 = 1 \andeqn t_1 t_1^* < e_{11} \otimes e_{11} \otimes 1.
\]
Apply the previous lemma
to $T,$ and call the resulting partial isometry $t_2.$
Define a (nonunital) \hm\  $\ph_0 : M_n \to A$ by
identifying $A$ with the corner
\[
A_0 = (e_{11} \otimes e_{11} \otimes 1) (\Mt{n} \Mt{n} A)
             (e_{11} \otimes e_{11} \otimes 1),
\]
and setting $\ph_0 (b) = t_1 (b \otimes 1 \otimes 1) t_1^*.$ Set
$t = t_1 t_2,$ which is an isometry in $A_0,$ and regard
it as an element of $A.$ Then, using the result of
the previous lemma,
we have $T (b) = t^* \ph_0 (b) t$ for all $b \in M_n.$

Let $p = t_1 t_1^*,$ which we also regard as an element of $A.$
Since $A$ is
purely infinite and simple and $1 - p \neq 0,$ there is a nonzero
\hm\  $\ph_1 : M_n \to (1 - p) A (1 - p).$ Define
$\ph : M_n \to A$ by $\ph (b) = \ph_0 (b) + \ph_1 (b).$ Then we still
have $T (b) = t^* \ph (b) t$ for all $b \in M_n.$

Use Lemma 1.5
to choose a partial isometry $s_0 \in A$ such that
$s_0^* s_0 = \ph (1)$ and
$\| s_0^* a s_0 - \ph (S (a)) \| < \ep$ for all $a \in F.$
Set $s = s_0 t.$ Then for $a \in F$ we have
\[
\| s^* a s - T (S (a)) \| = \| t^* s_0^* a s_0 t - t^* \ph (S (a)) t \|
  < \ep.
\]
Moreover, $s$ is an isometry because
$s_0^* s_0 \geq t_1 t_1^* \geq t t^*,$ and it is
not unitary because in fact $s_0^* s_0 > t_1 t_1^*.$
\end{pff}

The following lemma will be used to control the completely bounded
norms of perturbations of maps on finite dimensional operator spaces.
For now, we only need to know
that the map $W$ of the lemma satisfies
$\| W \|_{\cb},$ $\| W^{-1} \|_{\cb} < 1 + \ep$ for $\|a_l - b_l \|$
small enough, which is
contained in the proof of Proposition 2.6 of \cite{JP}.
The more explicit estimate will be relevant in Section 6.

\begin{lem}     
Let $A$ be a unital
\ca, let $a_1, \dots, a_m \in A$ be linearly independent,
and assume that $E = \spn (a_1, \dots, a_m)$ is unital and
selfadjoint (and hence an operator system in $A$).
Define
\[
M = \sup \left\{ \max_{1 \leq l \leq m} | \af_l | :
    \left\| \sum_{l = 1}^m \af_l a_l \right\| \leq 1 \right\}.
\]
Then for $b_1, \dots, b_m \in A,$ the linear map
$W : E \to \spn (b_1, \dots, b_m),$
given by $W (a_l) = b_l,$ satisfies
\[
\| W \|_{\cb} \leq 1 + m M \sum_{l = 1}^m \|a_l - b_l \|
\]
and, if $m M \sum_{l = 1}^m \|a_l - b_l \| < 1,$ then
\[
\| W^{-1} \|_{\cb} \leq
        \left( 1 - m M \sum_{l = 1}^m \|a_l - b_l \| \right)^{-1}.
\]
\end{lem}

\begin{pff}
Give $\C^m$ the norm
$\| (\af_1, \dots, \af_m ) \|_{\infty}
                 = \max_{1 \leq l \leq m} | \af_l |.$
Define $Q : E \to \C^m$ by sending $a_l$ to the
$l$-th standard basis vector $\xi_l,$ and define
$R : \C^m \to A$ by $R (\xi_l) = b_l - a_l.$ Then
$\| Q \| = M$ and $\| R \| \leq \sum_{l = 1}^m \|a_l - b_l \|.$
We have, using Lemma 2.3 of \cite{EH},
\[
\| R \circ Q \|_{\cb} \leq m \| R \circ Q \| \leq
        m M \sum_{l = 1}^m \|a_l - b_l \|.
\]
Since $W (a) = a + R (Q (a)),$ the
estimate on $\| W \|_{\cb}$ is now immediate. We further have, for
any $n$ and any $a \in \Mt{n} E,$
\[
\| (\id_{M_n} \otimes W )(a) \| \geq
       \| a \| - \| [\id_{M_n} \otimes (R \circ Q)] (a) \|
   \geq \| a \| (1 - \| R \circ Q \|_{\cb}).
\]
Therefore
\[
\| (\id_{M_n} \otimes W)^{-1} \| \leq
        \left( 1 - m M \sum_{l = 1}^m \|a_l - b_l \| \right)^{-1}
\]
for all $n.$
\end{pff}

We will need the following variant of an argument from one of the
proofs of \cite{Kr1}.

\begin{lem} 
Let $A$ be a unital \ca, let $E \subset A$ be an operator system,
let $H$ be a Hilbert space, and let $S : E \to L (H)$ be a unital
selfadjoint completely bounded map. Then there exists a
\ucp\  map $T : A \to L (H)$ such that
$\| T |_E - S \|_{\cb} \leq  \| S \|_{\cb} - 1.$
\end{lem}

\begin{pff}
Wittstock's generalization of the Arveson extension theorem
(see Theorem 7.2 of \cite{Pl}) provides a linear map
$Q : A \to L (H)$ such that $\| Q \|_{\cb} = \| S \|_{\cb}$ and
$Q |_E = S.$ Replacing $Q$ by
$x \mapsto \half (Q (x) + Q (x^*)^*),$ we may assume in addition that
$Q$ is selfadjoint. Now apply Proposition 1.19 of \cite{Ws}
(see also the original version in the proof of Theorem 4.1 of
\cite{Kr1}) to obtain a \ucp\  map $T : A \to L (H)$ such that
$\| T - Q \|_{\cb} \leq  \| Q \|_{\cb} - 1.$
Then also $\| T |_E - S \|_{\cb} \leq  \| S \|_{\cb} - 1.$
\end{pff}

The following lemma is the second main technical result of this section.

\begin{lem}     
Let $A$ be a separable unital exact \ca,
let $E$ be a \fd\  operator system in $A,$ and let $\ep > 0.$
For every $\dt < \ep$ there exists an integer
$n$ such that whenever $B_1$ and $B_2$ are separable unital \ca s, with
$B_2$ nuclear, and $V : E \to B_1$ and $W : E \to B_2$ are two
\ucp\  maps  such that $V$ is injective and $V^{-1} : V (E) \to E$
satisfies
$\| V^{-1} \otimes \id_{M_n}\| \leq 1 + \dt,$ then there is a
\ucp\  map $T : B_1 \to B_2$ such that $\|T \circ V - W \| < \ep.$
\end{lem}

\begin{pff}
Let $\rh = ( \ep - \dt ) / [2 (1 + \dt)] > 0.$ Since $A$ is exact,
it has a nuclear embedding in $L(H)$ for some Hilbert space
$H.$ (See \cite{Kr1}, \cite{Kr2}, or Theorem 9.1 of \cite{Ws}.)
We may thus assume that $A$ is a unital subalgebra of $L(H)$ with
the inclusion map nuclear.
Let $\{a_1, \dots, a_m\}$
be a basis for $E.$ Choose $\mu > 0$ small enough (using the previous
lemma) that if $b_1, \dots, b_m \in L(H)$ satisfy
$\| a_l - b_l \| < \mu,$ then the map $T (a_l) = b_l,$ from
$E$ to $\spn (b_1, \dots, b_m),$ satisfies
$\| T^{-1} \|_{\cb} < 1 + \rh.$ By nuclearity of the inclusion, there
are $n$ and \ucp\  maps $S_1 : E \to M_n$ and
$\tilde{S_2} : M_n \to L(H)$ such that the elements
$b_l = \tilde{S_2} \circ S_1 (a_l)$ satisfy $\| a_l - b_l \| < \mu.$
Let $T$ be as above, using this choice of $b_1, \dots, b_m,$ and let
$F = S_1 (E),$ which is an operator space in $M_n.$ Define
$S_2 : F \to E$ by $S_2 = T^{-1} \circ \tilde{S_2}.$ Then $S_2$ is
unital, $S_2 \circ S_1 = \id_E,$ and  $\|S_2\|_{\cb} < 1 + \rh.$ 
Moreover, from $S_1 (x^*) = S_1 (x)^*$ for $x \in E,$ we get
$S_2 (y^*) = S_2 (y)^*$ for $y \in F.$

Further choose, using the nuclearity of $B_2,$ an
integer $r$ and \ucp\  maps $W_1 : E \to M_r$ and
$W_2 : M_r \to B_2$ such that $\|W_2 \circ W_1 - W \| < \rh.$ Since
$F$ is a subspace of $M_n$ and $\| W_1 \circ S_2\|_{\cb} < 1 + \rh,$
Lemma 1.9
provides a \ucp\  map ${Q} : M_n \to M_r$ such that
$\| {Q} |_{F} - W_1 \circ S_2\| < \rh.$

Now consider $S_1 \circ V^{-1},$ which is a linear map from
$V(E)$ to $M_n.$ Since $S_1$ is \uacp, we have
\[
\| (S_1 \circ V^{-1}) \otimes \id_{M_n}\| \leq
  \| S_1 \|_{\cb} \| V^{-1} \otimes \id_{M_n}\| \leq 1 + \dt.
\]
By Proposition 7.9 of \cite{Pl}, this implies that
$\| S_1 \circ V^{-1} \|_{\cb} \leq 1 + \dt.$ Also, $S_1 \circ V^{-1}$
is unital and selfadjoint. Applying Lemma 1.9
again, we obtain a \ucp\  map $R : B_1 \to M_n$ such that
$\| R |_{V (E)} - S_1 \circ V^{-1}\| \leq \dt.$ We then have
\beqr
\lefteqn{
\| W_2 \circ {Q} \circ R |_{V (E)} - W \circ V^{-1}\|
              }     \\
 & \leq & \| W - W_2 \circ W_1 \| \| V^{-1}\| +
\| W_2 \| \|Q \circ R |_{V (E)} - W_1 \circ S_2 \circ S_1 \circ V^{-1}\|
                                                      \\
 & \leq & \rh (1 + \dt) + \| R |_{V (E)} - S_1 \circ V^{-1}\| +
             \| {Q} |_{F} - W_1 \circ S_2\| \|S_1 \circ V^{-1}\|
                                                      \\
 & <  & \rh (1 + \dt) +  \dt + \rh ( 1 + \dt ) = \ep.
\eeqr
Thus, $T = W_2 \circ {Q} \circ R$ is a \ucp\  map from $B_1$ to $B_2$
whose restriction to $V (E)$ differs in norm from $W \circ V^{-1}$ by
less than $\ep.$
\end{pff}

The following result is known in the
important special case in which the right hand side of the inequality
is zero. It is easy to give an example to show that the exponent
$\half$ can't be improved, but this also follows from the example after
the next (more significant) lemma. It turns out that it
is even impossible to improve the constant $\sqrt{2}.$

\begin{lem} 
Let $A$ be a unital \ca, let $u \in A$ be unitary, and
let  $s \in A$ be an isometry with range projection $e = s s^*.$
Then
\[
\| u - [e u e + (1 - e) u (1 - e)] \| \leq
       \inf \{ (2 \| s^* u s - v \|)^{1/2} : v \in A \,\,
                            {\mathrm{unitary}} \}.
\]
\end{lem}

\begin{pff}
We prove that if $v \in A$ is any other unitary, then
\[
\| u - [e u e + (1 - e) u (1 - e)] \| \leq \sqrt{2 \| s^* u s - v \|}.
\]
The element $s v s^*$ is a unitary in $e A e$ and
\[
\| e u e - s v s^*\| = \| s^* s u s^* s - s v s^*\| \leq
            \| s^* u s - v \|.
\]
So $\| (e u e)^* (e u e) - e \| \leq 2 \| s^* u s - v \|.$ Now
\[
e = e u^* u e = (e u e)^* (e u e) + [(1 - e) u e]^* [(1 - e) u e].
\]
Therefore
\[
\| [(1 - e) u e]^* [(1 - e) u e] \| \leq 2 \| s^* u s - v \|,
\]
so that $\|(1 - e) u e \| \leq \sqrt{2 \| s^* u s - v \|}.$
Similarly, using $u u^* = 1$ instead of $u^* u = 1,$ we get
$\| e u (1 - e) \| \leq \sqrt{2 \| s^* u s - v \|}.$
Since $e$ is orthogonal to $1 - e,$ it follows that
\[
\| u - [e u e + (1 - e) u (1 - e)] \| =
        \|(1 - e) u e +  e u (1 - e) \| \leq \sqrt{2 \| s^* u s - v \|}.
\]
\end{pff}

\begin{lem} 
Let $A$ be a unital \ca, let $s$ and $t$ be two isometries in $A,$
and let $D$ be a unital subalgebra of $A$ which is isomorphic to
$\OA{2}$ and such that every element of $D$ commutes with
$s$ and $t.$ Then there is a unitary $z \in A$ such that
whenever $u$ and $v$ are unitaries in $A$ commuting with every element
of $D,$ then
\[
\| z^* u z - v \| \leq
 11 \left[ \rule{0em}{3ex}
            \max (\| s^* u s - v \|, \|t^* v t - u \|) \right]^{1/2}.
\]
\end{lem}

\begin{pff}
Let $B$ be the relative commutant of $D$ in $A.$ Then $s$ and $t,$
along with all possible choices of $u$ and $v,$ are in $B.$
Since $\OA{2}$ is nuclear, there is a \hm\  from $\OT{B}$ to $A$
which is the identity on $B$ and sends $\OA{2}$ to $D.$
Therefore we may as well take $A = \OT{B},$ with $s, \, t \in B.$
We have to show that there is a unitary $z \in \OT{B}$ such that
whenever $u$ and $v$ are unitaries in $B,$ then
\[
\| z (1 \otimes u) z^* - 1 \otimes v \| \leq
 11 \left[ \rule{0em}{3ex}
            \max (\| s^* u s - v \|, \|t^* v t - u \|) \right]^{1/2}.
\]

The unitary $z$ will chop $1 \otimes u$ in pieces and reassemble them in
a different way. To construct it, we start by defining an assortment of
\pj s and partial isometries.
Define
\[
e_1 = s s^* \andeqn f_1 = t t^*,
\]
and further define
\[
e_2 = s f_1 s^* \leq e_1, \,\,\,\,\,\, f_2 = t e_1 t^* \leq f_1, \andeqn
              f_3 = t e_2 t^* \leq f_2.
\]
Then define two sets of \mops\  summing to $1$ by
\[
p_1 = 1 - e_1, \,\,\,\,\,\, p_2 = e_1 - e_2, \andeqn p_3 = e_2
\]
and
\[
q_1 = 1 - f_1, \,\,\,\,\,\, q_2 = f_1 - f_2, \,\,\,\,\,\,
        q_3 = f_2 - f_3, \andeqn q_4 = f_3.
\]
Next, construct partial isometries
\[
c_1 = p_2 s q_1, \,\,\,\,\,\, c_2 = p_1 t^* q_2, \,\,\,\,\,\,
            c_3 = p_2 t^* q_3, \andeqn c_4 = p_3 t^* q_4.
\]
One checks that
\[
c_1^* c_1 = q_1 \andeqn c_1 c_1^* = p_2,
\]
while
\[
c_j^* c_j = q_j \andeqn c_j c_j^* = p_{j - 1}
\]
for $j = 2, \, 3, \, 4.$ Let $s_1$ and $s_2$ be the standard
generators of $\OA{2}.$ Define
\[
z = s_1 \otimes c_1 + 1 \otimes c_2 + s_2 \otimes c_3 + 1 \otimes c_4,
\]
which is easily checked to be a unitary in $\OT{B}$ such that
\[
z (1 \otimes q_1) z^* = s_1 s_1^* \otimes p_2, \,\,\,
z (1 \otimes q_2) z^* = 1 \otimes p_1, \,\,\,
z (1 \otimes q_3) z^* = s_2 s_2^* \otimes p_2,
                            \,\,\, {\mathrm{and}} \,\,\,
z (1 \otimes q_4) z^* = 1 \otimes p_3.
\]

Now let $u,\, v \in B$ be unitaries.
Set $\dt = \max (\| s^* u s - v \|, \|t^* v t - u \|).$ Using
$s q_1 = p_2 s$ and $p_2 \leq s s^*,$ we obtain
\beqr
\| c_1 (q_1 v q_1) c_1^* - p_2 u p_2\| & = &
                         \| p_2 s v s^* p_2 - p_2 u p_2\|         \\
 & = & \| p_2 s v s^* p_2 - p_2 s s^* u s s^* p_2\|
                     \leq \| v - s^* u s\| \leq \dt.
\eeqr
Similarly, for $j = 2, \, 3, \, 4$ one gets
\[
\| c_j (q_j v q_j) c_j^* - p_{j - 1} u p_{j - 1}\|
         \leq \| t^* v t - u\| \leq \dt.
\]
One checks that
\beqr
\lefteqn{z [1 \otimes (q_1 v q_1 +
                  q_2 v q_2 + q_3 v q_3 + q_4 v q_4)] z^*}   \\
 & = & s_1 s_1^* \otimes c_1 (q_1 v q_1) c_1^*
   + 1 \otimes c_2 (q_2 v q_2) c_2^*
   +  s_2 s_2^* \otimes c_3 (q_3 v q_3) c_3^*
   + 1 \otimes c_4 (q_4 v q_4) c_4^*,
\eeqr
so that (using $s_1 s_1^* + s_2 s_2^* = 1$)
\[
\| z (1 \otimes v) z^* - 1 \otimes u\| \leq
        \dt + \|q_1 v q_1 + q_2 v q_2 + q_3 v q_3 + q_4 v q_4 - v \| +
                  \|p_1 u p_1 + p_2 u p_2 + p_3 u p_3 - u\|.
\]

We now have to estimate the last two terms in the last inequality above.
Recall that $e_1 = s s^*,$ so that
\[
\| u - [e_1 u e_1 + (1 - e_1) u (1 - e_1)]\| \leq 
     \sqrt{2 \| s^* u s - v \|} \leq \sqrt{2 \dt}
\]
by Lemma 1.11.
Since $e_2 = s t t^* s^*,$ we get
\[
\| e_2 u e_2 - s t u t^* s^* \| \leq \| t^* s^* u s t - u \|
        \leq \|s^* u s - v \| + \| t^* v t - u\| \leq 2 \dt.
\]
So
\[
\| u - [e_2 u e_2 + (1 - e_2) u (1 - e_2)]\| \leq \sqrt{4 \dt},
\]
again by Lemma 1.11.
Compressing by $e_1 \geq e_2,$ we get
\[
\| e_1 u e_1 - [e_2 u e_2 + (e_1 - e_2) u (e_1 - e_2)]\|
            \leq \sqrt{4 \dt}.
\]
Recalling the definitions of $p_1,$ $p_2,$ and $p_3,$ it follows that
\[
\|p_1 u p_1 + p_2 u p_2 + p_3 u p_3 - u\|
                    \leq \left(\sqrt{2} + \sqrt{4}\right) \sqrt{\dt}.
\]

A similar sequence of estimates, with one more step and using the
equations $f_1 = t t^*,$ $f_2 = t s s^* t^*,$ and
$f_3 = (t s t) (t s t)^*,$ gives
\[
\|q_1 v q_1 + q_2 v q_2 + q_3 v q_3 + q_4 v q_4 - v \|
         \leq \left(\sqrt{2} + \sqrt{4} + \sqrt{6}\right) \sqrt{\dt}.
\]
We note that necessarily $\dt \leq 2,$ so that $\dt \leq \sqrt{2\dt}.$
We can therefore put our estimates together to get
\[
\|z^* (1 \otimes u) z - 1 \otimes v\| =
   \| z (1 \otimes v) z^* - 1 \otimes u\| \leq
 \left[\sqrt{2} + \left(\sqrt{2} + \sqrt{4}\right) +
        \left(\sqrt{2} + \sqrt{4} + \sqrt{6}\right) \right] \sqrt{\dt}
  \leq 11 \sqrt{\dt}.
\]
\end{pff}

The exponent $\half$ in this lemma can't be improved, as we now show by
example, even if $z$ is allowed to depend on $u$ and $v.$
It follows from the proof of the lemma that the  exponent $\half$
can't be improved in the previous lemma either. (We will also see a
more indirect proof of this below: any improvement here would imply a
a corresponding improvement in the exponent in Proposition 4.13,
but Remark 6.11
shows that no improvement is possible there.)

\begin{exa} 
Let $D = M_2 \otimes \OA{2}.$ Let $s_1$ and $s_2$ be the two standard
generating isometries of $\OA{2},$ and set $p_j = s_j s_j^*.$
Let $\af \in \R,$ and define a unitary $w_{\af} \in M_2$ by
\[
w_{\af}  = \left( \begin{array}{cc} \cos (\af) & \sin (\af) \\
                        - \sin (\af) & \cos (\af) \end{array} \right).
\]
Use the fact that any two nonzero \pj s in $M_2 \otimes \OA{2}$ are
\mvn\  to choose $c \in M_2 \otimes \OA{2}$ such that
\[
c^* c = 1 \otimes p_2 \andeqn
c c^* = \left( \begin{array}{cc} 1 & 0 \\ 0 & 0 \end{array} \right)
                                       \otimes p_2.
\]
Then set
\[
u = w_{\af} \otimes 1, \,\,\,\,\,\,
v = w_{\af} \otimes p_1 + 1 \otimes p_2, \,\,\,\,\,\,
s = 1 \otimes p_1 + c, \andeqn t = 1 \otimes s_1.
\]
Then
\[
s^* u s = w_{\af} \otimes p_1 + \cos (\af) (1 \otimes p_2) \andeqn
t^* v t = u.
\]
Therefore
\[
\max (\| s^* u s - v \|, \|t^* v t - u \|) = 1 - \cos (\af).
\]
Also
\[
\spec (u) = \{\exp (i \af), \exp (- i \af)\} \andeqn
\spec (v) = \{\exp (i \af), 1, \exp (- i \af)\}.
\]
Therefore, if $A$ is any \ca\  at all which contains $D$ as a
unital subalgebra, and if $z$ is any unitary in $A,$ we have
\[
\|z^* u z - v\| \geq | \exp (i \af) - 1 | =
          \sqrt{2} (1 - \cos (\af))^{1/2}.
\]
\QED
\end{exa}

Our first application of the technical results of this section is the
following lemma and theorem.

\begin{lem} 
Let $A$ be a separable unital exact \ca, and let $B$ be a \kfalg.
Let $\ph,$ $\ps: A \to B$ be two injective unital \hm s. Then the
\hm s from $A$ to $\OA{2} \otimes B,$ given by
$a \mapsto 1 \otimes \ph (a)$ and $a \mapsto 1 \otimes \ps (a),$
are \ayue.
\end{lem}

\begin{pff}
Let $u_1, \dots, u_n \in A$ be unitaries, and let $\ep > 0.$
We prove that there is a unitary $z \in \OT{B}$ such that
\[
\| z ( 1 \otimes \ph (u_j)) z^* - 1 \otimes \ps (u_j) \| < \ep
\]
for $1 \leq j \leq n.$ Let
\[
E = \spn \{1, u_1, u_1^*, \dots, u_n, u_n^* \},
\]
which is a \fd\  operator system. Lemma 1.10
applies to $\ph |_E$ and $\ps |_E$ (with $\dt = 0$), and
provides \ucp\  maps $S,$ $T : B \to B$ such that
\[
\| S \circ \ph |_E - \ps |_E \| <
              \frac{1}{2} \left(\frac{\ep}{11} \right)^2 \andeqn
\| T \circ \ps |_E - \ph |_E \| <
              \frac{1}{2} \left(\frac{\ep}{11} \right)^2.
\]
In particular,
\[
\| S (\ph (u_j)) - \ps (u_j) \| <
             \frac{1}{2} \left(\frac{\ep}{11} \right)^2 \andeqn
\| T (\ps (u_j)) - \ph (u_j) \| < 
             \frac{1}{2} \left(\frac{\ep}{11} \right)^2
\]
for $1 \leq j \leq n.$ Proposition 1.7
then gives isometries $s,$ $t \in B$ such that
\[
\| s^* \ph (u_j) s - \ps (u_j) \| <
                    \left(\frac{\ep}{11} \right)^2 \andeqn
\| t^* \ps (u_j) t - \ph (u_j) \| < \left(\frac{\ep}{11} \right)^2
\]
for $1 \leq j \leq n.$
The existence of the required unitary $z \in \OT{B}$ now
follows from Lemma 1.12.
\end{pff}

As a corollary, we obtain:

\begin{thm}  
Let $A$ be a separable unital exact \ca. Then any two
injective unital \hm s from $A$ to $\OA{2}$ are
\ayue\  (Definition 0.5).
\end{thm}

\begin{pff}
Let $\ph,$ $\ps : A \to \OA{2}$ be injective and unital. Let
$\mu : \OT{\OA{2}} \to \OA{2}$ be an isomorphism (from Theorem 0.8),
and let $\bt : \OA{2} \to \OT{\OA{2}}$ be $\bt (a) = 1 \otimes a.$
Then $\mu \circ \bt$ is \ayue\  to $\id_{\OA{2}}$ by
Proposition 0.7.
Also, $\bt \circ \ph$ is \ayue\  to $\bt \circ \ps$ by
Lemma 1.14.
Thus $\ph$ is \ayue\  to $\mu \circ \bt \circ \ph,$ which
is \ayue\  to $\mu \circ \bt \circ \ps,$ which in turn
is \ayue\  to $\ps.$ 
\end{pff}

\section{Embedding in $\OA{2}$}

The main result of this section is that every separable unital exact
\ca\  can be unitally embedded in $\OA{2}.$ The algebra of bounded
sequences modulo sequences converging to zero
plays an important role in the proof,
so we begin by establishing notation concerning it.

\begin{ntn} 
For any \ca\  $D,$ we denote by $\LI{D}$ the set of bounded
sequences $d = (d_1, d_2, \dots)$ with values in $D.$ It is a \ca\  with
the obvious operations and norm.
For compatibility with the notation for ultrapowers below, we define
$c_{\infty} (D) = C_0 (\N) \otimes D \subset \LI{D}$ and
$D_{\infty} = \LI{D} / c_{\infty} (D).$ We denote by
$\pi_{\infty}^{(D)}$ (or $\pi_{\infty}$ when $D$ is clear from the
context) the quotient map $\LI{D} \to D_{\infty}.$
\end{ntn}

The technical results of the previous section imply fairly directly that
if $A$ is separable and exact, and has an embedding in
$(\OA{2})_{\infty}$ which lifts to a \ucp\  map to $\LI{\OA{2}},$
then $A$ has an embedding in $\OA{2}.$ In particular, this applies to
quasidiagonal exact \ca s. The rest of the section is devoted to
the extension from the quasidiagonal case to the general case.
It is possible to embed any separable \ca\  in a crossed product
of a quasidiagonal separable \ca\  by $\Z,$ and the method is to show
that crossed products by $\Z$ of exact \ca s embeddable in
$\OA{2}$ are again embeddable in $\OA{2}.$

\begin{lem} 
Let $A$ be a unital separable exact \ca. If there is a unital injective
\hm\  from $A$ to
$(\OA{2})_{\infty} = \LI{\OA{2}} / c_{\infty} (\OA{2})$ which has a
lifting to a \ucp\  map from $A$ to $\LI{\OA{2}},$
then there is an injective unital \hm\  from $A$ to $\OA{2}.$
\end{lem}

\begin{pff}
Let $\ph : A \to (\OA{2})_{\infty}$ be a unital injective \hm\  which
has a lifting as in the hypotheses of the lemma.
Let $u_1, u_2, \dots$ be a sequence of unitaries in $A$ whose linear
span is dense in $A.$ Define \fd\  operator systems $E_n$ by
$E_n = \spn \{1, u_1, u_1^*, \dots, u_n, u_n^*\}.$
Then $\C \cdot 1 = E_0 \subset E_1 \subset \cdots \subset A$ and
$\overline{\bigcup_{n = 0}^{\infty} E_n} = A.$
We first show that there is a unital injective
\hm\  $\ps : A \to (\OA{2})_{\infty}$ with a \ucp\   lifting
$a \mapsto V(a) = (V_1 (a), V_2 (a), \dots)$ from $A$ to $\LI{\OA{2}}$
with the following property: For each fixed $n,$ for
all sufficiently large $m$ (how large depends on $n$), the restriction
$V_m |_{E_n}$ is injective, and furthermore its inverse, defined on
$V_m (E_n),$ satisfies
$\limi{m} \| (V_m |_{E_n} )^{-1} \otimes \id_{M_k} \| = 1$ for all
$k \in \N.$

To do this, let $a \mapsto Q(a) = (Q_1 (a), Q_2 (a), \dots)$
be a lifting of $\ph$ to a \ucp\  map from $A$ to $\LI{\OA{2}}.$
Injectivity of $\ph$ does not imply that
$\limi{m} \| (Q_m \otimes \id_{M_k}) (a) \| = \| a \|,$
but we can remedy this by grouping the $Q_m$ together in blocks.
Injectivity of $\ph$ does imply that for every $N \in \N,$ the map
$a \mapsto \ph^{(N)} (a) =
                     \pi_{\infty} (Q_{N + 1} (a), Q_{N + 2} (a), \dots)$
is again an injective \hm. Therefore,
for every $N,$ $k \in \N$ and $a \in \Mt{k} A,$ we have
\[
\limi{m} \|((Q_{N + 1} \otimes \id_{M_k}) (a), \dots,
                 (Q_{N + m} \otimes \id_{M_k}) (a)) \| =
 \|(\ph^{(N)} \otimes \id_{M_k}) (a) \| = \| a \|.
\]
Since each $E_n$ is \fd, we can therefore
construct by induction a sequence
\[
0 = N_1 < N_2 < \cdots < N_m < N_{m + 1} < \cdots
\]
of integers such that
\[
\|((Q_{N_m + 1} \otimes \id_{M_k}) (a), \dots,
                 (Q_{N_{m + 1}} \otimes \id_{M_k}) (a)) \|
         \geq (1 - 2^{-m}) \| a \|
\]
for $k \leq m$ and $a \in \Mt{k} E_m.$
Let $\sm_m : \OA{2}^{N_{m + 1} - N_m} \to \OA{2}$
be any unital \hm.
Define $V_m : A \to \OA{2}$ by
\[
V_m (a) = \sm_m ((Q_{N_m + 1} (a), \dots, Q_{N_{m + 1}} (a)).
\]
Note that $V_m$ is \uacp\   because each $Q_j$ is, so that also
$V(a) = (V_1 (a), V_2 (a), \dots)$ defines a \ucp\  map
from $A$ to $\LI{\OA{2}}.$ By construction we have
$\limi{m} \| (V_m |_{E_n} )^{-1} \otimes \id_{M_k} \| = 1$ for each
fixed $k,$ $n \in \N.$
Taking $k = 1,$ we see that the linear map
$\ps = \pi_{\infty} \circ V$ is isometric, hence injective.
Since $\limi{j} (Q_j (ab) - Q_j (a) Q_j (b)) = 0$ for $a,$ $b \in A,$
we also obtain $\limi{m} (V_m (ab) - V_m (a)V_m (b)) = 0$
for $a,$ $b \in A.$ Therefore $\ps$ is a \hm, and its lifting
$V$ satisfies the required conditions.

Choose numbers $\dt_m > 0$ such that $\dt_0 \geq \dt_1 \geq \cdots$
and $2 \dt_m + 11 \sqrt{5 \dt_m} < 2^{-m}.$ Using Lemma 1.10,
choose positive integers $k(m)$ with
$k(0) \leq k(1) \leq \cdots$ and such that whenever
$V,$ $W : E_m \to \OA{2}$ are \ucp\  with $V$ injective and
$\|V^{-1} \otimes \id_{M_{k (m)}} \| \leq 1 + \dt_m,$ then there is a
\ucp\  map $T: \OA{2} \to \OA{2}$ such that
$\|T \circ V - W \| < 2 \dt_m.$
The conditions on $V$ imply that we can pass to a subsequence in the
variable $m$ in such a way that $V_m |_{E_n}$ is injective for
$n \leq m,$ and moreover we have the estimates
\[
\| (V_m |_{E_n} )^{-1} \otimes \id_{M_{k (m)}} \| \leq 1 + \dt_m,
  \,\,\,\,\,  \| V_m (u_n)^* V_m (u_n) - 1 \| < \dt_m,
       \,\,\,\,\, {\mathrm{and}} \,\,\,\,\,
         \| V_m (u_n) V_m (u_n)^* - 1 \| < \dt_m
\]
for all $m$ and all $n \leq m.$

Using these estimates and Lemma 1.10,
find \ucp\  maps $S_m,$ $T_m : \OA{2} \to \OA{2}$ such that
\[
\| T_m \circ V_m |_{E_m} - V_{m + 1} |_{E_m} \| \leq 2 \dt_m \andeqn
\| S_m \circ V_{m + 1} |_{E_m} - V_m |_{E_m} \| \leq 2 \dt_m.
\]
For $1 \leq j \leq m$ define unitaries
$x_m^{(j)} = V_m (u_j) [V_m (u_j)^* V_m (u_j)]^{-1/2}.$ Then
$\|x_m^{(j)} - V_m (u_j) \| \leq \dt_m$ by Lemma 0.1.
It follows that
\[
\| T_m (x_m^{(j)}) - x_{m + 1}^{(j)} \| \leq 4 \dt_m \andeqn
\| S_m (x_{m + 1}^{(j)}) - x_m^{(j)} \| \leq 4 \dt_m.
\]
Proposition 1.7
gives isometries $s_m,$ $t_m \in \OA{2}$
(depending on $m$) such that
\[
\| s_m^* x_m^{(j)} s_m - x_{m + 1}^{(j)} \| \leq 5 \dt_m \andeqn
\| t_m^* x_{m + 1}^{(j)} t_m - x_m^{(j)} \| \leq 5 \dt_m
\]
for $1 \leq j \leq m.$ Lemma 1.12
now gives unitaries
$z_m \in \OA{2} \otimes \OA{2}$ such that
\[
\|z_m (1 \otimes x_m^{(j)}) z_m^* - 1 \otimes x_{m + 1}^{(j)} \| \leq
         11 \sqrt{5 \dt_m}.
\]
Thus
\[
\|z_m (1 \otimes V_m (u_j) ) z_m^* - 1 \otimes V_{m + 1} (u_j) \| \leq
        2 \dt_m + 11 \sqrt{5 \dt_m} < 2^{-m}
\]
for $1 \leq j \leq m.$

Now define $y_n = z^*_1 z^*_2 \cdots z^*_n.$ Then the $y_n$ are
unitaries such that $\limi{n} y_n (1 \otimes V_n (u_j) ) y_n^*$
exists for all $j.$ It follows that the limit
$\ps_0 (a) = \limi{n} y_n (1 \otimes V_n (a) ) y_n^*$ exists for
all $a\in \bigcup_{n = 1}^{\infty} E_n.$ Furthermore, for each $n$
and $m,$
the map $V_m |_{E_n}$ is \uacp.
Therefore $\ps_0 |_{E_n}$ is \uacp, whence $\| \ps_0 |_{E_n} \| \leq 1.$
It follows that $\ps_0$ extends by continuity to a \ucp\  map
$\ps : A \to \OA{2} \otimes \OA{2}.$ Since
$\limi{m} (V_m (ab) - V_m (a) V_m (b)) = 0$ for all
$a \in \bigcup_{n = 1}^{\infty} E_n,$ it follows that $\ps$ is actually
a \hm. Finally, for $a \in \bigcup_{n = 1}^{\infty} E_n$ we have
$\| \ps (a) \| = \limi{n} \| V_m (a)\| = \| a \|,$ from which it
follows that $\ps$ is isometric and hence injective.

We thus have an injective unital \hm\   from $A$ to
$\OA{2} \otimes \OA{2}.$ The existence of an injective unital
\hm\   from $A$ to $\OA{2}$ now follows from
the isomorphism $\OA{2} \otimes \OA{2} \cong \OA{2}$ (Theorem 0.8).
\end{pff}

Recall that a separable \ca\   $A$ is called {\em quasidiagonal}
(weakly
quasidiagonal in some papers) if there is an injective representation
$\pi$ of $A$ on a separable Hilbert space $H$ and a sequence
$p_1 \leq p_2 \leq \cdots$ of finite rank \pj s on $H$ such that
$p_n \to 1$ in the strong operator topology and
$\limi{n} \| p_n \pi (a) - \pi (a) p_n \| = 0$ for all $a \in A.$

\begin{cor} 
Let $A$ be a separable unital exact quasidiagonal \ca. Then there
exists a unital injective \hm\  from $A$ to $\OA{2}.$
\end{cor}

\begin{pff}
By the previous lemma, it suffices to find a unital injective \hm\  from
$A$ to $\LI{\OA{2}} / c_0 (\OA{2})$ which has a lifting to a \ucp\  map
from $A$ to $\LI{\OA{2}}.$ Now $\OA{2}$ contains a unital copy of
every matrix algebra $M_k,$ so one can easily construct a unital
injective \hm\  from any product $\prod_{n = 1}^{\infty} M_{k(n)}$ to
$\LI{\OA{2}}.$ (The product $\prod_{n = 1}^{\infty} M_{k(n)}$ means,
of course, the \ca\  of all bounded sequences in the set-theoretic
product.) Therefore there is a unital injective \hm\  from
$\prod_{n = 1}^{\infty} M_{k(n)} / \bigoplus_{n = 1}^{\infty} M_{k(n)}$
to $\LI{\OA{2}} / c_0 (\OA{2}).$ So it suffices to
find a unital injective \hm\  from $A$ to
$\prod_{n = 1}^{\infty} M_{k(n)} / \bigoplus_{n = 1}^{\infty} M_{k(n)},$
with a lifting to a \ucp\  map
from $A$ to $\prod_{n = 1}^{\infty} M_{k(n)},$
for some sequence $k(1), k(2), \dots$ of positive integers. This
is easy to do, and is done in Proposition 3.1.3 and the preceding
remark in \cite{BK}.
\end{pff}

We now start our preparations for the general case. The following lemma
is a variant of a result of Effros and Haagerup
\cite{EH}, and will be used to construct the lifting in the
hypotheses of Lemma 2.2.

\begin{lem}     
Let $A$ and $B$ be separable unital \ca s, let $J$ be an ideal in
$A$ which is approximately injective in the sense of \cite{EH}
(definition before Lemma 3.3), and let $\ph : A \to B/J$ be an
injective \hm.
Let $H$ be a separable infinite dimensional Hilbert space, and
suppose that the induced map of algebraic tensor products
$A \talg L (H) \to ( B \talg L (H)) / ( J \talg L (H))$
extends continuously to a (necessarily injective) \hm\  %
\[
\overline{\ph} : A \tmin L (H) \to ( B \tmin L (H)) / ( J \tmin L (H)).
\]
Then there is a \ucp\  map $T : A \to B$ which lifts $\ph.$
\end{lem}

\begin{pff}
Let $\rh : B \to B/J$ be the quotient map.

We first reduce to the case $A = B/J$ and $\ph = \id_{B/J}.$
Let $B_0 = \rh^{-1} (\ph (A)) \subset B,$ and let $\rh_0 = \rh |_{B_0}.$
Since the minimal tensor product preserves inclusions, we have
$J \tmin L (H) \subset B_0 \tmin L (H) \subset B \tmin L (H).$
So the hypotheses of the lemma hold with $B_0$ in place of $B.$

We now assume $A = B/J.$
The desired conclusion will follow from Theorem 3.4 of \cite{EH},
provided we verify the hypothesis (b) there, that is, that for
every unital \ca\  $C,$ the
kernel of $\rh \otimes \id_C : B \tmin C \to (B / J) \tmin C$
is exactly $J \tmin C.$
The hypothesis of the lemma is that this is true for $C = L (H).$

If $C$ is separable, we may suppose $C \subset L (H),$
and use Proposition 2.6 of \cite{Ws}.
(Alternatively, combine the method below for reduction to
the separable case with Lemma 3.9 of \cite{Kr1}.)
For general $C,$ we need only show that
$\ker (\rh \otimes \id_C) \subset J \tmin C.$
Let $a \in \ker (\rh \otimes \id_C).$
Choose a separable subalgebra $C_0 \subset C$ such that
$a \in B \tmin C_0 \subset B \tmin C.$
Using $(B/J) \tmin C_0 \subset (B/J) \tmin C,$ we have
$\ker (\rh \otimes \id_C) =
   \ker (\rh \otimes \id_{C_0}) \cap B \tmin C_0,$
whence $a \in \ker (\rh \otimes \id_{C_0}).$
So $a \in J \tmin C_0 \subset J \tmin C$ by the separable case.
\end{pff}

We now turn to the crossed product part of the construction.

Let $\af : G \to \Aut (A)$ be an action of a discrete group $G$
on a \ca\  $A.$ By a {\em{covariant representation}} $(u, \ph)$
of the system $(G, A, \af)$ in a unital \ca\  $B,$ we mean the
obvious generalization of a covariant representation on a Hilbert
space. That is, $u : G  \to U(B)$ is a \hm\  from $G$ to the unitary
group $U(B)$ of $B,$ $\ph : A \to B$ is a \hm\  of \ca s, and
$u (g) \ph (a) u (g)^* = \ph (\af_g (a))$ for all $a \in A$ and
$g \in G.$

The following lemma is the easy generalization to this context of
standard results on regular representations of crossed products by
discrete amenable groups.

\begin{lem}  
Let $G$ be a discrete amenable group, and let
$\af : G \to \Aut (A)$ be an action of $G$ on a unital \ca\  $A.$
Let $(u, \ph)$ be a covariant representation of $(G, A, \af)$ in a
unital \ca\  $B,$ with $\ph$ injective.
For $g \in G,$ let $g$ also denote the corresponding elements of
the group \ca\  $C^* (G)$ and the crossed product
\ca\  $C^* (G, A, \af).$ Then there is an injective
\hm\  $\ps : C^* (G, A, \af) \to C^* (G) \otimes B$ determined
by $\ps (a) = 1 \otimes \ph (a)$ for $a \in A$ and
$\ps (g) = g \otimes u(g)$ for $g \in G.$
\end{lem}

\begin{pff}
The existence (and uniqueness) of $\ps$ is immediate from the
universal property of the crossed product. We have to prove
injectivity.

Let $\pi_0: B \to L(H_0)$ be an injective representation of
$B$ on a Hilbert space $H_0,$ and let $\ld$ be the regular
representation of $C^* (G)$ on the Hilbert space $l^2 (G).$ Then
$\sm = (\ld \otimes \pi_0) \circ \ps$ is a representation of
$C^* (G, A, \af)$
on the  Hilbert space $H = l^2 (G) \otimes H_0.$ We will prove it
is unitarily equivalent to the regular representation  $\pi$ of
$C^* (G, A, \af)$ on $H$ associated with the
injective representation $\pi_0 \circ \ph$ of $A.$
(See Section 7.7 of \cite{Pd}.) This will prove injectivity of
$\ps,$ since $G$ is amenable (so that
$\pi$ is injective by Theorem 7.7.5 of \cite{Pd}).

The two representations are given by the formulas
\[
(\pi (a) \xi) (g) = (\pi_0 \circ \ph \circ \af_g^{-1}) (a) (\xi (g))
   \andeqn
(\sm (a) \xi) (g) = (\pi_0 \circ \ph) (a) (\xi (g))
\]
for $a \in A,$ $\xi \in H$ (viewed as $l^2 (G, H_0)$), and $g \in G,$
and
\[
\pi (g) =  \ld (g) \otimes 1 \andeqn
         \sm (g) = \ld (g) \otimes \pi_0 ( u(g))
\]
for $g \in G.$
Let $v$ be the unitary on $l^2 (G, H_0)$ given by
$(v \xi) (g) = \pi_0 ( u(g)) (\xi (g))$ for $\xi \in l^2 (G, H_0)$
and $g \in G.$ Then one can check directly that
\[
v \pi (a) v^* = \sm (a) \andeqn v \pi (g) v^* = \sm (g)
\]
for $a \in A$ and $g \in G.$ So $\sm$ is unitarily equivalent
to $\pi,$ as desired.
\end{pff}

In applications of the following lemma, the approximate innerness
assumption will be derived from Lemma 1.14.

\begin{lem}     
Let $B$ a unital \ca, let $A$ be a subalgebra of $B$
which contains the identity, and let $\sm \in \Aut (A).$ Suppose
that $\sm$ is approximately inner in $B,$ that is,
there is a sequence $v_1, v_2, \dots$ of unitaries in $B$ such that
$\limi{n} v_n a v_n^* = \sm (a)$ for all $a \in A.$ Let $z$ be the
standard generator of $C(S^1)$ and let $u$ be the canonical unitary
in $C^* (\Z, A, \sm)$ which implements $\sm$ on $A.$ Then the maps
\[
a \mapsto 1 \otimes \pi^{B} (a, a, \dots)          \andeqn
u \mapsto z \otimes \pi^{B} (v_1, v_2, \dots)
\]
define an injective
\hm\  $\ph :
     C^* (\Z, A, \sm) \to C(S^1) \otimes [\LI{B} / c_{\infty} (B)].$
Moreover, for any unital \ca\ $C,$ this \hm\  extends continuously to
an injective \hm\  %
\[
C^* (\Z, A, \sm) \tmin C \to
  C(S^1) \otimes [(\LI{B} \tmin C) / (c_{\infty} (B) \tmin C)].
\]
\end{lem}

\begin{pff}
We first show that the last sentence in the lemma follows from the rest.
Representing everything on Hilbert spaces and forming the spatial
tensor and crossed products, we easily see that
$C^* (\Z, A, \sm) \tmin C$ = $C^* (\Z, A \tmin C, \sm \otimes \id_C).$
Moreover, clearly
\[
\limi{n} (v_n \otimes 1) x (v_n \otimes 1)^* = (\sm \otimes \id_C) (x)
\]
for all $x \in A \tmin C.$ (Check on the algebraic tensor product.) The
first part of the lemma (applied to both $A$ and $A \tmin C$) therefore
implies that $\ph$ extends continuously to
an injective \hm\  %
\[
\overline{\ph} : C^* (\Z, A, \sm) \tmin C \to
  C(S^1) \otimes [\LI{B \tmin C} / c_{\infty} (B \tmin C)].
\]
Now $\LI{B} \tmin C$ is a subalgebra of $\LI{B \tmin C}.$ (Represent
$B$ faithfully on a Hilbert space $H_1,$ and represent $\LI{B}$
faithfully on $l^2 (\N) \otimes H_1$ in the obvious way. Represent
$C$ faithfully on another Hilbert space $H_2,$ and compare the
spatial tensor products as represented on
$l^2 (\N) \otimes H_1 \otimes H_2.$) Since
\[
c_{\infty} (B) \tmin C
           = c_{\infty} (\C) \tmin B \tmin C = c_{\infty} (B \tmin C),
\]
the inclusion of $\LI{B} \tmin C$ in $\LI{B \tmin C}$ gives an injective
\hm\    %
\[
[\LI{B} \tmin C]/[c_{\infty} (B) \tmin C] \to
                            \LI{B \tmin C} / c_{\infty} (B \tmin C).
\]
One immediately checks that
the range of $\overline{\ph}$ is contained in the image of
\[
C(S^1) \otimes
   \left[  \rule{0em}{2.2ex}
                [\LI{B} \tmin C] / [c_{\infty} (B) \tmin C] \right].
\]
This gives the desired extension.

We now prove the first part of the lemma.
The hypotheses immediately imply that
\[
(v_1, v_2, \dots) \cdot (a, a, \dots) \cdot (v_1, v_2, \dots)^*
           - ( \sm (a), \sm (a), \dots) \in c_{\infty} (B)
\]
for all $a \in A.$ Therefore
\[
a \mapsto \pi^{B} (a, a, \dots)         \andeqn
u \mapsto \pi^{B} (v_1, v_2, \dots)
\]
define a \hm\  from $C^* (\Z, A, \sm)$ to $\LI{B} / c_{\infty} (B).$
Moreover, $a \mapsto \pi^{B} (a, a, \dots)$ is injective. So
\[
a \mapsto 1 \otimes \pi^{B} (a, a, \dots)          \andeqn
u \mapsto z \otimes \pi^{B} (v_1, v_2, \dots)
\]
define an injective \hm\  from $C^* (\Z, A, \sm)$ to
$C(S^1) \otimes [\LI{B} / c_{\infty} (B)]$ by Lemma 2.5.
\end{pff}

\begin{lem}     
Let $B$ a separable nuclear unital \ca, let $A$ be a subalgebra of $B$
which contains the identity, and let $\sm \in \Aut (A)$ be
approximately inner in $B.$ Then the
\hm\  $C^* (\Z, A, \sm) \to C(S^1) \otimes [\LI{B} / c_{\infty} (B)]$
of the previous lemma
has a lifting to a \ucp\  map
$C^* (\Z, A, \sm) \to C(S^1) \otimes \LI{B}.$
\end{lem}

\begin{pff}
This now follows immediately from Lemma 2.4,
using
$C^* (\Z, A, \sm)$ in place of $A,$ $C(S^1) \otimes \LI{B}$
in place of $B,$ and $C(S^1) \otimes c_{\infty} (B)$ in place of $J.$
The ideal $C(S^1) \otimes c_{\infty} (B)$ is approximately injective
because it is nuclear, and the extension to a \hm\  %
\[
C^* (\Z, A, \sm) \tmin L(H) \to
[C(S^1) \otimes \LI{B} \tmin L(H)] /
                          [C(S^1) \otimes c_{\infty} (B) \tmin L(H)]
\]
is obtained from the previous lemma with $C = L(H).$
\end{pff}

\begin{thm}      
Let $A$ be a separable unital exact \ca. Then there exists an
injective unital \hm\  from $A$ to $\OA{2}.$
\end{thm}

\begin{pff}
The cone $C_0 ([0, 1)) \otimes A$ is quasidiagonal by Theorem 5
of \cite{Vc}.
The \ca\  $B_0 = (C_0 (\R) \otimes A)\unit$ is the unitization
of a subalgebra of $C_0 ([0, 1)) \otimes A,$ hence also quasidiagonal.
It is still exact by Proposition 7.1 (iii) and (vi) of \cite{Kr2}.
(Also compare with Remark 4.4 (4) of \cite{Ws}.)
Therefore Corollary 2.3
provides a unital embedding
$\ph_0 : B_0 \to \OA{2}.$
Let $B = C^* (\Z, B_0, \ta),$ where the action $\ta$ is by
translation on $\R$ and is trivial on $A.$
We produce an embedding of $B$ in
$(\OA{2})_{\infty} = \LI{\OA{2}} / c_{\infty} (\OA{2})$ which has a
lifting to a \ucp\  map from $B$ to $\LI{\OA{2}}.$

Let $\ta_1$ be the automorphism of $B_0$ which generates the action,
and let $\ps_0 = \ph_0 \circ \ta_1 : B_0 \to \OA{2}.$
Let $\mu : \OT{\OA{2}} \to \OA{2}$ be an isomorphism,
obtained from Theorem 0.8.
Define
$\ph,$ $\ps : B_0 \to \OA{2}$ by $\ph (a) = \mu (\ph_0 (a) \otimes 1)$
and $\ps (a) = \mu (\ps_0 (a) \otimes 1).$ Lemma 1.14
implies that $\ps$ is \ayue\  to $\ph.$ Thus,
using the embedding $\ph$ of $B_0$ in $\OA{2},$ the
automorphism $\ta_1$ is approximately inner in $\OA{2}$ in the sense of
Lemma 2.6.
That lemma therefore provides an injective \hm\  from
$B$ to $C(S^1) \otimes (\OA{2})_{\infty},$ which has a lifting to a
\ucp\  map from $B$ to $C(S^1) \otimes \LI{\OA{2}}$ by
Lemma 2.7.
It is easy to find an
injective \hm\  from $C(S^1)$ to the $2^{\infty}$ UHF algebra
$D,$ and it follows from Corollary 7.5 of \cite{Rr1} that
$D \otimes \OA{2} \cong \OA{2}.$ 
We thus obtain an injective composite \hm %
\[
B \longrightarrow C(S^1) \otimes (\OA{2})_{\infty}
  \longrightarrow D \otimes (\OA{2})_{\infty} \longrightarrow
  (D \otimes \OA{2})_{\infty} \stackrel{\cong}{\longrightarrow}
   (\OA{2})_{\infty},
\]
with a \ucp\  lifting given by
\[
B \longrightarrow C(S^1) \otimes \LI{\OA{2}}
  \longrightarrow D \otimes \LI{\OA{2}} \longrightarrow
  \LI{D \otimes \OA{2}} \longrightarrow \LI{\OA{2}}.
\]

The crossed product $B = C^* (\Z, B_0, \ta)$ is still exact, by
Proposition 7.1 (v) of \cite{Kr2}. Lemma 2.2
therefore provides
an injective unital \hm\  $\gm : B \to \OA{2}.$
Now $B$ contains as a subalgebra
$C^* (\Z, C_0 (\R) \otimes A)) \cong C(S^1) \otimes \Kt A,$
and this subalgebra in turn contains an isomorphic copy $A_0$
of $A.$ Let $p \in B$ be the identity of $A_0.$ Then
\[
\gm |_{A_0} : A_0 \to \gm (p) \OA{2} \gm (p)
\]
is a unital embedding of $A$ in $\gm (p) \OA{2} \gm (p) \cong \OA{2},$
as desired.
\end{pff}

\section{Tensor products with $\OA{2}$ and $\OI$}

In this section, we prove that if $A$ is a simple separable unital
nuclear \ca, then $\OT{A} \cong \OA{2},$ and that if, in addition,
$A$ is purely infinite, then $\OIA{A} \cong A.$
The key technical point is that if $A$ is \kfive, and if $\om \in \SCC,$
then the relative commutant of the image of $A$ in the ultrapower
$A_{\om},$ the algebra of bounded sequences in $A$ modulo those that
vanish at $\om,$ is again purely infinite and simple.
Once we have this simplicity result,
the rest of the proof that $\OIA{A} \cong A$ is done by
essentially the same methods as those of \cite{Rr3}. (We actually
prove, for future use elsewhere, a somewhat more general statement.
This disguises the similarity with \cite{Rr3} a little.)

We begin by establishing notation for ultrapowers.

\begin{ntn} 
We adopt the following notation regarding ultrapowers and associated
objects. Let $D$ be a \ca. Then (as in Notation 2.1)
$\LI{D}$ denotes
the \ca\  of bounded sequences with values in $D.$
For any $d = (d_1, d_2, \dots) \in \LI{D},$ the function
$n \mapsto \| d_n \|$ is bounded and thus defines a \ct\  function
on the
Stone-\v{C}ech compactification $\bt \N$ of $\N.$ Therefore for each
$\om \in \SCC,$ the limit $\lim_{n \to \om} \| d_n \|$ exists.
In particular, it makes sense to define $\lim_{n \to \om} d_n = 0$
to mean $\lim_{n \to \om} \| d_n \| = 0.$
We then define a closed ideal $c_{\om} (D)$ in $\LI{D}$ and the
corresponding quotient by
\[
c_{\om} (D) = \{d \in \LI{D} : \lim_{n \to \om} d_n = 0 \}   \andeqn
           D_{\om} = \LI{D} / c_{\om} (D).
\]
Denote the quotient map by $\pi_{\om}^{(D)} : \LI{D} \to D_{\om}.$
If $D$ is clear from the context, we sometimes write $\pi_{\om}.$
Note that $\| \pi_{\om}^{(D)} (d) \| = \limo{n} \| d_n \|.$

We will regard $D$ as a subalgebra of $D_{\om}$ via the diagonal
embedding of $D$ in $\LI{D}.$

If $A$ is a subalgebra of a
\ca\  $B,$ we denote by $A' \cap B$ the relative commutant of
$A$ in $B.$ In particular, with the identification above, we denote by
$D' \cap D_{\om}$ the relative commutant of $D$ in $D_{\om}.$
\end{ntn}

\begin{rmk} 
If $F$ is a \fd\  \ca, then $\pi_{\om}^{(F)}$ is an isomorphism.
More generally, if $F$ is \fd\  and $D$ is arbitrary, then
$\pi_{\om}^{(F \otimes D)}$ defines an isomorphism
$(F \otimes D)_{\om} \to F \otimes D_{\om}.$
\end{rmk}

The following lemma contains the main part of the proof that
$A' \cap A_{\om}$ is simple.

\begin{lem}    
Let $A$ be a \kfalg, and let $\om \in \SCC.$
Let $a,$ $b \in A' \cap A_{\om}$ be selfadjoint with
$\spec (b) \subset \spec (a).$ Then there is a nonunitary isometry
$s \in A' \cap A_{\om}$
such that $ss^*$ commutes with $a$ and $s^* a s = b.$
\end{lem}

\begin{pff}
Scaling both $a$ and $b$ by the same factor, we may assume that
$\| a \|,$ $\|b \| \leq \pi / 2.$ Let $v = \exp ( ia),$ and let
$X = \spec (v),$ which is a subset of $S^1$ intersected with the
right halfplane. Let $z\in C(X)$ be the standard generating
unitary, $z (\zt) = \zt.$
Then the assignments
\[
z \otimes 1 \mapsto v  \andeqn  1 \otimes x \mapsto
    \pi_{\om} (x, x, \dots)
\]
define a unital \hm\  $\ph: C(X) \otimes A \to A_{\om},$
and similarly the assignments
\[
z \otimes 1 \mapsto \exp (ib) \andeqn  1 \otimes x \mapsto
    \pi_{\om} (x, x, \dots)
\]
define a unital \hm\  $\ps: C(X) \otimes A \to A_{\om}.$

We prove that $\ph$ is injective. (Note that $\ps$ need not be
injective.) To see this, note that, since
$A$ is simple, $\ker (\ph)$ must have the
form $C_0 (U) \otimes A$ for some open subset $U \subset X.$
If $U \neq \emptyset,$ then there is a nonzero
\ct\  function $f \in C_0 (U).$ This gives
$0 = \ph (f \otimes 1) = f (v),$ contradicting the fact that
$f$ is a nonzero element of $C( \spec (v) ).$

Since $C(X) \otimes A$ is nuclear, there are, by Theorem 0.3,
\ucp\  maps $V,$ $W : C(X) \otimes A \to \LI{A}$ which lift
$\ph$ and $\ps.$ Then $V$ has the form
$V (x) = (V_1 (x), V_2 (x),  \dots)$ for \ucp\  maps
$V_m : C(X) \otimes A \to A,$ and similarly
$W (x) = (W_1 (x), W_2 (x),  \dots)$ with $W_m$ \uacp.
We next want to apply Lemma 1.10,
and for this we need
information on the injectivity of $V_m$ on \fd\  subspaces.

Choose a sequence
$u_1, u_2, \dots$ of unitaries in $A$ whose linear span is dense,
and let $E_n \subset C(X) \otimes A$ be given by
\[
E_n = \spn \{1, z \otimes 1, z^*\otimes 1, 1 \otimes u_1,
        1 \otimes u_1^*, \dots, 1 \otimes u_n, 1 \otimes u_n^*\}.
\]
Temporarily fix $n$ and $k.$ For
$x \in E_n \otimes M_k \subset C(X) \otimes A \otimes M_k,$
we have (using Remark 3.2
and the definition of $c_{\om} (A)$
in the first step and injectivity of $\ph \otimes \id_{M_k}$ in
the second)
\[
\limo{m} \| (V_m \otimes \id_{M_k}) (x)\| =
   \| (\ph \otimes \id_{M_k}) (x)\| = \| x \|.
\]
Since $E_n$ is finite dimensional, it follows that there is a
neighborhood $U$ of $\om$ in $\bt \N$ such that for all
$m \in U \cap \N,$ the map $V_m |_{E_n}$ is invertible; moreover,
$\limo{m} \| (V_m |_{E_n})^{-1} \otimes \id_{M_k}\| = 1.$
This holds for all $n$ and $k.$

Using Lemma 1.10,
choose positive integers $k(m)$ with
$k(0) \leq k(1) \leq \cdots$ and such that whenever
$\tilde{V},$ $\tilde{W} : E_n \to A$ are \ucp\  with $\tilde{V}$
injective and
$\|\tilde{V}^{-1} \otimes \id_{M_{k (m)}} \| \leq 1 + 1 / m,$
then there is a \ucp\  map $T: A \to A$ such that
$\|T \circ \tilde{V} - \tilde{W} \| < 2 / m.$
Choose a decreasing sequence of neighborhoods
$U_1 \supset U_2 \supset \cdots$ of $\om$ in $\bt \N$ such that
for all $m \in U_n \cap \N$ we have
\[
\| (V_m |_{E_n})^{-1} \otimes \id_{M_{k (m)}} \| \leq 1 + 1 / m.
\]
Replacing $U_n$ by $U_n - \{1, 2, \dots, n\},$ we may assume that
$\N \cap \bigcap_{n + 1}^{\infty} U_n = \emptyset.$
(Note that we can't have $\bigcap_{n + 1}^{\infty} U_n = \{\om\},$
since $\om$ doesn't have a countable neighborhood base. However, it is
certainly true that if $f : \N \to \C$ is a function such that
$\limi{n} \sup_{m \in U_n} | f(m) | = 0,$ then $\limo{n} f(n) = 0.$)
Lemma 1.10
now provides \ucp\  maps $T_m : A \to A$ such that
$\| T_m \circ V_m |_{E_n} - W_m |_{E_n} \| \leq 2 / m$
for $m \in (U_n - U_{n + 1}) \cap \N.$ By Proposition 1.7
(combined with the compactness of the closed unit ball of $E_n$), there
are nonunitary isometries $s_m \in A$ for
$m \in (U_n - U_{n + 1}) \cap \N$ such that
\[
\| s_m^* V_m (x) s_m - W_m (x) \| \leq 3 \| x \| / m
\]
for all $x \in E_n.$ Since
\[
E_n \subset E_{n + 1} \subset \cdots  \andeqn
  \N \cap U_n \cap U_{n + 1} \cap \cdots = \emptyset,
\]
this estimate in fact holds for all
$m \in U_n \cap \N$ and $x \in E_n.$

Define $s = \pi_{\om} (s_1, s_2, \cdots).$ Clearly $s$ is an isometry
in $A_{\om}.$ It is not unitary since
\[
\| 1 - s s^* \| =  \limo{m}  \| 1 - s_m s_m^* \| = 1.
\]
For $m \in U_n \cap \N$ we have
\[
\| s_m^* u_n s_m - u_n \| \leq
      \| V_m (1 \otimes u_n) - u_n \| +  \| W_m (1 \otimes u_n) - u_n \|
                       + 3 / m.
\]
Fix $n$ and let $m \to \om.$ The last term on the right certainly
converges to $0.$ The first two terms do so as well because,
recalling our identification of $A$ as a subalgebra of $A_{\om},$ we
have
$\pi_{\om} (V (1 \otimes u_n)) = \pi_{\om} (W (1 \otimes u_n)) = u_n.$
Thus $\limo{m} \| s_m^* u_n s_m - u_n \| = 0$ for each fixed $n.$
It follows that $s^* u_n s = u_n$ for all $n.$ Using $z \otimes 1$
in place of $1 \otimes u_n,$ we also obtain $s^* v s = \exp (i b).$
Since $u_n$ and $\exp (ib)$ are unitary, Lemma 1.11
implies that $s s^*$ commutes with $u_n$ and $v.$

Since $s s^*$ commutes with $u_n$ and $s^* u_n s = u_n,$ we have
\[
s u_n = s s^* u_n s = u_n s s^* s = u_n s.
\]
Since $u_1, u_2, \dots$ span a dense subspace of $A,$ this implies
that $s \in A' \cap A_{\om}.$

Since $s s^*$ commutes with $v,$ it also commutes with $v^*,$ and
it follows that $x \mapsto s^* x s = s^* [s s^* x s s^*] s$ is a
\hm\  from the unital \ca\  generated by $v$ to $A_{\om}.$ So
$s^* f(v) s = f(s^* v s)$ for every continuous function $f$ on
$\spec (v).$ Taking $f = - i \log,$ we obtain $s^* a s = b.$
This completes the proof.
\end{pff}

\begin{prp}   
Let $A$ be a \kfalg, and let $\om \in \SCC.$
Then $A' \cap A_{\om}$ is unital, simple, and purely infinite.
\end{prp}

\begin{pff}
Obviously $A' \cap A_{\om}$ is unital.
We show that every nonzero hereditary subalgebra $B$ of
$A' \cap A_{\om}$ contains a \pj\ $e \neq 1$ which is \mvn\  to $1.$
(This clearly implies that $A' \cap A_{\om}$ is simple and that
every nonzero hereditary subalgebra of $A' \cap A_{\om}$
contains an infinite \pj.)
So choose $c \in B$ selfadjoint with $1 \in \spec (c).$
Apply the previous lemma with $a = c$ and $b = 1$ to obtain a
nonunitary isometry
$s \in A' \cap A_{\om}$ such that $s^* c s = 1$ and the projection
$e = s s^*$ commutes with $c.$ By construction, $e$  is \mvn\  to $1.$
Furthermore, $ece = s (s^* c s) s^* = e,$ and from $ec = ce$ we get
$cec = (ece)^2 = e^2 = e.$ Therefore $e \in B.$
\end{pff}

We note that the proposition can be proved without knowing that the
isometry of Lemma 3.3 can be chosen to be nonunitary.
One still gets from the proof above that each
$A' \cap A_{\om}$ is simple and either purely
infinite or isomorphic to $\C.$ It follows from \cite{AP2} that
$A$ has a nontrivial central sequence, that is, a sequence
$a \in \LI{A}$ such that $\pi_{\infty} (a) \in A' \cap A_{\infty}$
but for which there is no sequence $z \in \LI{Z(A)}$ satisfying
$\limi{n} \|a_n - z_n\| = 0.$ From this it is possible to deduce that
there is at least one $\om_0 \in \SCC$ such that
$A' \cap A_{\om_0} \not\cong \C.$ One can then show
that $A' \cap A_{\infty}$ contains a unital copy of $\OI.$ 
For every other $\om \in \SCC,$
the image of this subalgebra in $A' \cap A_{\om}$ is
nontrivial, so also $A' \cap A_{\om} \not\cong \C.$

\begin{dfn}    
Let $A$ and $B$ be \sep\  unital \ca s. An
{\em asymptotically central inclusion} of $A$ in $B$ is a sequence
of unital injective \hm s $\ph_n : A \to B$ such that
$\| \ph_n (a) b - b \ph_n (a) \| \to 0$ for all $a \in A$ and
$b \in B.$
\end{dfn}

The rest of the proof that $\OT{A} \cong \OA{2}$ was inspired by a
talk by Mikael R\o rdam.
The following lemma is the main remaining part.

\begin{lem}       
Let $A$ be a simple separable unital nuclear \ca\  which has an
asymptotically central inclusion of $\OA{2}.$ Then $A \cong \OA{2}.$
\end{lem}

\begin{pff}
We have a unital \hm\  $\ph : \OA{2} \to A$ by assumption, and
a unital \hm\  $\ps : A \to \OA{2}$ by Theorem 2.8. Furthermore,
$\ps \circ \ph$ is \ayue\  to $\id_{\OA{2}}$ by Proposition 0.7.
In the rest of the proof, we show that any two unital endomorphisms
of $A$ are \ayue. In particular, we then have
$\ph \circ \ps$ \ayue\  to $\id_A,$ so that $A \cong \OA{2}$ by
Lemma 0.6.

First observe that $A$ is purely infinite. Indeed,
clearly $A$ is infinite. Moreover, $\OA{2}$ is approximately
divisible in the sense of \cite{BKR}, by Proposition 7.7 of \cite{Rr1}.
Therefore so is $A.$ So $A$ is purely infinite by
Theorem 1.4 (a) of \cite{BKR}.

Now let $\gm : A \to A$ be a unital endomorphism; we show that
$\gm$ is \ayue\  to $\id_A.$ Certainly $\gm$ is nuclear, unital, and
completely positive, so by Proposition 1.7
there are
isometries $v_n \in A$ such that $\limi{n} v_n^* a v_n = \gm (a)$
for all $a \in A.$ Choose any $\om \in \SCC,$ and
set $v = \pi^{(A)}_{\om} (v_1, v_2, \dots ),$
which is an isometry in $A_{\om}.$ Regarding $A$ as a subalgebra
of $A_{\om}$ as in Notation 3.1,
we have
$v^* a v = \gm (a)$ for all $a \in A.$ It now follows from
Lemma 1.11
that $v v^*$ commutes with every unitary in $A,$
whence $v v^* \in A' \cap A_{\om}.$

Let $F_1 \subset F_2 \subset \cdots$ be finite selfadjoint subsets of
$A$ whose union is dense in $A,$ and such that $v_n v_n^* \in F_n.$
Let $s_1$ and $s_2$ be the standard generating isometries of $\OA{2}.$
The existence of an asymptotically central inclusion of $\OA{2}$ in $A$
provides unital \hm s $\sm_n : \OA{2} \to A$ such that
$\| \sm_n (x) a - a \sm_n (x) \| < \frac{1}{n}$ for $a \in F_n$
and $x$ in the generating set $\{s_1, s_1^*, s_2, s_2^*\}$ of $\OA{2}.$
The definition $\sm (x) = \pi^{(A)}_{\om} (\sm_1 (x), \sm_2 (x), \dots)$
yields a unital \hm\  $\sm : \OA{2} \to A' \cap A_{\om}$ whose range
also commutes with $v v^*.$ We then calculate in
$K_0 (A' \cap A_{\om}):$
\[
[ v v^*] = [ \sm (s_1) v v^* \sm (s_1)^* + \sm (s_2) v v^* \sm (s_2)^*]
   = 2 [ v v^*],
\]
so that $[ v v^*] = 0$ in $K_0 (A' \cap A_{\om}).$ Similarly
$[1] = 0$ in $K_0 (A' \cap A_{\om}).$ Since $A' \cap A_{\om}$ is
purely infinite and simple, it follows from Theorem 1.4 and
Proposition 1.5 of \cite{Cu1}
that there is $w \in A' \cap A_{\om}$ such that $w^* w = 1$ and
$w w^* = v v^*.$ Then $u = w^* v$ is a unitary in $A_{\om};$
moreover, for $a \in A$ we have
\[
u^* a u = v^* w a w^* v = v^* a v = \gm (a).
\]

To show that $\gm$ is \ayue\  to $\id_A,$ let $F \subset A$ be
finite, with $\|a\| \leq 1$ for $a \in F,$ and let and $\ep > 0.$
Choose $(b_1, b_2, \dots) \in \LI{A}$ such that
$\pi^{(A)}_{\om} (b_1, b_2, \dots) = u.$
\Wolog, $\|b_n\| \leq 1$ for all $n.$
There is a neighborhood $U$ of $\om$ in $\bt \N$ such that, for every
$n \in U \cap \N,$ the unitary
$u_n = b_n (b_n^* b_n)^{-1/2}$ satisfies
$\|u_n - b_n\| \leq \frac{\ep}{3},$ and also
$\| b_n^* a b_n - \gm (a) \| < \frac{\ep}{3}$ for all $a \in F.$
Choose one such $n;$
then $\| u_n^* a u_n - \gm (a) \| < \ep$ for all $a \in F.$
\end{pff}

\begin{thm}      
Let $A$ be a simple separable unital nuclear \ca.
Then $\OT{A} \cong \OA{2}.$
\end{thm}

\begin{pff}
Let $B = \bigotimes_1^{\infty} \OA{2},$ which we think of as
${\dirlim} \bigotimes_1^n \OA{2}.$ Obviously there is an
asymptotically central inclusion of $\OA{2}$ in $B.$ So there is also
an asymptotically central inclusion of $\OA{2}$ in $B \otimes A.$
The previous lemma therefore implies that $B \cong \OA{2}$ and
$B \otimes A \cong \OA{2}.$ So $\OT{A} \cong \OA{2}.$
\end{pff}

\begin{rmk}      
It is actually possible to get this far (except for Theorem
1.15, which we haven't used yet) with only a unital (necessarily
injective) \hm\  $\OA{2} \otimes \OA{2} \to \OA{2}.$ Then one
would get $\OA{2} \otimes \OA{2} \cong \OA{2}$ as a corollary
to the previous theorem. However, there doesn't seem to be a
way to get such a \hm\  which is simpler than going through much
of R\o rdam's proof of the isomorphism.
\end{rmk}

We now turn to the proof that $\OIA{A} \cong A.$ For use elsewhere, we
prove a somewhat more general statement, in which $\OI$ is replaced
by a subalgebra of $B \subset A' \cap A_{\om}.$ We will eventually
take $B = \bigotimes_1^{\infty} \OI,$
but for now we merely assume that it is separable,
that it contains the unit of $A' \cap A_{\om},$  and that the two
obvious maps from $B$ to $B \tmin B$ are \ayue.

\begin{lem}      
(Compare with Propositions 2.7 and 2.8 of \cite{EfR}.)
Let $B$ be a separable unital \ca. Suppose that the two maps
$\af,$ $\bt : B \to B \tmin B,$ given by
\[
\af (b) = b \otimes 1  \andeqn  \bt (b) = 1 \otimes b,
\]
are \ayue. Then $B$ is simple and nuclear.
\end{lem}

\begin{pff}
We prove nuclearity first. This argument is taken from \cite{EfR}.
By hypothesis there is a sequence
$u_1, u_2, \dots$ of unitaries in $B \tmin B$ such that
$\limi{n} u_n ( b \otimes 1) u_n^* = 1 \otimes b$ for all $b \in B.$
Choose $c_n$ in the algebraic tensor product $B \talg B$ with
$\| c_n \| \leq 1$ and $\limi{n} \| c_n - u_n \| = 0.$ Then also
$\limi{n} c_n ( b \otimes 1) c_n^* = 1 \otimes b$ for all $b \in B.$
Choose any state $\om$ on $B,$ and define $T_n : B \to B$ by
\[
T_n (b) = (\om \otimes \id_B) (c_n ( b \otimes 1) c_n^*).
\]
Note that $\om \otimes \id_B : B \tmin B \to B$ is well defined, unital,
and completely positive. (See, for example, Proposition IV.4.23 (i) of
\cite{Tk}.) Since $\| c_n \| \leq 1,$ the maps $T_n$ are thus
completely positive contractions.

For fixed $n,$ write $c_n = \sum_{j = 1}^m x_j \otimes y_j$ with
$x_j,$ $y_j \in B.$ Then
\[
T_n (b) =  \sum_{j, \, k = 1}^m \om (x_j b x_k^*) y_j y_k^*,
\]
so that $T_n$ has finite rank (at most $m^2$). We further have
\[
\| T_n (b) - b \| =
   \| (\om \otimes \id_B) (c_n ( b \otimes 1) c_n^*)
                                 - (\om \otimes \id_B) (1 \otimes b) \|
  \leq \| c_n ( b \otimes 1) c_n^* - 1 \otimes b \|,
\]
which converges to $0$ as $n \to \infty.$ Thus, we have shown that
$\id_B$ is a pointwise norm limit of completely positive contractions.
So $B$ is nuclear. (It is not necessary to require that the
approximating maps be unital. See the comment after Theorem 0.2.)

Now we prove simplicity. Suppose $B$ has a nontrivial ideal $J.$
As in the proof of Proposition 2.7 of \cite{EfR}, if $b \in J$
then $b \otimes 1 \in J \tmin B$ but $b \otimes 1 \not\in B \tmin J.$
However, with $u_n$ as in the first paragraph of the proof, we have
\[
u_n^* (1 \otimes b) u_n \in B \tmin J \andeqn
   \limi{n} u_n^* (1 \otimes b) u_n = b \otimes 1,
\]
which implies $b \otimes 1 \in B \tmin J.$ This contradiction shows
that $B$ is simple.
\end{pff}

We can now write simply $\otimes$ instead of $\tmin$ for tensor
products involving $B.$

\begin{lem}      
Let $A,$ $B,$ and $C$ be separable unital \ca s, and let $\om \in \SCC.$
Let $S (b) = (S_1 (b), S_2 (b), \dots)$ and
$T (c) = (T_1 (c), T_2 (c), \dots)$ define
\ucp\  maps from $B$ and $C$ respectively to $\LI{A}$ such that
$\pi_{\om} \circ S$ and $\pi_{\om} \circ T$ are unital \hm s whose
images lie in $A' \cap A_{\om}.$
For any finite subsets $F \subset A,$ $G \subset B,$ and $H \subset C,$
and any $k$ and $\ep > 0,$ there is a neighborhood $U$ of $\om$ such
that for every $n \in U \cap \N,$ we have
\[
\| T_n (c) S_k (b) - S_k (b) T_n (c) \| < \ep, \,\,\,\,\,\,
\| T_n (c) a - a T_n (c) \| < \ep, \andeqn
\| T_n (c_1 c_2) - T_n (c_1) T_n (c_2)  \| < \ep,
\]
for all $a \in F,$ $b \in G,$ and $c,$ $c_1,$ $c_2 \in H.$
\end{lem}

\begin{pff}
This is immediate.
\end{pff}

The following lemma is closely related to Lemmas 2 and 3 of \cite{Rr3}.

\begin{lem}       
Let $A$ be a separable unital \ca, let $\om \in \SCC,$ and let
$B \subset A' \cap A_{\om}$ be a unital subalgebra which satisfies
the hypotheses of Lemma 3.9.
Then there is a unital
\hm\  $\ph : B \otimes A \to A$ such that the map
$a \mapsto \ph (1 \otimes a)$ is \ayue\  to $\id_A.$
\end{lem}

\begin{pff}
Since $B$ is nuclear (by Lemma 3.9),
its inclusion in
$A_{\om}$ lifts to a \ucp\  map $Q : B \to \LI{A}.$ We write this
map as $Q (b) = (Q_1 (b), Q_2 (b), \dots)$ for \ucp\  maps
$Q_n : B \to A.$

Choose finite selfadjoint subsets
\[
F_1 \subset F_2 \subset \cdots \subset A \andeqn
      G_1 \subset G_2 \subset \cdots \subset B
\]
whose unions are dense in $A$ and $B.$
Using the hypothesis on $B,$ choose unitaries $u_k \in B \otimes B$
such that $\| u_k (b  \otimes  1) u_k^* - 1 \otimes b \| < 2^{-k}$
for $b \in G_k.$ Choose $b_k  \in B \otimes B$ in the algebraic
tensor product and so close to $u_k$ that the unitary
$z = b_k (b_k^* b_k)^{-1/2}$ satisfies
$\| z (b  \otimes  1) z^* - 1 \otimes b \| < 2 \cdot 2^{-k}$
for $b \in G_k.$ Write
$b_k = \sum_{j = 1}^l d_j^{(1)} \otimes d_j^{(2)}$ (suppressing the
dependence on $k$ on the right).
There are then a finite set  $G_k' \subset B$
(containing $G_k$ and all $d_j^{(i)}$) and $\ep_k > 0$ such
that whenever $S,$ $T : B \to A$ are
\ucp\  maps such that
\[
\| S (bc) - S(b) S(c) \| < \ep_k, \,\,\,\,\,\,
     \| T(bc) - T (b) T (c) \| < \ep_k,
  \,\,\,\,\,\, \| S (b) T (c) - T (c) S (b) \| < \ep_k,
\]
\[
\| S (b) a - a S (b) \| < \ep_k, \andeqn \| T (b) a - a T (b) \| < \ep_k
\]
for $a \in F_k$ and $b,$ $c \in G_k',$ then there is a unitary
$v \in A$ (close to $\sum_{j = 1}^l S (d_j^{(1)}) T (d_j^{(2)})$)
such that
\[
\| v S (b) v^* - T (b) \| < 3 \cdot 2^{-k} \andeqn
  \| v a v^* - a \| <  2^{-k}
\]
for $a \in F_k$ and $b \in G_k.$
Constructing the $G_k'$ and $\ep_k$ in order, we may assume that
$G_k' \subset G_{k + 1}'$ and $\ep_k > \ep_{k + 1}$ for all $k,$
and that $\ep_k \to 0.$

Now use the previous lemma to
choose inductively $n (1) < n (2) < \cdots$ such that
\[
\| Q_{n (1)} (bc) - Q_{n (1)}(b) Q_{n (1)}(c) \| < \ep_1 \andeqn
\| Q_{n (1)} (b) a - a Q_{n (1)} (b) \| < \ep_1
\]
for $a \in F_1$ and $b,$ $c \in G_1',$ and
\[
\| Q_{n (k + 1)} (bc) - Q_{n (k + 1)}(b) Q_{n (k + 1)}(c) \|
                               < \ep_{k + 1}, \,\,\,\,\,\,
\| Q_{n (k + 1)} (b) a - a Q_{n (k + 1)} (b) \| < \ep_{k + 1},
\]
and
\[
\| Q_{n (k + 1)} (b) Q_{n (k)} (c) - Q_{n (k)} (c) Q_{n (k + 1)} (b) \|
                               < \ep_k
\]
for $a \in F_k$ and $b,$ $c \in G_{k + 1}'.$
The previous paragraph gives a unitary $v_k \in A$ such that
\[
\| v_k Q_{n (k + 1)} (b) v_k^* - Q_{n (k)} (b) \| < 3 \cdot 2^{-k}
     \andeqn \| v_k a v_k^* - a \| <  2^{-k}
\]
for $a \in F_k$ and $b \in G_k.$

Define a unitary
$w_n \in A$ by $w_k = v_1 v_2 \cdots v_{k - 1}.$ Since
$\sum_{k = 1}^{\infty} 2^{-k} < \infty$ and
$\bigcup_{n = 1}^{\infty} F_k$ is dense in $A,$
one checks that $\af (a) = \limi{k} w_k a w_k^*$ exists for all
$a \in A.$ Clearly $\af$ is a unital \hm\  from $A$ to $A$ which is
\ayue\  to $\id_A.$ Also,
\[
\|w_{k + 1}  Q_{n (k + 1)} (b) w_{k + 1}^* - w_k Q_{n (k)} (b) w_k^*\|
          < 3 \cdot 2^{-k}
\]
for $b \in G_k,$ so $\bt (b) = \limi{n} w_k Q_{n (k)} (b) w_k^*$
exists for $b \in B.$ Clearly $\bt (bc) = \bt (b) \bt (c)$ for
$b,$ $c \in \bigcup_{n = 1}^{\infty} G_k',$ so $\bt$ is a unital
\hm\   from $B$ to $A.$
Moreover, for $a \in F_{k}$ and $b \in G_{k},$ we have
\[
\| (w_k a w_k^*) (w_k Q_{n (k)} (b) w_k^*) -
                  (w_k Q_{n (k)} (b) w_k^*) (w_k a w_k^*) \|
   = \| a Q_{n (k)} (b) - Q_{n (k)} (b) a \| < \ep_k,
\]
and $\limi{k} \ep_k = 0.$ It follows that the ranges of $\af$
and $\bt$ commute. Since $B \tmax A = B \otimes A,$ the desired
\hm\  $\ph : B \otimes A \to A$ can be defined by the formula
$\ph (b \otimes a) = \af (a) \bt (b).$
\end{pff}

\begin{prp}   
Under the hypotheses of Lemma 3.11,
we have $B \otimes A \cong A.$
\end{prp}

\begin{pff}
Since $B$ is nuclear, we write $\otimes$ rather than $\tmin$ for
tensor products involving $B.$

We follow the proof of the theorem following Lemma 3 in \cite{Rr3}.
Let $\bt : B \otimes A \to A$ be the \hm\  of Lemma 3.11,
and
let $\af : A \to B \otimes A$ be the \hm\  $\af (a) = 1 \otimes a.$
We know from Lemma 3.11
that $\bt \circ \af$ is \ayue\  to $\id_A.$ We show that
$\af \circ \bt$ is \ayue\  to $\id_{B \otimes A}.$ By Lemma 0.6,
this will imply that $B \otimes A \cong A.$

By the definition of \aue, there is a sequence $w_1, w_2, \dots$ of
unitaries in $A$ such that
\[
\limi{n} \| w_n \bt ( \af (a)) w_n^* - a\| = 0
\]
for all $a \in A.$
The hypothesis on $B$ provides a sequence $v_1, v_2, \dots$ of
unitaries in $B \otimes B$ such that
\[
\limi{n} \| v_n (x \otimes 1) v_n^* - 1 \otimes x \| = 0
\]
for all $x \in B.$ 
Note that for $x, \, y \in B,$ the elements $\af ( \bt (x \otimes 1))$
and $y \otimes 1$ of $B \otimes A$ commute. Therefore there is a
\hm\  $\sm : B \otimes B \to B \otimes A$ satisfying
\[
\sm (x \otimes 1) = \af ( \bt (x \otimes 1))  \andeqn
\sm (1 \otimes y) = y \otimes 1.
\]
Define $u_n = (1 \otimes w_n) \sm (v_n).$ Using the fact that
$\sm (v_n)$ commutes with $\af ( \bt (1 \otimes a))$ for
$a \in A$ and $1 \otimes w_n$ commutes with
$x \otimes 1 = \limi{n} \sm (v_n) \af ( \bt (x \otimes 1)) \sm (v_n)^*$
for $x \in B,$ we obtain as in \cite{Rr3}
$\limi{n} u_n \af ( \bt (c)) u_n^* = c$ for all $c \in B \otimes A.$
\end{pff}

\begin{lem}
Let $A$ be a separable unital \ca, let $\om \in \SCC,$ and let
$B \subset A' \cap A_{\om}$ be a separable
nuclear unital subalgebra. Let $C$ be
a separable nuclear unital \ca, and let $\ph : C \to A' \cap A_{\om}$ be
a unital \hm. Then there is a unital \hm\  $\ps : C \to A' \cap A_{\om}$
whose image lies in the relative commutant of $B.$
\end{lem}

\begin{pff}
Since $B$ and $C$ are nuclear, there are \ucp\  maps
$S : B \to \LI{A}$ and $T : C \to \LI{A}$
such that $\pi_{\om} \circ S$ is the inclusion of $B$ and
$\pi_{\om} \circ T = \ph.$ Write
$S (b) = (S_1 (b), S_2 (b), \dots)$ and
$T (c) = (T_1 (c), T_2 (c), \dots).$
Choose finite subsets
\[
F_1 \subset F_2 \subset \cdots \subset A, \,\,\,\,\,\,
      G_1 \subset G_2 \subset \cdots \subset B, \andeqn
      H_1 \subset H_2 \subset \cdots \subset C
\]
whose unions are dense.
Using Lemma 3.10, choose $n (1) \, n (2), \dots \in \N$ such that
\[
\| T_{n (k)} (c) S_k (b) - S_k (b) T_{n (k)} (c) \| < \frac{1}{k},
        \,\,\,\,\,\,
\| T_{n (k)} (c) a - a T_{n (k)} (c) \| < \frac{1}{k},
\]
and
\[
\| T_{n (k)} (c_1 c_2) - T_{n (k)} (c_1) T_{n (k)} (c_2)  \| <
          \frac{1}{k}
\]
for all $a \in F_k,$ $b \in G_k,$ and $c,$ $c_1,$ $c_2 \in H_k.$
Then define $\ps (c) = \pi_{\om} (T_{n (1)} (c), T_{n (2)} (c), \dots).$
\end{pff}

\begin{thm}   
Let $A$ be a \kfalg . Then $\OIA{A} \cong A.$
\end{thm}

\begin{pff}
Let $B = \bigotimes_1^{\infty} \OI,$ which we think of as
the direct limit over $n$ of $B_n = \bigotimes_1^n \OI,$ with
maps $b \mapsto b \otimes 1.$ We apply Proposition 3.12 with this
$B.$ Clearly $B$ is separable, unital, and nuclear.
We need to check two more conditions.

First, we must embed $B$ as a unital subalgebra of $A' \cap A_{\om}.$
Now $A' \cap A_{\om}$ is
purely infinite simple by Proposition 3.4, so certainly contains
a unital copy of $\OI.$ Using the previous lemma and induction,
we obtain unital \hm s $\ph_n : B_n \to A' \cap A_{\om}$ such that
$\ph_{n + 1} (b \otimes 1) = \ph_n (b)$ for $b \in B_n.$ Taking
direct limits gives a unital \hm\  $\ph : B \to A' \cap A_{\om}.$
This \hm\  is injective because $B$ is simple.

The second condition to check is that the two maps
$\af (b) = b \otimes 1$ and $\bt (b) = 1 \otimes b,$ from $B$ to
$B \otimes B,$ are \ayue.

For $F \subset \OI,$ let $F^{(n)} \subset B_n$ be the set of all
$1 \otimes \cdots \otimes 1 \otimes b \otimes 1 \otimes \cdots
                                     \otimes 1,$
with $b \in F,$ and where $b$ is in the $k$-th tensor factor for some
$k$ with $1 \leq k \leq n.$ Also let $\ps_n : B_n \to B$ be the
inclusion. It suffices to show that for each finite $F \subset \OI,$
each $n,$ and each $\ep > 0,$ there exists a unitary $v \in B \otimes B$
such that
$\| v (\ps_n (x) \otimes 1 ) v^* - 1 \otimes \ps_n (x) \| < \ep$
for all $x \in F^{(n)}.$

Now clearly $B$ has an asymptotically central inclusion of $\OI.$
It follows that $B$ is approximately
divisible in the sense of \cite{BKR}. Since $B$ is simple and infinite,
it is purely infinite by Theorem 1.4 (a) of \cite{BKR}.
So Theorem 3.3 of \cite{LP} implies that the maps from $\OI$ to
$B \otimes B,$ given by $x \mapsto \ps_1 (x) \otimes 1$ and
$x \mapsto 1 \otimes \ps_1 (x),$ are \ayue. Approximating unitaries
in $B$ by ones in the terms of the direct system, we find that
for $F$ and $\ep$ as above there is $N$ and $u \in B_N \otimes B_N$ such
that, with $1_m$ denoting the identity of $\bigotimes_1^m \OI,$
we have
\[
\| u [ (b \otimes 1_{N - 1}) \otimes 1_N ] u^*
                    - 1_N \otimes  (b \otimes 1_{N - 1}) \| < \ep
\]
for all $b \in F.$
Taking the tensor product of this with itself $n$ times, we find
that
$v_0 = u \otimes \cdots \otimes u \in \bigotimes_1^n (B_N \otimes B_N)$
satisfies
\[
\| v_0 (1 \otimes \cdots \otimes 1 \otimes
           [ (b \otimes 1_{N - 1}) \otimes 1_N ]
               \otimes 1 \otimes \cdots \otimes 1) v_0^*
  - 1 \otimes \cdots \otimes 1 \otimes
           [ 1_N \otimes  (b \otimes 1_{N - 1})]
                \otimes 1 \otimes \cdots \otimes 1 \| < \ep
\]
for all $b \in F$ and with $(b \otimes 1_{N - 1}) \otimes 1_N$ in any
of the $n$ factors $B_N \otimes B_N.$ Rearranging this using
associativity of the tensor product, we obtain a unitary
$v \in B_{Nn} \otimes B_{Nn}$ such that
\[
\| v [ (x \otimes 1_{(N - 1)n}) \otimes 1_{Nn}] v^*
                    - 1_{Nn} \otimes (x \otimes 1_{(N - 1)n}) \| < \ep
\]
for all $x \in F^{(n)}.$ Embedding $B_{Nn}$ in $B$ completes the
proof that $\af$ and $\bt$ are \ayue.

We now apply Proposition 3.12, obtaining $B \otimes A \cong A.$
Taking in particular $A = \OI,$ we obtain $\OI \cong B \otimes \OI.$
But
clearly $B \otimes \OI \cong B.$ So we can replace $B$ by $\OI,$ getting
$\OIA{A}  \cong A.$
\end{pff}

\section{Exact continuous fields}

This section consists of various preparatory results on (exact)
continuous fields of \ca s which will be used in the next section to
obtain continuous embeddings of continuous fields in $\OA{2}.$
We start with several general results, on tensor products, pullbacks,
and representations, and then go on to apply the results of
Section 1 to continuous fields whose section algebras are exact.
We obtain a discrete version of continuous embedding: Fibers over
nearby points have embeddings in $\OA{2}$ which are close in a
suitable sense. In the next section, we show how to make continuously
varying choices of the embeddings.

Continuous fields of \ca s are taken to be as defined in Chapter 10
of \cite{Dx}. For notation,
if $A$ is a continuous field over $X,$ then we
denote by $A (x)$ the fiber over $x \in X$ and
by $\Gm (A)$ the set of all continuous sections of $A.$ We briefly
recall the axioms for $\Gm (A)$ (\cite{Dx}, 10.1.2 and 10.3.1):

(1) Each $A (x)$ is a \ca.

(2) The section space $\Gm (A)$ is a subspace of the set-theoretic
product
$\prod_{x \in X} A (x)$ which is closed under addition, scalar
multiplication, multiplication, and adjoint.

(3) The set
$\{ a (x) : a \in \Gm (A) \}$ is dense in $A (x)$ for all $x.$
(Proposition 10.1.10 of \cite{Dx} shows that, in the presence of the
other axioms, this is equivalent to requiring that
$\{ a (x) : a \in \Gm (A) \} = A (x)$ for all $x.$)

(4) For every $a\in \Gm (A)$ the function $x \mapsto \| a (x) \|$ is
continuous.

(5) The section space
$\Gm (A)$ is closed under local uniform approximation. That is,
if $a$ is a section, and if for every $x_0 \in X$ and $\ep > 0$ there
exists
a neighborhood $U$ of $x_0$ and a continuous section $b$ such that
$\| a (x) - b (x) \| < \ep$ for $x \in U,$ then $a$ is continuous.

It follows from Axiom (5) that the pointwise scalar product of
a continuous function on $X$ and a continuous section of $A$ is
again  a continuous section of $A.$ (See Proposition 10.1.9 of
\cite{Dx}.)

More general bundles, in which (4) is weakened to merely require that
the norms of sections be upper semicontinuous,
are in some ways more natural. (See \cite{HK} for a general discussion,
mostly in the context of Banach spaces.)
However, the results of this and the
next section have no possibility of being true for them.

We denote by $\ev_x$ the evaluation map
$a \mapsto a (x)$ from $\Gm (A)$ to $A (x).$ We say that $A$ is
{\em unital} if each $A (x)$ is unital and the section
$x \mapsto 1_{A (x)}$ is continuous.
Note that it is possible to have every $A (x)$ unital but the section
$x \mapsto 1_{A (x)}$ discontinuous.

We start with several general results on continuous fields.

\begin{prp}   
Let $X$ be a \cpt\  space, let $A$ be a continuous field of
\ca s over $X,$ and let $B$ be a nuclear \ca. Then there
is a continuous field $A \otimes B$ of \ca s over $X,$
such that $(A \otimes B) (x) = A (x) \otimes B$ and
$\Gm (A \otimes B) = \Gm (A) \otimes B.$
\end{prp}

\begin{pff}
This is (i) implies (v) of Theorem 3.2 of \cite{KW}.
\end{pff}

The cases we need are $B = M_n$ and $B = K,$ which can be handled a
little more directly.

\begin{lem}   
Let $X$ be a topological space, and let $A$ be a continuous field of
\ca s over $X.$ Then
there is a continuous field $A\unit$ of \ca s over $X$ such that
$A\unit (x)$
is the unitization $A (x)\unit$ for every $x \in X,$ and such
that the continuous sections of $A\unit$ are the sections of the form
$a (x) = a_0 (x) + \ld (x) \cdot 1_{A (x)\unit}$ for $a_0$ a
continuous section of $A$ and $\ld : X \to \C$ continuous.
\end{lem}

\begin{pff}
We define $A\unit (x) = A (x)\unit$ for each $x.$ (Recall that we add
a new identity to $A (x)$ even if it already has one.) We take the
continuous sections of $A\unit$ to be as in the statement. Axioms (1),
(2), and (3) are obvious. Axiom (5) is easily checked using the fact
that $A (x)\unit = A (x) \oplus \C$ as a Banach space. This leaves
axiom (4).

Considering $a (x)^* a(x),$ we see that it suffices to consider sections
$a (x) = a_0 (x) + \ld (x) \cdot 1_{A (x)\unit}$ with $a (x) \geq 0$
for all $x.$ In this case, the corresponding
function $\ld$ is nonnegative and $a_0$ is selfadjoint.
Let $f : \C \to \R$ be $f (t) = \max( \re (t), 0).$
By Proposition 10.3.3 of
\cite{Dx}, the section $x \mapsto f (a_0 (x))$ is continuous.
Now 
\[
\| a (x) \| = \ld (x) + \sup \spec (a_0 (x))  =
              \ld (x) + \| f (a_0 (x)) \|,
\]
which is continuous.
\end{pff}

\begin{lem}   
Let $f : X \to Y$ be a continuous map of topological spaces, and let
$A$ be a continuous field over $Y.$ Then there is a
continuous field $f^* (A)$ over $X$ such that $f^* (A) (x) = A(f(x))$
for all $x \in X,$ and such that the continuous sections of
$f^* (A)$ are the locally uniform limits (in the sense implicit in
Axiom (5)) of sections of the form $x \mapsto a(f(x))$ for
$a \in \Gm (A).$
\end{lem}

\begin{pff}
It is immediate from the corresponding axioms for $A$
that $f^* (A)$ satisfies axioms (1) and (3). Axiom (5)
for $f^* (A)$
follows from the construction of $\Gm (f^* (A))$
as the set of locally uniform limits
of a set of sections. Axioms (2) and (4) hold for the set of
sections of the form $x \mapsto a(f(x))$ for
$a \in \Gm (A),$ and are preserved under passage to
locally uniform limits, so also hold for $\Gm (f^* (A)).$
\end{pff}

\begin{lem}   
Let $X$ be a locally compact Hausdorff space, let $U \subset X$ be open,
and let $A_0$ be a unital continuous field of \ca s over $U.$ Then
there is a unital continuous field $A$ of \ca s over $X$ such that
$A (x) = A_0 (x)$ for $x \in U$ and $A (x) = \C$ for $x \not\in U,$
and such that a section $a$ is continuous if and only if
there is a continuous function $\ld : X \to \C$ and a continuous section
$a_0$ of $A_0$ for which $x \mapsto \| a_0 (x) \|$ vanishes at
infinity on $U,$
such that $a (x) = \ld (x) \cdot 1_{A (x)}$ for $x \not\in U$ and
$a (x) = \ld (x) \cdot 1_{A (x)} + a_0 (x)$ for $x \in U.$
\end{lem}

\begin{pff}
One checks directly that the given set of sections satisfies the
definition of a unital continuous field of \ca s.
\end{pff}

The following definition and lemma do not work very well for more
general bundles (for which the norm of a section is only required to
be upper semicontinuous).

\begin{dfn}   
Let $X$ be a topological space, and let $A$ be a continuous field of
\ca s over $X.$ A {\em representation} of
$A$ in a \ca\  $D$ is a family $\ph = ( \ph_x )_{x \in X}$ of \hm s
$\ph_x : A(x) \to D$ which is continuous in the following sense:
for every continuous section $a$ of $A,$ the function
$x \mapsto \ph_x (a(x))$ is continuous from $X$ to $D.$
The representation is called {\em injective} if every $\ph_x$ is
injective.
\end{dfn}

In this terminology, a continuous field is Hilbert continuous
(Definition 3.3 of \cite{Rf}) if it has an injective representation
in some $L (H).$

\begin{lem}   
Let $X$ be a topological space, let $A$ be a continuous field of
\ca s over $X,$ and let $\ph = ( \ph_x )_{x \in X}$ be a
representation of $A$ in some \ca\  $D.$ Suppose $\ph_x$ is
injective for every $x$ is some dense subset $S$ of $X.$ Then
$\ph_x$ is injective for every $x \in X.$
\end{lem}

\begin{pff}
Let $x \in X,$ and suppose $\ph_x$ is not injective. Choose an element
$a \in \ker (\ph_x)$ with $\|a\| = 1.$
By Proposition 10.1.10 of \cite{Dx}, there
is a continuous section $b$ of $A$ such that $b (x) = a.$
By continuity of $x \mapsto \| b (x) \|$ and $x \mapsto \ph_x (b (x)),$
there is a neighborhood $U$ of $x$ such that $\| b (y) \| > 3/4$
and $\| \ph_y (b (y)) \| < 1/4$ for $y \in U.$ Choose $y \in U \cap S.$
Then $b (y)$ is an element of $A (y)$ such that
$\| \ph_y (b (y)) \| < \| b (y) \|,$ so $\ph_y$ is not injective. This
contradiction proves the lemma.
\end{pff}

We will work with continuous fields
satisfying the exactness conditions given
in the following theorem, essentially in \cite{KW}.

\begin{thm}   
Let $X$ be a compact metric space, and let
$A$ be a continuous field of \ca s over $X,$ with
$\Gm (A)$ separable. Then the following
conditions are equivalent:

(1) The section algebra $\Gm (A)$ is an exact \ca.

(2) Every fiber $A (x)$ is an exact \ca, and the identity maps
of the fibers
\[
\id_{A (x)} : A (x) \to \Gm (A) / \{a \in \Gm (A) : \ev_x (a) = 0 \}
 = A (x)
\]
are {\em{locally liftable,}} that is,
for every \fd\  operator system $E  \subset A (x)$ the inclusion
of $E$ in $A (x)$
has a \ucp\  lifting to a map from $E$ to $\Gm (A).$

(3) Every fiber $A (x)$ is an exact \ca, and for every \ca\  $B,$ the
tensor product bundle $A \tmin B$ (as described in the introduction
to \cite{KW}) is a continuous field.

(4) Every fiber $A (x)$ is an exact \ca, and for a separable infinite
dimensional Hilbert space $H,$ the
tensor product bundle $A \tmin L (H)$ is a continuous field.
\end{thm}

\begin{pff}
Note that if $\Gm (A)$ is exact, then the fibers $A (x),$ being
quotients of $\Gm (A),$ are exact by Proposition 7.1 (ii) of
\cite{Kr2}.
So the equivalence of the first three conditions (and some others)
is Theorem 4.6 of \cite{KW}.
That (3) implies (4) is immediate.
The proof that (4) implies (2) is derived from Theorem 3.2 of \cite{EH}
in the same way that Lemma 2.4 is derived from Theorem 3.4 of \cite{EH}.
\end{pff}

\begin{lem}   
Let $X$ be a compact metric space, and let
$A$ be a continuous field of \ca s over $X,$ with
$\Gm (A)$ separable. Let $A\unit$ be the unitization (as in
Lemma 4.2).
Then $\Gm (A)$ is exact if and only if $\Gm (A\unit)$ is exact.
\end{lem}

\begin{pff}
There is a split short exact sequence
\[
0 \longrightarrow \Gm (A) \longrightarrow \Gm (A\unit)
    \longrightarrow C (X) \longrightarrow 0.
\]
The result therefore follows from Proposition 7.1 (vi) of \cite{Kr2}.
\end{pff}

\begin{prp}   
Let $X$ be a compact metric space, and let
$A$ be a continuous field of \ca s over $X$ with nuclear fibers $A (x),$
and with $\Gm (A)$ separable. Then $A$ satisfies the conditions of the
previous theorem.
\end{prp}

\begin{pff}
By Theorem 0.3,
the maps $\id_{A (x)}$ are in fact
liftable, so condition (2) of Theorem 4.7
holds.
\end{pff}

The rest of this section is devoted to the proof that if $\Gm (A)$
is separable and exact, then nearby fibers have nearby embeddings
in $\OA{2}.$

\begin{lem}   
Let $X$ be a compact metric
space, let $A$ be a unital continuous field of \ca s over $X$
such that $\Gm (A)$ is separable and exact,
and let $a_1, \dots, a_m$ be continuous sections of $A.$
Let $x_0 \in X.$ Then for every
$\ep > 0$ there exists a neighborhood $U$ of $x_0$ in $X$ such that
for all $x \in U$ there are injective unital \hm s
$\ph_x : A (x_0) \to \OA{2}$ and
$\ps_x : A (x) \to \OA{2},$ and \ucp\  maps $S_x  : A (x_0) \to \OA{2}$
and $T_x  : A (x) \to \OA{2},$ satisfying
\[
\| S_x (a_l (x_0)) - \ps _x (a_l (x)) \| < \ep \andeqn
                  \| T_x (a_l (x)) - \ph_x (a_l (x_0)) \| < \ep
\]
for $1 \leq l \leq m.$
\end{lem}

\begin{pff}
By considering the real and imaginary parts of the given sections
(and using $\ep / 2$ in place of $\ep$), we may assume \wolog\  that
$a_1, \dots, a_m$ are selfadjoint. Similarly, we may assume
that $\|a_l \| \leq 1$ for all $l.$

We next reduce to the case in which $1, a_1 (x_0), \dots, a_m (x_0)$
are linearly independent. In the general case, we may number the
$a_l$ in such a way that $1, a_1 (x_0), \dots, a_{m_0} (x_0)$
are linearly independent and $a_{m_0 + 1} (x_0), \dots, a_m (x_0)$
are linear combinations of $1, a_1 (x_0), \dots, a_{m_0} (x_0).$
Assume the lemma holds for $a_1, \dots, a_{m_0}.$
Set $a_0 = 1.$ For $m_0 + 1 \leq l \leq m,$ write
\[
a_l (x_0) = \sum_{j = 0}^{m_0} \af_{jl} a_j (x_0),
\]
with $\af_{jl} \in \C,$ and then define
$\tilde{a}_l = \sum_{j = 0}^{m_0} \af_{jl} a_j.$
Then the $\tilde{a}_l$ are also continuous sections, and
$\tilde{a}_l (x_0) = a_l (x_0).$ Set
\[
\af = 1 + \max_{m_0 + 1 \leq l \leq m} \sum_{j = 0}^{m_0} | \af_{jl} |.
\]
Choose $U$ so small that the
conclusion of the lemma holds for $a_1, \dots, a_{m_0},$ with
$\ep / (2 \af)$
in place of $\ep,$ and also so small that
$\| \tilde{a}_l (x) - a_l (x) \| < \ep/ 2$ for $x \in U$ and
$m_0 + 1 \leq l \leq m.$ The resulting \hm s $\ph_x$ and $\ps_x,$ and
\ucp\  maps $S_x$ and $T_x,$
then satisfy
\[
\| S_x (a_l (x_0)) - \ps_x (a_l (x)) \| < \frac{\ep}{2 \af} \andeqn
\| T_x (a_l (x)) - \ph_x (a_l (x_0)) \| < \frac{\ep}{2 \af}
\]
for $1 \leq l \leq m_0.$ Hence, for $m_0 + 1 \leq l \leq m$ we have
\[
\| S_x (a_l (x_0)) - \ps_x (a_l (x)) \| \leq
  \sum_{j = 0}^{m_0} | \af_{jl} |
                              \| S_x (a_j (x_0)) - \ps_x (a_j (x)) \| +
  \| \ps_x \| \| \tilde{a}_l (x) - a_l (x) \| < \ep.
\]
Similarly,
\[
\| T_x (a_l (x)) - \ph_x (a_l (x_0)) \| \leq
         \| T_x \| \| a_l (x) - \tilde{a}_l (x) \|
+ \sum_{j = 0}^{m_0} | \af_{jl} | \| T_x (a_j (x)) - \ph_x (a_j (x_0))\|
    < \ep.
\]
This proves the reduction.

We now assume that $1, a_1, \dots, a_m$ are selfadjoint, have norm at
most $1,$ and are linearly independent at $x_0.$ Set
$E = \spn (1, a_1 (x_0), \dots, a_m (x_0));$ then
$E$ is a finite dimensional operator system contained in $A (x_0).$
By local liftability (Theorem 4.7 (2)),
there is a \ucp\  map
$S : E \to \Gm (A)$ such that $\ev_x \circ S = \id_E.$
Let $b_l = S (a_l (x_0)),$ so that $b_l$ is a continuous section
of $A$ satisfying $b_l (x_0) = a_l (x_0).$ Choose an open set
$U_0 \subset X,$ containing $x_0,$ such that
$\| b_l (x) - a_l (x) \| < \ep / 4$ for all $x \in U_0.$
Define $S^{(0)}_x : E \to A (x)$ by $S^{(0)}_x = \ev_x \circ S.$
Then $S^{(0)}_x$ is \uacp, and
$\| S^{(0)}_x (a_l (x_0)) - a_l (x) \| < \ep / 4$
for all $x \in U_0.$
Use Theorem 2.8
to find an injective unital
\hm\   $\ps_x : A (x) \to \OA{2}.$ Since $\OA{2}$ is nuclear,
Proposition 0.4
provides a \ucp\  map
$S_x : A (x_0) \to \OA{2}$ such that
$\| S_x (a_l (x_0)) - \ps_x (S^{(0)}_x (a_l (x_0))) \| < \ep / 4.$
This gives $\| S_x (a_l (x_0)) - \ps_x ( a_l (x)) \| < \ep / 2,$
for all $x \in U_0.$

We still need $T_x.$
Choose, using Lemma 1.10,
an integer $n$ such that whenever
$V : E \to A(x)$ and $W : E \to \OA{2}$ are two \ucp\  maps  such that
$V$ is injective and $\| V^{-1} \otimes \id_{M_n}\| \leq 1 + \ep / 4,$
there is a \ucp\  map
$Q : A (x) \to \OA{2}$ such that $\|Q \circ V - W \| < \ep / 2.$
We now claim that there is an open set
$U_1 \subset X,$ containing $x_0,$ such that
$S^{(0)}_x$ is injective and
$\| (S^{(0)}_x)^{-1} \otimes \id_{M_n}\| \leq 1 + \ep / 4$
for all $x \in U_1.$ Let $\{e_{ij} : 1 \leq i, \, j \leq n\}$ be
a system of matrix units for $M_n.$ Set $N = (m + 1) n^2,$ and
consider the set of $N$ sections
\[
\{c_l : 1 \leq l \leq N\} =
\{e_{ij} \otimes 1 : 1 \leq i, \, j \leq n\} \cup
\{e_{ij} \otimes a_l : 1 \leq l \leq m, \,\, 1 \leq i, \, j \leq n\}
\]
of the continuous field $M_n \otimes A.$
(See Lemma 4.1.)
Note that
$c_1 (x_0), \dots, c_N (x_0)$ are linearly independent. Define
a compact subset $S \subset \C^N$ by
\[
S = \{(\ld_1, \dots, \ld_N) \in \C^N :
   \| \ld_1 c_1 (x_0) + \cdots + \ld_N c_N (x_0) \| = 1\}.
\]
If the claim is false, there is a sequence $(y_k)$ of distinct points in
$X$ such that $y_k \to x_0,$ and such that there are elements
$(\ld_1^{(k)}, \dots, \ld_N^{(k)}) \in S$ satisfying
\[
\| \ld_1^{(k)} c_1 (y_k) + \cdots + \ld_N^{(k)} c_N (y_k) \| <
   \frac{1}{1 + \ep/4}.
\]
Passing to a subsequence, we may assume that
$\ld_l = \limi{k} \ld_l^{(k)}$ exists for each $l.$ By the Tietze
extension theorem, there are continuous functions $f_l$ on $X$
such that $f_l (y_k) = \ld_l^{(k)}$ for all $k$ and
$f_l (x_0) = \ld_l.$ Then $c = f_1 c_1 + \cdots + f_N c_N$ is a
continuous section with $\|c (x_0) \| = 1$ and
$\liminf_{x \to x_0} \|c(x) \| \leq \frac{1}{1 + \ep / 4},$ a
contradiction. This proves the claim.

Set $U = U_0 \cap U_1.$
Let $x \in U.$ Choose (using Theorem 2.8)
some injective unital
\hm\  $\ph_x : A (x_0) \to \OA{2}.$ From above, we have
$\| S^{(0)}_x (a_l (x_0)) - a_l (x) \| < \ep / 4 < \ep / 2.$
Furthermore, $V = S^{(0)}_x : E \to A (x)$ and
$W = \ph_x |_E : E \to \OA{2}$ are \ucp\  maps  such that
$V$ is injective and $\| V^{-1} \otimes \id_{M_n}\| \leq 1 + \ep / 4.$
Therefore, by the choice of $n$ and Lemma 1.10,
there is a
\ucp\  map $T_x : A (x) \to \OA{2}$ such that
$\| T_x \circ S^{(0)}_x |_E - \ph_x |_E \| < \ep / 2.$
It follows that
\[
\| T_x (a_l (x)) - \ph_x (a_l (x_0)) \| \leq
  \| a_l (x) - S^{(0)}_x (a_l (x_0)) \| +
              \| T_x \circ S^{(0)}_x (a_l (x_0)) - \ph_x (a_l (x_0)) \|
< \frac{\ep}{2} + \frac{\ep}{2} = \ep,
\]
as desired.
\end{pff}

\begin{prp}   
Let $X,$ $A,$ and $a_1, \dots, a_m$ be as in Lemma 4.10.
Assume we are given, for each $x \in X,$ an injective unital
\hm\  $\io_x : A (x) \to \OA{2}.$
Define $\rh_0 : X \times X \to [0, \infty)$ by
\[
\rh_0 (x, y) = \inf_T  \max_{1 \leq l \leq m}
         \| T (a_l (x)) - \io_y (a_l (y)) \|,
\]
where the infemum is taken over all \ucp\  maps $T : A(x) \to \OA{2}.$
Then

(1) $\rh_0 (x, y)$ does not depend on the choice of the \hm s $\io_x.$

(2) $\rh_0$ is continuous on $X \times X.$

(3)
$\rh_0 (x, x) = 0$ and $\rh_0 (x, z) \leq \rh_0 (x, y) + \rh_0 (y, z)$
for $x,$ $y,$ $z \in X.$
\end{prp}

\begin{pff}
Part (1) is immediate from the fact (Theorem 1.15)
that
any two injective unital \hm s from $A (y)$ to $\OA{2}$ are \ayue.

We next prove (3). We get $\rh_0 (x, x) = 0$ by taking $T = \io_x.$
For the triangle inequality, let $x,$ $y,$ $z \in X,$ and let $\ep > 0.$
\Wolog\  $\| a_l (y) \| \leq 1$ for all $l.$
Choose \ucp\  maps $S : A (x) \to \OA{2}$ and
$T : A (y) \to \OA{2}$ such that
\[
\max_{1 \leq l \leq m} \| S (a_l (x)) - \io_y (a_l (y)) \|
          \leq \rh_0 (x, y) + \frac{\ep}{3} \andeqn
\max_{1 \leq l \leq m} \| T (a_l (y)) - \io_z (a_l (z)) \|
          \leq \rh_0 (y, z) + \frac{\ep}{3}.
\]
Apply Lemma 1.10
with $A = A (y),$ $B_1 = B_2 = \OA{2},$
$E = \spn (1, a_1 (y), \dots, a_m (y)),$ $\dt = 0,$ $V = \io_y |_E$
(so that $\| V^{-1} \|_{\cb} = 1$), and $W = T |_E,$ to obtain
a \ucp\  map $R : \OA{2} \to \OA{2}$ such that
$\|R \circ \io_y |_E - T |_E \| < \ep / 3.$ Then
$R \circ S : A (x) \to \OA{2}$ is \uacp\  and satisfies
\beqr
\lefteqn{
\| (R \circ S) (a_l (x)) - \io_z (a_l (z)) \|    }  \\
     & \leq & \| R \| \| S (a_l (x)) - \io_y (a_l (y)) \| +
         \| (R \circ \io_y) (a_l (y)) - T (a_l (y)) \| +
         \| T (a_l (y)) - \io_z (a_l (z)) \|    \\
  & < & \left(\rh_0 (x, y) + \frac{\ep}{3}\right) + \frac{\ep}{3} +
               \left(\rh_0 (y, z) + \frac{\ep}{3}\right)
      = \rh_0 (x, y) + \rh_0 (y, z) + \ep.
\eeqr
This shows that $\rh_0 (x, z) \leq \rh_0 (x, y) + \rh_0 (y, z) + \ep,$
and we let $\ep \to 0.$

For continuity (part (2)), let $x_0,$ $y_0 \in X$ and let $\ep > 0.$
Use the previous lemma to
choose neighborhoods $U$ of $x_0$ and $V$ of $y_0$ such that for
$x \in U,$ both $\rh_0 (x_0, x)$ and $\rh_0 (x, x_0)$ are less than
$\ep / 2,$  and similarly $\rh_0 (y_0, y),$ $\rh_0 (y, y_0) < \ep / 2$
for $y \in V.$ Then for $x \in U$ and $y \in V,$
\[
\rh_0 (x, y) \leq \rh_0 (x, x_0) + \rh_0 (x_0, y_0) + \rh_0 (y_0, y)
       < \rh_0 (x_0, y_0) + \ep,
\]
and similarly $\rh_0 (x, y) > \rh_0 (x_0, y_0) - \ep.$
\end{pff}

\begin{rmk}   
Let $X,$ $A,$ and $a_1, \dots, a_m$ be as in Lemma 4.10,
and suppose that all the fibers of $A$ are nuclear. (In this case,
exactness of $\Gm (A)$ is automatic, by Proposition 4.9.)
Then the definition of the function $\rh_0$ of the previous
proposition can be reformulated to look much more like a distance:
\[
\rh_0 (x, y) = \inf_T  \max_{1 \leq l \leq m}
         \| T (a_l (x)) - a_l (y) \|,
\]
where the infemum is taken over all \ucp\  maps $T : A(x) \to A (y).$

To see this, let $\tilde{\rh}_0 (x, y)$ denote the expression on the
right hand side. Obviously, $\rh_0 (x, y) \leq \tilde{\rh}_0 (x, y).$
For the reverse inequality, let $S : A(x) \to \OA{2}$ satisfy
\[
\max_{1 \leq l \leq m} \| S (a_l (x)) - \io_y (a_l (y)) \|
                    < \rh_0 (x, y) + \frac{\ep}{2}.
\]
Since $A (y)$ is nuclear, there are $n$ and \ucp\  maps
$Q_0 : A (y) \to M_n$ and $R : M_n \to A (y)$ such that
$\| R (Q_0 (a_l (y))) - a_l (y) \| < \frac{\ep}{2}$ for
$1 \leq l \leq m.$ The Arveson extension theorem
(Theorem 6.5 of \cite{Pl}) gives a \ucp\  map
$Q : \OA{2} \to M_n$ such that
$Q |_{\io_y (A (y))} = Q_0 \circ \io_y^{- 1}.$ Set
$T = R \circ Q \circ S,$ giving
\beqr
\lefteqn{
\tilde{\rh}_0 (x, y)  \leq
    \max_{1 \leq l \leq m} \| T (a_l (x)) - a_l (y) \|  } \\
  & \leq & \| R \| \| Q \|
            \max_{1 \leq l \leq m} \| S (a_l (x)) - \io_y (a_l (y)) \|
+ \max_{1 \leq l \leq m} \|(R \circ Q \circ \io_y)(a_l (y)) - a_l (y) \|
 < \rh_0 (x, y) + \ep.
\eeqr
\end{rmk}

The square root in the following proposition comes from the one in
Lemma 1.12.
It can't be removed, as follows from Remark 6.11.

\begin{prp}   
Let $X$ and $A$ be as in Lemma 4.10,
and let $u_1, \dots, u_m$ be continuous unitary sections of $A.$
Define $\rh : X \times X \to [0, \infty)$ by
\[
\rh (x, y) = \sup_{\ph, \, \ps} \inf_{v \in U( \OA{2})}
        \max_{1 \leq l \leq m} \| v \ph (u_l (x)) v^* - \ps (u_l (y))\|,
\]
where the supremum runs over the (nonempty) sets of
injective unital \hm s
$\ph : A(x ) \to \OA{2}$ and $\ps : A (y) \to \OA{2}.$
Then $\rh$ is a continuous pseudometric on $X$ (that is, a metric
except that possibly $\rh (x, y)$ could be zero with $x \neq y$).
Moreover, if $\rh_0$ is as in the previous proposition,
using $u_1, \dots, u_m$ in place of $a_1, \dots, a_m,$ then
\[
\rh (x, y) \leq 11 \sqrt{ \max (\rh_0 (x, y), \rh_0 (y, x))}.
\]
\end{prp}

\begin{pff}
It follows from Theorem 2.8
that every $A (x)$ posesses an
injective unital \hm\  to $\OA{2},$ and from Theorem 1.15 
that any two such \hm s are \ayue. We can therefore rewrite the
definition of $\rh$ as follows: For each $x \in X,$ choose and fix
an injective unital \hm\  $\ph_x : A(x) \to \OA{2}.$ Then
\[
\rh (x, y) = \inf_{v \in U( \OA{2})}
   \max_{1 \leq l \leq m} \| v \ph_x (u_l (x)) v^* - \ph_y (u_l (y))\|.
\]
{}From this formula, it is obvious that $\rh$ is a pseudometric.
(Note that $\rh (x, y)$ can be at most $2$ for any $x$ and $y.$)

It remains to prove the estimate
\[
\rh (x, y) \leq 11 \sqrt{ \max (\rh_0 (x, y), \rh_0 (y, x))}.
\]
(Continuity of $\rh$ follows from this estimate and the previous
corollary, using the fact that $\rh$ is a pseudometric.)
Equivalently, we prove that for all $\ep > 0,$ there
is $v  \in U( \OA{2})$ such that
\[
\| v \ph (u_l (x)) v^* - \ps (u_l (y))\| <
    \ep + 11 \sqrt{ \max (\rh_0 (x, y), \rh_0 (y, x))}
\]
for $1 \leq l \leq m.$

Fix $x$ and $y,$ and let $\ep > 0.$ Choose $\ep_0 > 0$ so small that
\[
2 \ep_0 + 11 \sqrt{ \max (\rh_0 (x, y), \rh_0 (y, x)) + 3 \ep_0}
   \leq \ep + 11 \sqrt{ \max (\rh_0 (x, y), \rh_0 (y, x))}.
\]
The definition of $\rh_0$ gives \ucp\  maps
\[
S_0 : \ph_x (A(x)) \to \OA{2} \andeqn
    T_0 : \ph_y (A(y)) \to \OA{2}
\]
such that
\[
\| S_0 (\ph_x (u_l (x))) - \ph_y (u_l (y)) \| < \rh_0 (x, y) + \ep_0
    \andeqn
\| T_0 (\ph_y (u_l (y))) - \ph_x (u_l (x)) \| < \rh_0 (y, x) + \ep_0
\]
for $1 \leq l \leq m.$
Since $\OA{2}$ is nuclear,
it follows from Proposition 0.4
 that there are \ucp\  maps
\[
S : \OA{2} \to \OA{2} \andeqn T : \OA{2} \to \OA{2}
\]
such that
\[
\| S (\ph_x (u_l (x))) - \ph_y (u_l (y)) \| < \rh_0 (x, y) + 2 \ep_0
       \,\,\, {\mathrm{and}} \,\,\,
\| T (\ph_y (u_l (y))) - \ph_x (u_l (x)) \| < \rh_0 (y, x) + 2 \ep_0.
\]
Proposition 1.7
provides isometries $s,$ $t \in \OA{2}$ such that
\[
\| s^* \ph_x (u_l (x)) s - \ph_y (u_l (y)) \| < \rh_0 (x, y) + 3 \ep_0
       \,\,\, {\mathrm{and}} \,\,\,
\| t^* \ph_y (u_l (y)) t - \ph_x (u_l (x)) \| < \rh_0 (y, x) + 3 \ep_0.
\]
Now Lemma 1.12
provides a unitary $z \in \OT{\OA{2}}$ such that
\[
\| z ( 1 \otimes \ph_x (u_l (x))) z^* - 1 \otimes \ph_y (u_l (y)) \| <
         11 \sqrt{ \max (\rh_0 (x, y), \rh_0 (y, x)) + 3 \ep_0}.
\]
Let $\mu : \OT{\OA{2}} \to \OA{2}$ be an isomorphism
(Theorem 0.8).
Then $a \mapsto \mu ( 1 \otimes a )$ is \ayue\   to $\id_{\OA{2}}$
(by Proposition 0.7),
so there is a unitary $w \in \OA{2}$ such that
\[
\| w \mu ( 1 \otimes \ph_x (u_l (x))) w^* - \ph_x (u_l (x)) \| < \ep_0
  \andeqn
\| w \mu ( 1 \otimes \ph_y (u_l (y))) w^* - \ph_y (u_l (y)) \| < \ep_0
\]
for $1 \leq l \leq m.$ Set $v = w \mu (z) w^*.$ Then one checks that
\beqr
\| v \ph_x (u_l (x)) v^* - \ph_y (u_l (y)) \| & < &
    2 \ep_0 + 11 \sqrt{ \max (\rh_0 (x, y), \rh_0 (y, x)) + 3 \ep_0} \\
  &  \leq & \ep + 11 \sqrt{ \max (\rh_0 (x, y), \rh_0 (y, x))}.
\eeqr
\end{pff}

The following two definitions will simplify the notation and
terminology in several lemmas in the next section.

\begin{dfn}   
Let $X$ be a topological space, and let $A$ and $B$ be two continuous
fields of \ca s over $X$ with given continuous sections
$a_1, \dots, a_m$ of $A$ and $b_1, \dots, b_m$ of $B.$ If $\ph$ and
$\ps$ are representations of $A$ and $B$ in a \ca\  $D,$ then we
define the {\em sectional distance} $\ds (\ph, \ps)$ between
$\ph$ and $\ps$ (with respect to $a_1, \dots, a_m$ and
$b_1, \dots, b_m$) to be
\[
\ds (\ph, \ps) = \sup_{x \in X} \max_{1 \leq l \leq m}
          \| \ph_x (a_l (x)) - \ps_x (b_l (x)) \|.
\]
We suppress the sections $a_l$ and $b_l$ in the notation, since they
will always be clear from the context. We use the same notation for
representations $\ph^{(1)}$ and $\ph^{(2)}$ of
two different restrictions $A |_{X \times \{y_1\} }$ and
$A |_{X \times \{y_2\} }$ of a single continuous field over
$X \times Y$ with a single set of sections $a_1, \dots, a_m:$
\[
\ds ( \ph^{(1)}, \ph^{(2)} ) =
          \sup_{x \in X} \max_{1 \leq l \leq m}
      \| \ph^{(1)}_{x} (a_l (x, y_1)) - \ph^{(2)}_{x} (a_l (x, y_2)) \|.
\]
Sometimes $\ph^{(1)}$ and $\ph^{(2)}$ will be restrictions
$\ph |_{X \times \{y_1\} }$ and $\ph |_{X \times \{y_2\} }$ of the same
representation $\ph,$ and the obvious notation will also be used in this
case.
\end{dfn}

\begin{dfn}   
Let $X$ and $Y$ be compact metric spaces, with metric $d_Y$ on $Y.$
Let $A$ be a unital continuous field of \ca s over
$X \times Y,$ with
continuous unitary sections $u_1, \dots, u_m$ such that for each
$(x, y) \in X \times Y$ the elements
$u_1 (x, y), \dots, u_m (x, y)$ generate $A (x, y)$ as a \ca.
Let $\rh : [0, \infty) \to  [0, \infty)$ be a nondecreasing function
with $\lim_{t \to 0} \rh (t) = \rh (0) = 0.$
We say that $A$ is {\em $(X, \rh)$-embeddable}
(with respect to the sections $u_1, \dots, u_m$) in a
unital \ca\  $D$ if:

(1) For every $y \in Y$ there is an injective
unital representation of $A |_{X \times \{y\} }$ in $D.$

(2) If $\ph^{(1)}$ and $\ph^{(2)}$ are injective unital representations
of $A |_{X \times \{y_1\} }$ and $A |_{X \times \{y_2\} }$ in $D,$
then for $\ep > 0$
there exists a continuous unitary function $w : X \to D$ such that
the representation $w \ph^{(1)} w^*,$ given by
$a \mapsto w(x) \ph^{(1)}_x (a) w(x)^*$ for $a \in A (x),$ satisfies
\[
\ds (w \ph^{(1)} w^*, \ph^{(2)} ) < \ep + \rh ( d_Y (y_1,  y_2) ).
\]
\end{dfn}

\begin{rmk}   
Let $Y$ be a compact metric space, and let $A$ be a unital continuous
field of \ca s over $Y,$ with $\Gm (A)$ separable and exact.
We can apply the terminology of the
previous definition to $A$ by taking $X$ to be a one point space, so
that $Y = X \times Y.$
Proposition 4.13
now asserts that if
$m$ continuous unitary sections are given which generate $A (y)$
for each $y,$ then $A$ is $(X, \rh)$-embeddable in $\OA{2}$ for a
suitable $\rh.$ If $\rh_0$ is as in the statement of the
proposition, then
\[
\rh (r) = \sup_{d_Y (y_1, y_2) \leq r}
              11 \sqrt{ \max (\rh_0 (y_1, y_2), \rh_0 (y_2, y_1))}
\]
will work. (The relation $\lim_{t \to 0} \rh (t) = 0$ follows from
the continuity of $\rh_0$ and the compactness of $Y.$)
\end{rmk}

\section{Continuous embedding of exact continuous fields}

Let $A$ be a continuous field of \ca s over a compact metric space $X,$
with $\Gm (A)$ separable, exact, and unital. In the previous section, we
have shown that the fibers of $A$ over nearby points of $X$ have
embeddings in $\OA{2}$ which are close in a suitable sense.
In this section, we use the methods of Haagerup and R\o rdam \cite{HR}
to make, over a sufficiently nice space $X,$ a continuous selection of
these embeddings, so as to obtain a continuous representation  of $A$ in
$\OA{2}.$
It is not clear how general $X$ can be in our argument, but certainly
any compact manifold or finite CW complex is covered.
The methods of Blanchard \cite{Bl} cover more general spaces, but our
approach has the advantage of giving better information about how close
the embeddings of nearby fibers really are. We illustrate this for
$X = [0, 1],$ by showing that if the function $\rh_0$ of
Proposition 4.11
(which gives an abstract distance between
fibers) is $\Lip^{\af},$ then there is a $\Lip^{\af/2}$
representation of $A$ in $\OA{2}.$ In the next section, we will apply
our results to give a $\Lip^{1/2}$ representation of the field
of rotation algebras in $\OA{2},$ which is as good
as the representation of this field in $L (H)$ in \cite{HR}.

The first five lemmas of this section are essentially one dimensional,
with a parameter space carried along.
They enable us to treat the case $X = [0, 1]^n$ by induction on $n.$
Four of these lemmas are simply modifications of lemmas in \cite{HR},
the modifications being the parameter space, the need to make do
with approximate commutativity in some places where \cite{HR} has
exact commutativity, and the need to handle more general distance
estimates than $\Lip^{1/2}.$

The following definition is useful for describing our version of
Lemma 5.1 of \cite{HR}.

\begin{dfn}   
Let $X$ be a topological space, let $E$ be a Banach space, and let
$[ \af, \bt ]$ be an interval in $\R.$ A function
$\xi : X \times [ \af, \bt ] \to E$ will be called
{\em smooth in the $[ \af, \bt ]$ direction} if for every
$n$ the $n$-th derivative
$\displaystyle{\frac{d^n}{dt^n}} \xi (x, t)$ exists for every
$(x, t) \in X \times [ \af, \bt ],$ and is jointly continuous in
$x$ and $t.$ The function $\xi$ is
{\em piecewise smooth in the $[ \af, \bt ]$ direction} if there is a
partition $\af = t_0 < t_1 < \cdots < t_n = \bt$ such that
$\xi |_{X \times [t_{j - 1}, t_j]}$ is smooth
in the $[t_{j - 1}, t_j]$ direction for $1 \leq j \leq n.$
\end{dfn}

It might be more appropriate to allow the break points $t_j$ in
the definition of piecewise smoothness to depend continuously
on $x \in X,$ but the definition we give suffices for our purposes.

\begin{lem}   
Let $A$ be a unital \ca, and let $B \subset A$ be a unital
subalgebra with $B \cong \OA{2}.$ Let $X$ be a topological space,
and let $v : X \to U( B' \cap A)$ be a continuous function from
$X$ to the unitary group of $B' \cap A.$
Then there is a function $u : X \times [0, 1] \to U(B)$ which is smooth
in the $[0, 1]$ direction and such that for all $x \in X$ we have:

(1) $u(x, 0) = 1$ and $u (x, 1) = v (x).$

(2) $\left\| {\frac{d}{dt}} u (x, t) \right\| \leq 9$
  for all $t \in [0, 1].$

(3) $\| [u(x, t), a ] \| \leq 4 \| [v (x), a] \|$
   for all $t \in [0, 1]$ and $a \in B' \cap A.$

(4) $\left\| \left[
             {\frac{d}{dt}} u (x, t), a \right] \right\|
       \leq 9 \| [v (x), a] \|$
   for all $t \in [0, 1]$ and $a \in B' \cap A.$

(5) $\left\|{\frac{d}{dt}} \left[ u(x, t) a u(x, t)^* \right] \right\|
         \leq 45 \| [v (x), a] \|$
   for all $t \in [0, 1]$ and $a \in B' \cap A.$
\end{lem}

\begin{pff}
The proof of Lemma 5.1 of \cite{HR} works essentially as is, using
$B' \cap A$ in place of $M$ and $B \cong \OA{2}$ in place of $M',$
and just carrying along the extra parameter $x.$
The homomorphisms used there become $\pi, \, \rh : M_3 \to B$
and $\tilde{\pi}, \, \tilde{\rh} : M_3 \otimes (B' \cap A) \to A.$
We take $w$ and $h$ as there, and the element called $v$ there
becomes $z (x) = \diag (v(x)^*, 1, v(x) ).$ We define
\[
u (x, t) =
    \tilde{\pi} \left( \rule{0em}{3ex}
                         \exp (i t h) z(x) \exp (-i t h) z(x)^* \right)
\cdot  \tilde{\rh} \left( \rule{0em}{3ex}
                         \exp (i t h) z(x) \exp (-i t h) z(x)^* \right).
\]
The proofs of the estimates are exactly the same as in \cite{HR}.
\end{pff}

\begin{lem}   
Let $X$ be a compact Hausdorff space, let $v : X \to U (\OA{2})$ be
continuous, let $S \subset \OA{2}$ be compact, and let $\ep > 0.$
Then there is a function $u : X \times [0, 1] \to U (\OA{2})$
which is piecewise smooth in the $[0, 1]$ direction,
and such that for all $x \in X$ we have:

(1) $u(x, 0) = 1$ and $u (x, 1) = v (x).$

(2) $\left\| {\frac{d}{dt}} u (x, t) \right\|
            \leq 9 + \ep$ for all $t \in [0, 1].$

(3) $\| [u(x, t), a ] \| \leq 4 \| [v (x), a] \| + \ep$
   for all $t \in [0, 1]$ and $a \in S.$

(4) $\left\| \left[
             {\frac{d}{dt}} u (x, t), a \right] \right\|
    \leq 9 \| [v (x), a] \| + \ep$ for all $t \in [0, 1]$ and $a \in S.$

(5) $\left\| {\frac{d}{dt}} [u(x, t) a u(x, t)^*] \right\|
          \leq 45 \| [v (x), a] \| + \ep$
   for all $t \in [0, 1]$ and $a \in S.$
\end{lem}

\begin{pff}
Let $R = \sup_{a \in S} \|a\|.$ Choose $\ep_1 > 0$ so small that
\[
(10 + 8 R) \ep_1 + \ep_1^2 < \ep.
\]
Choose $0 < \dt < 1$ so small that
\[
9 \left( \frac{1}{1 - \dt} - 1 \right) < \ep_1  \andeqn
9 \left( \frac{1}{1 - \dt} - 1 \right) \cdot 2 R < \frac{\ep_1}{2}.
\]
Choose $0 < \ep_0 < 1$ so small that the quantities
\[
6 \ep_0, \,\,\,\,\,\, 6 R \ep_0, \,\,\,\,\,\,
\frac{18}{1 - \dt} \cdot 6 \ep_0, \,\,\,\,\,\,
\frac{2}{\dt} \arcsin \left( \frac{3 \ep_0}{2} \right), \,\,\,\,\,\,
{\mathrm{and}} \,\,\,\,\,\,
 \frac{4}{\dt} \arcsin \left( \frac{3 \ep_0}{2} \right) R
\]
are all less than $\ep_1.$
Choose a finite set $F \subset S$ such that every element of $S$
is within $\ep_0$ of an element of $F,$ and a finite subset
$G \subset U (\OA{2})$  such that every $v (x),$ for $x \in X,$ is
within $\ep_0$ of an element of $G.$

Let $\ph : \OT{\OA{2}} \to \OA{2}$ be an isomorphism (Theorem   0.8).
Then by Proposition 0.7
the \hm\  $a \to \ph (1 \otimes a)$ is \ayue\  to
$\id_{\OA{2}}.$ Therefore there is a unitary $w \in \OA{2}$ such
that $\|w \ph (1 \otimes b) w^* - b \| < \ep_0$ for
$b \in F \cup G.$ Replacing $\ph$ by $w \ph ( \cdot ) w^*,$
we may assume that $\| \ph (1 \otimes b) - b \| < \ep_0$ for
$b \in F \cup G.$ If now $a \in S,$ then consideration of
$b \in F$ with $\| a - b \| < \ep_0$ shows that
$\| \ph (1 \otimes a) - a \| < 3 \ep_0.$ Similarly,
$\| \ph (1 \otimes v(x)) - v(x) \| < 3 \ep_0$ for all $x \in X.$

Use the previous lemma, with $A = \OT{\OA{2}}$ and $B = \OT{\C},$
to choose a smooth (in the $[0, 1]$ direction)
function $u_0 : X \times [0, 1] \to U (\OT{\OA{2}})$
for the unitary $x \mapsto 1 \otimes v (x).$
Define $h (x) = - i \log (v (x)^* \ph (1 \otimes v (x))),$ and note that
$\| v (x)^* \ph (1 \otimes v (x)) - 1 \| < 3 \ep_0,$ so that
$\| h (x) \| < 2\arcsin ( \frac{3 \ep_0}{2} ).$ Then set
$u_1 (x, t) = \exp (i t h (x)) \ph (1 \otimes v (x)).$ Now define
\[
u (x, t) = \left\{ \begin{array}{ll}
    \ph \left( u_0 \left( x, \frac{t}{1 - \dt} \right) \right) &
                          \hspace{3em} 0 \leq t \leq 1 - \dt \\
    u_1 \left( x, \frac{1}{\dt} (t - 1) + 1\right)  &
                          \hspace{3em} 1 - \dt \leq t \leq 1.
  \end{array} \right.
\]
Clearly $(x, t) \mapsto u (x, t)$ is a continuous and piecewise smooth
(in the sense of Definition 5.1)
unitary path from $1$ to $v (x).$

For $0 \leq t \leq 1 - \dt,$ we now estimate the quantities in parts
(2) through (4) of the conclusion. We have
\[
\left\| {\frac{d}{dt}} u (x, t) \right\| =
 \frac{1}{1 - \dt}
                 \left\| {\frac{d}{dt}} u_0 (x, t) \right\|
  \leq \frac{9}{1 - \dt} \leq 9 + \ep_1.
\]
Next, for $a \in S,$ we have
\beqr
\| [u (x, t), a] \| & \leq & 2 \| \ph (1 \otimes a) - a \| +
      \left\| \ph \left( \left[ u_0 \left( x, \frac{t}{1 - \dt} \right),
                            1 \otimes a \right] \right) \right\|
                \\
 & \leq & 6 \ep_0 + 4 \| [v (x), a] \| \leq \ep_1 + 4 \| [v (x), a] \|
\eeqr
and
\beqr
\left\| \left[{\frac{d}{dt}} u (x, t), a  \right] \right\|
    & \leq & 2 \left\| {\frac{d}{dt}} u (x, t) \right\|
                                         \| \ph (1 \otimes a) - a \| +
         \left\| \ph \left( \left[ {\frac{d}{dt}}
          \left( u_0 \left(x, \frac{t}{1 - \dt} \right) \right),
                            1 \otimes a \right] \right) \right\|
           \\
 & \leq  & 6 \ep_0 \cdot \frac{9}{1 - \dt}
                        + \frac{1}{1 - \dt} \cdot 9  \| [v (x), a] \|
           \\
 & \leq & \frac{\ep_1}{2} +
  9 \left( \frac{1}{1 - \dt} - 1 \right) \cdot 2 R + 9  \| [v (x), a] \|
  \leq    \ep_1 + 9  \| [v (x), a] \|.
\eeqr

Estimating these same quantities for $1 - \dt \leq t \leq 1$ instead,
we obtain
\[
\left\| {\frac{d}{dt}} u (x, t) \right\| \leq
   \frac{1}{\dt} \| h \| \leq
    \frac{2}{\dt} \arcsin \left( \frac{3 \ep_0}{2} \right)
   < \ep_1 \leq 9 + \ep_1,
\]
\beqr
\| [u (x, t), a] \| & \leq & 2R \| u(x, t) - v (x) \| + \| [v (x), a] \|
  \leq  2R \| \ph (1 \otimes v (x)) - v (x) \| + \| [v (x), a] \|
           \\
  & \leq & 6R \ep_0 + \| [v (x), a] \| \leq \ep_1 + 4 \| [v (x), a] \|,
\eeqr
and
\[
\left\| \left[{\frac{d}{dt}} u (x, t), a\right] \right\| =
     \frac{1}{\dt} \| [ih (x), a] \| \leq \frac{2}{\dt} \| h \| \| a \|
  \leq \frac{4}{\dt} \arcsin \left( \frac{3 \ep_0}{2} \right) R
  \leq \ep_1 + 9  \| [v (x), a] \|.
\]

We thus obtain in both cases, for $a \in S,$
\[
\left\| {\frac{d}{dt}} u (x, t) \right\| \leq 9 + \ep_1,
      \,\,\,
          \| [u (x, t), a] \| \leq 4 \| [v (x), a] \| + \ep_1,
   \,\,\, {\mathrm{and}} \,\,\,
  \left\| \left[{\frac{d}{dt}} u (x, t), a \right] \right\|
                            \leq 9 \| [v (x), a] \| + \ep_1.
\]
Estimates (2), (3), and (4) follow since $\ep_1 < \ep.$
Moreover, following the reasoning at the end of the proof of Lemma 5.1
of \cite{HR}, we obtain for $a \in S$ (using also $\ep_0 \leq 1$)
\beqr
\lefteqn{
\left\| {\frac{d}{dt}} [ u(x, t) a u(x, t)^* ] \right\|
                                 } \\
         & \leq &
   \left\| \left[
      {\frac{d}{dt}} u (x, t), a\right] \right\|
                                               \| u(x, t)^* \|
 + \| [u (x, t), a] \|
            \left\|  {\frac{d}{dt}} u (x, t)^* \right\|
                  \\
  & \leq & 9 \| [v (x), a] \| + \ep_1 +
      (4 \| [v (x), a] \| + \ep_1) (9 + \ep_1)
   \leq 45 \| [v (x), a] \| + (10 + 4 \| [v (x), a] \|) \ep_1 + \ep_1^2
                  \\
 & \leq & 45 \| [v (x), a] \| + (10 + 8 R) \ep_1 + \ep_1^2
          < 45 \| [v (x), a] \| + \ep,
\eeqr
which proves (5).
\end{pff}

\begin{lem}   
Let $X$ be a compact metric space. Let $A$ and $B$ be two continuous
fields of \ca s over $X \times [0, 1],$ and let $\af$ and $\bt$
be representations of 
$A$ and $B$ in $\OA{2}.$ Let $u_1, \dots, u_m$ be continuous unitary
sections of $A,$ let $v_1, \dots, v_m$ be continuous unitary sections
of $B,$ and suppose there is $r > 0$ such that for every $t \in [0, 1]$
there is a continuous unitary function $z : X \to U (\OA{2})$ satisfying
$\ds (z (\af |_{X \times \{t\} }) z^*, \bt |_{X \times \{t\} } ) < r.$
Then there is a continuous unitary
function $w : X \times [0, 1] \to U ( \OA{2})$ such that
$\ds (w \af w^*, \bt) < 10 r.$
Moreover, given continuous functions $c_0$, $c_1 : X \to U (\OA{2})$
such that
$\ds (c_i (\af |_{X \times \{i\} }) c_i^*, \bt |_{X \times \{i\} }) < r$
for $i = 0,$ $1,$ we may choose $w$ to satisfy $w (x, i) = c_i (x).$
\end{lem}

\begin{pff}
Choose a partition
$0 = t_0 < t_1 < \cdots < t_n = 1$ of $[0, 1]$ such that
\[
\ds (\af |_{X \times \{ t \} }, \af |_{X \times \{ t_j \} } )
                                        < \frac{r}{21} \andeqn
\ds (\bt |_{X \times \{ t \} }, \bt |_{X \times \{ t_j \} } )
                                                    < \frac{r}{21}
\]
for $t \in [t_{j - 1}, t_j].$
By hypothesis, there are continuous unitary functions
$z_0, \dots, z_n : X \to U (\OA{2})$ such that
$\ds (z_j (\af |_{X \times \{ t_j \} }) z_j^*,
                         \bt |_{X \times \{ t_j \} }) < r.$
If unitaries $c_i$ are given, take $z_0 = c_0$ and $z_n = c_1.$
Now estimate
\beqr
\lefteqn{
\| [ z_{j - 1} (x)^* z_j (x), \af_{x, t_j} (u_l (x, t_j)) ] \|
  \leq
\ds ( (z_{j - 1}^* z_j) (\af |_{X \times \{ t_j \} })
                                              (z_{j - 1}^* z_j)^*,
                     \af |_{X \times \{ t_j \} } )      }
                         \\
 & \leq &  \ds (z_j (\af |_{X \times \{ t_j \} }) z_j^*,
                                      \bt |_{X \times \{ t_j \} })
 + \ds (\bt |_{X \times \{ t_j \} }, \bt |_{X \times \{ t_{j - 1} \} })
                \\
 & & \mbox{}
 + \ds (z_{j - 1}^* (\bt |_{X \times \{ t_{j - 1} \} }) z_{j - 1},
                                  \af |_{X \times \{ t_{j - 1} \} })
 + \ds (\af |_{X \times \{ t_{j - 1} \} }, \af |_{X \times \{ t_j \} })
                              \\
 & < & 2 r + \frac{2 r}{21}.
\eeqr
Combining this with the estimates at the beginning of the proof, we get
\[
\| [ z_{j - 1} (x)^* z_j (x), \af_{x, t} (u_l (x, t)) ] \|
      < 2 r + \frac{4 r}{21}
\]
for $t \in [t_{j - 1}, t_j],$ $x \in X,$ and $1 \leq l \leq m.$

By Lemma 5.3,
there are continuous unitary functions
$\tilde{z} : X \times [t_{j - 1}, t_j] \to \OA{2}$ with
\[
\tilde{z} ( x, t_{j - 1}) = 1 \andeqn
                  \tilde{z} (x, t_j) = z_{j - 1} (x)^* z_j (x),
\]
such that (from part (3) of the conclusion)
\beqr
\ds ((\tilde{z} |_{X \times \{t\}}) (\af |_{X \times \{t\}})
                                    (\tilde{z} |_{X \times \{t\}})^*,
                          \af |_{X \times \{t\}}) & = &
\sup_{x \in X} \| [\tilde{z} (x, t), \af_{x, t} (u_l (x, t)) ] \|
          \\
 & < & 4 \left(2 r + \frac{4 r}{21} \right) + \frac{r}{21}
    = 8 r + \frac{17 r}{21}
\eeqr
for $t \in [t_{j - 1}, t_j]$ and $1 \leq l \leq m.$
(We can reparametrize the interval because we don't use the estimates
on the derivatives.)
Now define $w (x, t) = z_{j - 1} (x) \tilde{z} (x, t)$
for $t \in [t_{j - 1}, t_j].$ Then $w$ is continuous, and for each $t,$
\beqr
\lefteqn{ \ds ((w |_{X \times \{t\}}) (\af |_{X \times \{t\}})
                                          (w |_{X \times \{t\}})^*,
                           \bt |_{X \times \{t\}}) } \\
 & \leq  & \ds ((\tilde{z} |_{X \times \{t\}}) (\af |_{X \times \{t\}})
                                     (\tilde{z} |_{X \times \{t\}})^*,
                        \af |_{X \times \{t\}})
   + \ds (\af |_{X \times \{t\}}, \af |_{X \times \{ t_{j - 1} \} })
                                     \\
 & & \mbox{}
   + \ds (\bt |_{X \times \{t\}}, \bt |_{X \times \{ t_{j - 1} \} })
   + \ds (z_{j - 1} (\af |_{X \times \{ t_{j - 1} \} }) z_{j - 1}^*,
                                  \bt |_{X \times \{ t_{j - 1} \} })
        \\
 & < & \left( 8 r + \frac{17 r}{21}  \right)
                               + \frac{2 r}{21} + \frac{2 r}{21} + r 
                                                     = 10 r,
\eeqr
as desired.
Furthermore, $w (x, 0) = z_0 (x)$ and $w (x, 1) = z_n (x)$ for all $x.$
\end{pff}

The next lemma is an analog of Lemma 5.2 of \cite{HR}. There is one
additional twist, namely the number $n',$ which is necessary in the
absence of a Lipschitz condition.

\begin{lem}   
Let $X$ be a compact metric space, let $A$ be a unital continuous field
of \ca s over
$X \times [0, 1],$ and let $u_1, \dots, u_m$ be
continuous unitary sections of $A$ such that for each
$(x, t) \in X \times [0, 1]$ the elements
$u_1 (x, t), \dots, u_m (x, t)$ generate $A (x, t)$ as a \ca.
Assume that $A$ is $(X, \rh)$-embeddable (Definition 4.15)
in $\OA{2}$ for some $\rh,$ using the usual metric on $[0, 1].$
Let $\ph^{(0)}$ and $\ph^{(1)}$ be injective unital representations
of $A |_{X \times \{t_0\} }$ and $A |_{X \times \{t_1\} }$ in $\OA{2},$
with
\[
\ds (\ph^{(0)} |_{X \times \{t_0\} }, \ph^{(1)} |_{X \times \{t_1\} })
       < d_0
\]
for some $d_0 > \rh ( t_1 - t_0).$ Let $0 < n' \leq n$ be integers, and
set $s_j = t_0 + j (t_1 - t_0) / n.$ Then there are
injective unital representations $\gm^{(j)}$ of
$A |_{X \times \{s_j\} }$ such that $\gm^{(0)} = \ph^{(0)},$
$\gm^{(n)} = \ph^{(1)},$ and
\[
\ds (\gm^{(j - 1)}, \gm^{(j)} ) <
    91 \rh \left(\frac{t_1 - t_0}{n} \right) + \frac{90}{n'} \cdot d_0
                     \andeqn
\ds (\gm^{(j)}, \ph^{(0)}) <
    91 n' \rh \left(\frac{t_1 - t_0}{n} \right) + 91 d_0
\]
for all $j.$
\end{lem}

\begin{pff}
The proof is easy if $\rh \left(\frac{t_1 - t_0}{n} \right) = 0.$ So
assume $\rh \left(\frac{t_1 - t_0}{n} \right) > 0.$ Choose $\ep > 0$
such that
\[
\ep < \frac{1}{2 n' + 1} \min \left(
  \rh \left(\frac{t_1 - t_0}{n} \right), d_0 - \rh (t_1 - t_0)  \right).
\]
Choose $n''$ and integers $0 = j (0) < j (1) < \cdots < j (n'') = n$
such that $n' \leq j (r) - j (r - 1) < 2 n'$ for all $r.$
We construct the $\gm^{(j)}$ for $j (r - 1) < j \leq j (r)$ in blocks,
by induction on $r.$

We start with $r = 1.$  Using both the existence of
injective unital representations and the sectional distance estimates
for them in the definition of $(X, \rh)$-embeddability, construct
injective unital representations $\bt^{(j)}$ of
$A |_{X \times \{s_j\} }$ for $0 \leq j \leq j (1)$ such that
\[
\bt^{(0)} = \ph^{(0)} \andeqn \ds ( \bt^{(j - 1)}, \bt^{(j)} ) <
      \ep + \rh \left(\frac{t_1 - t_0}{n} \right).
\]
(This is done by induction. Take $\bt^{(0)} = \ph^{(0)}.$ Choose
any injective unital representation $\bt$ of $A |_{X \times \{s_1\} },$
and find $z : X \to U (\OA{2})$ such that
$\ds (\bt^{(0)}, z \bt z^* ) <
                    \ep + \rh \left(\frac{t_1 - t_0}{n} \right).$
Set $\bt^{(1)} = z \bt z^*.$ Then choose an
injective unital representation $\bt$ of $A |_{X \times \{s_2\} },$
etc.)

For $0 \leq i, \, j \leq j (1)$ we then have
\[
\ds (\bt^{(i)}, \bt^{(j)} ) <
    | i - j| \left(\ep + \rh \left(\frac{t_1 - t_0}{n} \right) \right).
\]
The approximate unitary equivalence part of the hypotheses implies that
there is a continuous function $w_1 : X \to U (\OA{2})$ such that
\[
\ds (w_1 \bt^{(j (1))} w_1^* , \bt^{(0)} ) <
      \ep + \rh \left(j (1) \cdot \frac{t_1 - t_0}{n} \right) \leq
      \ep + \rh (t_1 - t_0).
\]
Combining the last two inequalities, we obtain, for
$0 \leq j \leq j (1),$
\beqr
\sup_{x \in X} \max_{1 \leq l \leq m}
                           \| [ w_1 (x), \bt^{(j)}_x (u_l (x, s_j)) ] \|
 & \leq & \ds (\bt^{(j)}, \bt^{(j (1))} ) + \ds (\bt^{(j)}, \bt^{(0)} )
           + \ep + \rh (t_1 - t_0)  \\
 & < & 2 n' \rh \left(\frac{t_1 - t_0}{n} \right) + \rh (t_1 - t_0).
\eeqr
Now apply Lemma 5.3
to obtain a function
$\tilde{w}_1 : X \times [0, 1] \to U (\OA{2})$ which is piecewise
smooth in the $[0, 1]$ direction and satisfies
\[
\tilde{w}_1 (x, 0) = 1, \,\,\,\,\,\,  \tilde{w}_1 (x, 1) = w_1 (x),
  \andeqn
\left\| \displaystyle{\frac{d}{dt}}
       \left[ \rule{0em}{2.2ex}
                   \tilde{w}_1 (x, t) a \tilde{w}_1 (x, t)^* \right] \right\|
           \leq 45 \| [w_1 (x), a] \| + \ep
\]
for all $x$ and $t,$ and for $a$ in the compact set
\[
S = \{ \bt^{(j)} (u_l (x, s_j)) :1 \leq l \leq m, \,\,
         0 \leq j \leq j (1), \,\, x \in X\}.
\]

Define
\[
\gm^{(j)}_x (a) = \tilde{w}_1 \left(x, \frac{j}{j(1)} \right) \cdot
     \bt^{(j)}_x (a) \cdot \tilde{w}_1 \left(x, \frac{j}{j(1)} \right)^*
\]
for $a \in A (x, s_j)$ and $0 \leq j \leq j (1).$
This immediately gives $\gm^{(0)} = \bt^{(0)} = \ph^{(0)}$ and
\[
\ds (\gm^{(j (1))}, \ph^{(0)} ) < \ep + \rh (t_1 - t_0).
\]
Moreover, using the derivative estimate from the previous paragraph,
as well as $n' \leq j(1) < 2 n',$
we obtain, for $1 \leq l \leq m,$ $1 \leq j \leq j (1),$ and
$x \in X,$
\beqr
\lefteqn{
\| \gm^{(j)}_x (u_l (x, s_j)) - \gm^{(j - 1)}_x (u_l (x, s_{j - 1})) \|
                                             }  \\
 & \leq &
\| \bt^{(j)}_x (u_l (x, s_j)) - \bt^{(j - 1)}_x (u_l (x, s_{j - 1})) \|
 + \int_{(j - 1) / j(1)}^{j / j(1)}
  \left\| \frac{d}{dt}
            \left[\tilde{w}_1 (x, t) \bt^{(j)}_x (u_l (x, s_j))
                               \tilde{w}_1 (x, t)^* \right] \right\| dt
               \\
 & < & \ep + \rh \left(\frac{t_1 - t_0}{n} \right) +
\frac{1}{j(1)} \left[ 45 \left(
     2 n' \rh \left(\frac{t_1 - t_0}{n} \right) + \rh (t_1 - t_0)
                                   \right) + \ep \right]  \\
 & \leq & 2 \ep + 91 \rh \left(\frac{t_1 - t_0}{n} \right)
             + \frac{45}{n'} \rh (t_1 - t_0).
\eeqr
Thus
\[
\ds (\gm^{(j - 1)}, \gm^{(j)} ) <
     2 \ep + 91 \rh \left(\frac{t_1 - t_0}{n} \right)
               + \frac{45}{n'} \rh (t_1 - t_0) <
  91 \rh \left(\frac{t_1 - t_0}{n} \right) + \frac{90}{n'} \cdot d_0.
\]

The induction step for $r < n''$ is essentially the same.
We start with
$\gm^{(j (r - 1))}$ satisfying
\[
\ds (\gm^{(j (r - 1))}, \ph^{(0)} ) < \ep + \rh (t_1 - t_0).
\]
Construct, as at the beginning of the initial step,
injective unital representations $\bt^{(j)}$ of
$A |_{X \times \{s_j\} }$ for $j (r - 1) \leq j \leq j (r)$ such that
\[
\bt^{(j (r - 1))} = \gm^{(j (r - 1))} \andeqn
  \ds ( \bt^{(j - 1)}, \bt^{(j)} ) <
      \ep + \rh \left(\frac{t_1 - t_0}{n} \right).
\]
(This redefines $\bt^{(j (r - 1))}$ from the previous step, but we are
done with the earlier one.)
We obtain the same estimate on $\ds (\bt^{(i)}, \bt^{(j)} )$ as before,
and then the \aue\  part of the hypotheses provides
$w_{r} : X \to U (\OA{2})$ such that
\[
\ds (w_r \bt^{(j (r))} w_r^* , \ph^{(0)} ) <
      \ep + \rh \left(j (r) \cdot \frac{t_1 - t_0}{n} \right) \leq
      \ep + \rh (t_1 - t_0).
\]
It follows that
\[
\ds (w_r \bt^{(j (r))} w_r^* , \bt^{(j (r - 1))} ) <
      2 (\ep + \rh (t_1 - t_0)).
\]
The commutator estimate in the initial step becomes
\beqr
\sup_{x \in X} \max_{1 \leq l \leq m}
                      \| [ w_{r} (x), \bt^{(j)}_x (u_l (x, s_j)) ] \|
 & \leq & \ds (\bt^{(j)}, \bt^{(j (r))} ) +
                                \ds (\bt^{(j)}, \bt^{( j (r - 1))} )
           + 2 (\ep + \rh (t_1 - t_0))  \\
 & < & 2 n' \rh \left(\frac{t_1 - t_0}{n} \right) + 2 \rh (t_1 - t_0)
\eeqr
for $j (r - 1) \leq j \leq j (r).$ Now apply Lemma 5.3
as before,
obtaining $\tilde{w}_{r} : X \times [0, 1] \to U (\OA{2}).$ (In the
definition of $S,$ now use $j (r - 1) \leq j \leq j (r).$)

Define
\[
\gm^{(j)}_x (a) =
     \tilde{w}_{r} \left(x, \frac{j - j(r - 1)}{j(r) - j(r - 1)} \right)
         \cdot \bt^{(j)}_x (a) \cdot
   \tilde{w}_{r} \left(x, \frac{j - j(r - 1)}{j(r) - j(r - 1)} \right)^*
\]
for $a \in A (x, s_j)$ and $j (r - 1) \leq j \leq j (r).$
This gives the same $\gm^{( j (r - 1))}$ that we already have, and
also immediately gives
\[
\ds (\gm^{(j (r))}, \ph^{(0)} ) < \ep + \rh (t_1 - t_0).
\]
Set $k = j(r) - j(r - 1) \geq n'.$ The analog of the second last
estimate in the initial step is then:
\beqr
\lefteqn{
\| \gm^{(j)}_x (u_l (x, s_j)) - \gm^{(j - 1)}_x (u_l (x, s_{j - 1})) \|
            }      \\
 & \leq &
\| \bt^{(j)}_x (u_l (x, s_j)) - \bt^{(j - 1)}_x (u_l (x, s_{j - 1})) \|
                   \\
 & & \mbox{}
 + \int_{(j - 1 - j (r - 1)) / k}^{(j - j (r - 1)) / k}
  \left\| \frac{d}{dt} \left[
     \tilde{w}_{r} (x, t) \bt^{(j)}_x (u_l (x, s_j))
                              \tilde{w}_{r} (x, t)^* \right] \right\| dt
                             \\
 & < & \ep + \rh \left(\frac{t_1 - t_0}{n} \right) +
\frac{1}{k} \left[ 45 \left(
     2 n' \rh \left(\frac{t_1 - t_0}{n} \right) + 2 \rh (t_1 - t_0)
                                   \right) + \ep \right]  \\
 & \leq & 2 \ep + 91 \rh \left(\frac{t_1 - t_0}{n} \right)
             + \frac{90}{n'} \rh (t_1 - t_0).
\eeqr
Thus again
\[
\ds (\gm^{(j)}, \gm^{(j - 1)} ) <
     2 \ep + 91 \rh \left(\frac{t_1 - t_0}{n} \right)
               + \frac{90}{n'} \rh (t_1 - t_0) <
  91 \rh \left(\frac{t_1 - t_0}{n} \right) + \frac{90}{n'} \cdot d_0.
\]

For the final step (with $r = n''$), we do things slightly differently.
Choose $\bt^{(j)}$ for $j (n'' - 1) \leq j \leq j (n'') = n$ as before,
but then choose $w_{n''}$ to satisfy
\[
\ds (w_{n''} \bt^{(n)} w_{n''}^* , \ph^{(1)} ) < \ep.
\]
(We have used $\ph^{(1)}$ in place of $\ph^{(0)}.$
Both $\bt^{(n)}$ and $\ph^{(1)}$ are injective representations
of $A |_{X \times \{t_1\}}$ in $\OA{2}.$)
This now gives
\beqr
\ds (w_{n''} \bt^{(n)} w_{n''}^* , \bt^{(j (n'' - 1))} )
  & \leq & \ds (w_{n''} \bt^{(n)} w_{n''}^* , \ph^{(1)}) +
      \ds (\ph^{(1)}, \ph^{(0)} ) +
      \ds ( \ph^{(0)}, \bt^{(j (n'' - 1))} )   \\
  &   < & 2 \ep + \rh (t_1 - t_0) + d_0,
\eeqr
so
\beqr
\sup_{x \in X} \max_{1 \leq l \leq m}
                      \| [ w_{n''} (x), \bt^{(j)} (u_l (x, s_j)) ] \|
 & \leq & \ds (\bt^{(j)}, \bt^{(n)} ) +
                                \ds (\bt^{(j)}, \bt^{( j (n'' - 1))} )
           + 2 \ep + \rh (t_1 - t_0) + d_0  \\
 & < & 2 n' \rh \left(\frac{t_1 - t_0}{n} \right) +
                                             \rh (t_1 - t_0) + d_0.
\eeqr
Choose $\tilde{w}_{n''}$ as in the induction step, and define
$\gm^{(j)}$ as
there for $j (n'' - 1) \leq j < n.$ Define $\gm^{(n)} = \ph^{(1)}.$
Set $k = n - j (n'' - 1) \geq n'.$ Then the estimate of
$\ds (\gm^{(j)}, \gm^{(j - 1)} )$ in the induction step becomes,
for $j (n'' - 1) \leq j < n,$
\beqr
\lefteqn{
\ds (\gm^{(j - 1)}, \gm^{(j)} )
      }    \\
 & < & \ep + \rh \left(\frac{t_1 - t_0}{n} \right) +
    \frac{1}{j (n'') - j (n'' - 1)} \left[ 45 \left(
         2 n' \rh \left(\frac{t_1 - t_0}{n} \right) 
                          + \rh (t_1 - t_0) + d_0 \right) + \ep \right]
                                    \\
 & \leq & 2 \ep + 91 \rh \left(\frac{t_1 - t_0}{n} \right)
             + \frac{45}{n'} \cdot (\rh (t_1 - t_0) + d_0).
\eeqr
For $j = n,$ the corresponding estimate requires one extra term,
namely $\ds (w_{n''} \bt^{(n)} w_{n''}^* , \ph^{(1)}) < \ep,$
so
\[
\ds (\gm^{(n)}, \gm^{(n - 1)} ) <
  3 \ep + 91 \rh \left(\frac{t_1 - t_0}{n} \right) +
                           \frac{45}{n'} \cdot (\rh (t_1 - t_0) + d_0).
\]
In either case, we still have
\[
\ds (\gm^{(j - 1)}, \gm^{(j)} ) <
    91 \rh \left(\frac{t_1 - t_0}{n} \right) + \frac{90}{n'} \cdot d_0.
\]

Our inductive construction is now complete. We have
\[
\ds (\gm^{(j - 1)}, \gm^{(j)} ) <
    91 \rh \left(\frac{t_1 - t_0}{n} \right) + \frac{90}{n'} \cdot d_0
\]
for all $j$ with $1 \leq j \leq n,$ and also, since
$\ep + \rh (t_1 - t_0) < d_0,$ we have
$\ds (\gm^{(j (r))}, \ph^{(0)} ) < d_0$
for all $r.$ Now let $j$ be arbitrary. Choose $r$ such that
$| j - j(r) | \leq n' - 1.$ Then
\beqr
\ds (\gm^{(j)}, \ph^{(0)} ) & \leq &
   \ds (\gm^{(j)}, \gm^{(j (r))}) + \ds (\gm^{(j (r))}, \ph^{(0)} )  \\
 &  \leq & (n' - 1) \left[ 91 \rh \left(\frac{t_1 - t_0}{n} \right)
                                       + \frac{90}{n'} \cdot d_0 \right]
        + d_0
 < 91 n' \rh \left(\frac{t_1 - t_0}{n} \right) + 91 d_0,
\eeqr
as desired.
\end{pff}

The following lemma is essentially the induction step in the main
part of the proof of the theorem of this section.
It is the analog of Lemma 5.3 and Theorem 5.4 of \cite{HR}.

\begin{lem}   
Let $X$ and $Y$ be compact metric spaces, with metric $d_Y$ on $Y.$
Let $A$ be a unital continuous field of \ca s over
$X \times [0, 1] \times Y,$ and suppose that there are
continuous unitary sections $u_1, \dots, u_m$ such that for each
$\xi \in X \times [0, 1] \times Y$ the elements
$u_1 (\xi), \dots, u_m (\xi)$ generate $A (\xi)$ as a \ca.
If $A$ is $(X, \rh)$-embeddable in $\OA{2},$ with respect to the metric
$d_{[0, 1] \times Y} ( (t_1, y_1), (t_2, y_2) ) =
              |t_2 - t_1| + d_Y (y_1, y_2),$
then $A$ is $(X \times [0, 1], 10 \rh)$-embeddable in $\OA{2},$ with
respect to the metric $d_Y.$
\end{lem}

\begin{pff}
Assuming the existence part of
$(X \times  [0, 1], 10 \rh)$-embeddability, the \aue\  part follows
directly from Lemma 5.4.
We therefore prove existence.

We can assume \wolog\  that $\rh ( r ) > 0$ for $r  > 0.$
Fix $y_0 \in Y;$ we construct an injective representation of
$A |_{X \times [0, 1] \times \{ y_0\} }$ in $\OA{2}.$
Choose positive integers $n_1,\, n_2, \dots$ such that $n_k \geq 360$
and the integers $N_k = n_1 n_2 \cdots n_k$ (with $N_0 = 1$) satisfy
\[
\rh \left( \frac{1}{N_{k + 1}} \right) \leq \frac{1}{2} \cdot
            \rh \left( \frac{1}{N_k} \right)
\]
for all $k.$ Set $n' = 360.$ 
Define $d_0 = 91 \rh (1),$ and inductively define
\[
d_{k + 1}
  = 91 \rh \left( \frac{1}{N_{k + 1}} \right) + \frac{90}{n'} \cdot d_k.
\]
Note that, by induction, we have
$\rh \left( \frac{1}{N_k} \right) \leq 2^{-k} \rh (1).$ We can now
inductively estimate $d_k.$ First, $d_0 \leq 2 \cdot 91 \cdot \rh (1).$
Furthermore, if $d_k \leq 2^{-k+1} \cdot 91 \cdot \rh (1),$ then
\[
d_{k + 1}
  = 91 \rh \left( \frac{1}{N_{k + 1}} \right) + \frac{1}{4} \cdot d_k
  \leq 2^{-k-1} \cdot 91 \cdot \rh (1) +
                \frac{1}{4} \cdot 2^{-k+1} \cdot 91 \cdot \rh (1)
 = 2^{-k}\cdot 91 \cdot \rh (1).
\]
Therefore
$d_k \leq 2^{-k+1} \cdot 91 \cdot \rh (1)
                         \leq 2^{-k} \cdot 182 \cdot \rh (1)$
for all $k.$

We now construct injective unital representations $\ph^{(t)}$ of
$A |_{X \times \{t\} \times \{ y_0\} }$ in $\OA{2}$ for $t$ in the set
\[
S = \{ j/N_k : k \geq 0, \,\, 0 \leq j \leq N_k\}.
\]
These representations are required to satisfy
\[
\ds (\ph^{( (j - 1)/N_k)}, \ph^{(j/N_k)} ) < d_k
\]
for $1 \leq j \leq N_k$ and
\[
\ds (\ph^{( j / N_k + r / N_{k + 1})}, \ph^{(j / N_k)} )
              < d_k + n' d_{k + 1}
\]
for $0 \leq j \leq N_k - 1$ and $0 \leq r \leq n_{k + 1}.$

The construction
is by induction on $k.$ By hypothesis, there are
injective unital representations $\ph^{(t)}$ of
$A |_{X \times \{t\} \times \{ y_0\} }$ in
$\OA{2}$ for $t = 0$ and $t = 1.$ Since $\rh (1) < d_0,$ the
\aue\  part of the hypothesis of $(X, \rh)$-embeddability implies
that there is a continuous function $v : X \to \OA{2}$
such that $\ds (\ph^{(0)}, v \ph^{(1)} v^*) < d_0.$ Replacing
$\ph^{(1)}$ by $v \ph^{(1)} v^*,$ we may assume that
$\ds (\ph^{(0)}, \ph^{(1)}) < d_0.$

For the induction step, assume that $\ph^{(j / N_k)}$ has been
constructed for $0 \leq j \leq N_k.$ Apply the previous
lemma to
each interval $[j/N_k, (j + 1)/N_k],$ and call the resulting
representations $\ph^{( j / N_k + r / N_{k + 1})}$
for $0 \leq r \leq n_{k + 1}.$ The estimates in the conclusion of
that lemma give
\[
\ds (\ph^{( (i - 1) /N_{k + 1})}, \ph^{(i/N_{k + 1})} )
  < 91 \rh \left( \frac{1}{N_{k + 1}} \right) + \frac{90}{n'} \cdot d_k
  = d_{k + 1}
\]
for $j n_{k + 1} < i \leq (j + 1) n_{k + 1}$ and
\[
\ds (\ph^{(i/N_{k + 1})}, \ph^{(j n_{k + 1}/N_{k + 1})} )
  < 91 n' \rh \left( \frac{1}{N_{k + 1}} \right) + 91 d_k
  = d_k + n' d_{k + 1}
\]
for $j n_{k + 1} \leq i \leq (j + 1) n_{k + 1},$ as desired. This
completes the inductive construction.

Note now that for $0 \leq r \leq n_{k + 1}$ we have
\beqr
\ds (\ph^{( j / N_k + r / N_{k + 1})}, \ph^{((j + 1) / N_k)} )
  & \leq & \ds (\ph^{( j / N_k + r / N_{k + 1})}, \ph^{(j / N_k)} )
       + \ds (\ph^{(j/N_k)}, \ph^{( (j + 1) / N_k)} )     \\
  & < & 2 d_k + n' d_{k + 1}.
\eeqr
Let $0 \leq j_1 < j_2 \leq N_{k + 1},$ and let $i_1$ and $i_2$
be the smallest and largest integers respectively that satisfy
\[
\frac{j_1}{N_{k + 1}} \leq \frac{i_1}{N_k}  \andeqn
\frac{i_2}{N_k} \leq \frac{j_2}{N_{k + 1}}.
\]
If $i_1 < i_2,$ then
\beqr
\lefteqn{
\ds (\ph^{( j_1 / N_{k + 1})}, \ph^{( j_2 / N_{k + 1})} )
    }   \\
  & \leq & \ds (\ph^{( j_1 / N_{k + 1})}, \ph^{( i_1 / N_k)} )
      + \ds (\ph^{( i_1 / N_k)}, \ph^{( i_2 / N_k)} )
      + \ds (\ph^{( i_2 / N_k)}, \ph^{( j_2 / N_{k + 1})} ) \\
  & < & 3 d_k + 2 n' d_{k + 1}
          + \ds (\ph^{( i_1 / N_k)}, \ph^{( i_2 / N_k)} ).
\eeqr
Otherwise, we have $i_1 = i_2$ or $i_1 = i_2 + 1.$ In either case,
\[
\ds (\ph^{( j_1 / N_{k + 1})}, \ph^{( j_2 / N_{k + 1})} )
 \leq \ds (\ph^{( j_1 / N_{k + 1})}, \ph^{( i_2 / N_k)} )
       + \ds (\ph^{( i_2 / N_k)}, \ph^{( j_2 / N_{k + 1})} )
  \leq 3 d_k + 2 n' d_{k + 1}.
\]
(We actually
get $2 d_k + 2 n' d_{k + 1}$ if it happens that $i_1 = i_2 + 1.$)
The second estimate necessarily holds whenever
\[
\frac{j_2}{N_{k + 1}} - \frac{j_1}{N_{k + 1}} < \frac{1}{N_k}.
\]
An
induction argument therefore shows that if $0 \leq t_1 \leq t_2 \leq 1$
are in $S$ and satisfy $t_2 - t_1 < \frac{1}{N_k},$ then
(using $d_k \leq 2^{-k} \cdot 182 \cdot \rh (1)$ and $n' = 360$)
\[
\ds (\ph^{(t_1)}, \ph^{(t_2)} ) <
   \sum_{s = k}^{\infty} ( 3 d_s + 2 n' d_{s + 1}) \leq
         (6 + 2 n') \cdot 182 \cdot \rh (1)
                   \sum_{s = k + 1}^{\infty} 2^{-s}
   \leq 133,000 \rh (1) \cdot 2^{-k}.
\]

Consider now, for each fixed $l,$ the function from $S$ to
$C(X, \OA{2})$ which sends $t$ to
$x \mapsto \ph^{(t)}_x (u_l (x, t, y_0)).$ The estimate of the previous
paragraph implies that this function is uniformly continuous, and
therefore extends by continuity to a function defined on all of
$[0, 1],$ whose value at $t$ we denote by $x \mapsto w_l (x, t).$
We now want to extend the map
$u_l (x_0, t_0, y_0) \mapsto w_l (x_0, t_0)$ to
a \hm\  $\ps_{(x_0, t_0)} : A (x_0, t_0, y_0) \to \OA{2}.$
Let $p$ be a polynomial in $2 m$ noncommuting variables, and suppose
\[
p (u_1 (x_0, t_0, y_0), u_1 (x_0, t_0, y_0)^*, \dots,
              u_m (x_0, t_0, y_0), u_m (x_0, t_0, y_0)^*) = 0.
\]
Then
\[
t \mapsto p (u_1 (x_0, t, y_0), u_1 (x_0, t, y_0)^*, \dots,
              u_m (x_0, t, y_0), u_m (x_0, t, y_0)^*)
\]
is a continuous section of
$A |_{ \{ x_0\} \times  [0, 1] \times \{y_0\}}$ which vanishes at
$(x_0, t_0, y_0).$ Considering a sequence $(t_k)$
in $S$ which converges to
$t_0,$ and using the fact that $\ps_{(x_0, t)}$ is the restriction
of a \hm\  for $t \in S,$ we see that
\beqr
\lefteqn{
p (w_1 (x_0, t_0), w_1 (x_0, t_0)^*, \dots, w_m (x_0, t_0),
             w_m (x_0, t_0)^* )  }   \\
 & = & \limi{k} \ps_{(x_0, t_k)} (
           p (u_1 (x_0, t_k, y_0), u_1 (x_0, t_k, y_0)^*, \dots,
              u_m (x_0, t_k, y_0), u_m (x_0, t_k, y_0)^*) ) = 0.
\eeqr
It follows that $u_l (x_0, t_0, y_0) \mapsto w_l (x_0, t_0)$
extends to a \hm\ from the $*$-subalgebra generated by
$u_1 (x_0, t_0, y_0), \dots, u_m (x_0, t_0, y_0)$ to $\OA{2}.$
A similar approximation argument shows that this \hm\  is
a contraction. Therefore it extends to the \ca\  generated by
these elements, which is $A (x_0, t_0, y_0)$ by hypothesis.
One checks directly that for a polynomial $p_{x, t}$
in $2 m$ noncommuting variables and with coefficients varying
continuously with $(x, t),$ the function
\beqr
(x, t) & \mapsto &
         \ps_{x, t} (p_{x, t} (u_1 (x, t, y_0), u_1 (x, t, y_0)^*,
               \dots, u_m (x, t, y_0), u_m (x, t, y_0)^*) )
                                    \\
 &  & \mbox{} = p_{x, t} (w_1 (x, t), w_1 (x, t)^*, \dots, w_m (x, t),
                                                 w_m (x, t)^* )
\eeqr
is continuous. A partition of unity argument shows that sections
of the form
\[
(x, t) \mapsto p_{x, t} (u_1 (x, t, y_0), u_1 (x, t, y_0)^*, \dots,
                      u_m (x, t, y_0), u_m (x, t, y_0)^*)
\]
are dense in the set of all continuous sections of
$A |_{ X \times  [0, 1] \times \{y_0\}},$ and a standard argument
now shows that $(x, t) \mapsto \ps_{(x, t)} (a (x, t))$
is continuous for any continuous section $a$ of
$A |_{ X \times  [0, 1] \times \{y_0\}}.$ Therefore $\ps$ is a
representation of $A |_{ X \times  [0, 1] \times \{y_0\}}$ in
$\OA{2}.$ The maps $\ps_{(x, t)}$ are injective for $t \in S$ by
construction, so $\ps$ is injective by Lemma 4.6.

This completes the proof of the existence part of
$(X \times  [0, 1], 10 \rh)$-embeddability.
\end{pff}

\begin{thm}   
Let $X \subset [0, 1]^n$ be a compact subset, and suppose there is
an open set $U \subset [0, 1]^n$ which contains $X$ and a continuous
retraction $f : U \to X$ such that $f |_X = \id_X.$ Let $A$
be a continuous field of \ca s over $X,$ and assume that the
section algebra $\Gm (A)$ is separable and exact. Then
$A$ has an injective representation in $\OA{2},$ which can be taken
unital if $A$ is unital. Moreover, if $A$ is unital then any
two injective unital representations $\ph^{(1)}$ and $\ph^{(2)}$ of $A$
are
\ayue\  in the sense that, for any finite set $F$ of sections and any
$\ep > 0,$ there is a continuous unitary function $u : X \to \OA{2}$
such that
$\| u (x) \ph^{(1)}_x (a (x)) u (x)^* - \ph^{(2)}_x (a (x)) \| < \ep$
for all $a \in F$ and $x \in X.$
\end{thm}

\begin{pff}
We do the existence part first.

Unitizing $A$ as in Lemma 4.2,
we reduce to the unital case.
(Exactness of the section algebra is preserved by Lemma 4.8.)
Theorem 8 of \cite{OZ} shows that $\Kt \Gm (A)$ is singly generated.
So $(\Kt \Gm (A))\unit$ has
a finite generating set, which we may take to consist of
selfadjoint sections $a_1, \dots, a_l$ satisfying $\| a_l \| < \pi.$

Let $B$ be the continuous field over $[0, 1]^n$ obtained by first
forming the unitized tensor product $(\Kt A)\unit$ following
Lemmas 4.1
and 4.2,
then
constructing $f^* ( (\Kt A)\unit)$ over $U$  following Lemma 4.3,
and finally extending over $[0, 1]^n$ with fiber $\C$ at points
not in $U$ as in Lemma 4.4.
Note that $B |_X \cong (\Kt A)\unit.$

We must show that $\Gm (B)$ is exact. First,
$\Gm (\Kt A) \cong \Kt \Gm (A)$ is exact by
Proposition 7.1 (iii) of \cite{Kr2}.
Then $\Gm ((\Kt A)\unit)$ is exact, by Lemma 4.8
again.
We now show that $\Gm (B)$ is exact by showing that the
maps $B (x_0) \to \Gm (B)/\{ b \in \Gm (B) : b (x_0) = 0 \}$ are locally
liftable (Condition (2) of Theorem 4.7).
If $x_0 \not\in U,$ then $B (x_0) = \C,$ and this map lifts to
a \hm\  to $\Gm (B).$  If $x_0 \in U,$ let $E \subset B (x_0)$ be
a \fd\  operator system. Note that
$B (x_0) = ((\Kt A)\unit) (f (x_0)),$ so there is (by (1) implies (2)
of Theorem 4.7)
a \ucp\  map
$T_0 : E \to \Gm ((\Kt A)\unit)$ which lifts the map
\[
E \longrightarrow \Gm ((\Kt A)\unit) / \{ a \in \Gm ((\Kt A)\unit) :
                                a (f (x_0)) = 0 \}.
\]
Choose any state $\om$ on $E,$ and choose a \ct\  function
$h_0 : [0, 1]^n \to [0, 1]$ such that $h_0$ vanishes outside
$U$ and $h_0 (x_0) = 1.$ Then define $T : E \to \Gm (B)$ by
\[
T (b) (x) =
     h_0 (x) T_0 (b) (f (x)) + (1 - h_0 (x)) \om (b) \cdot 1_{B (x)}.
\]
This map lifts $E \to \Gm (B)/\{ b \in \Gm (B) : b (x_0) = 0 \},$
and is readily checked to be \uacp.
This completes the proof that $\Gm (B)$ is exact.

Choose a continuous function
$h : [0, 1]^n \to [0, 1]$ such that $h (x) > 0$ exactly when
$x \in U.$ Then the functions $x \mapsto h (x) a_l ( f (x))$ (taken
to be $0$ for $x \not\in U$) are continuous sections of $B.$
Together with $1,$ their values at $x$ generate $B (x)$ for all $x.$
Therefore (since $\| a_l \| < \pi$)
the sections $u_l (x) = \exp (i h (x) a_l ( f (x))),$ for
$1 \leq l \leq m,$ also generate each $B (x).$

For $0 \leq k \leq n$
write $[0, 1]^n$ as $Y_k \times Z_k$ with $Y_k = [0, 1]^k$ and
$Z_k = [0, 1]^{n - k}.$ From Remark 4.16
it follows that
$B$ is $(Y_0, \rh)$-embeddable in $\OA{2}$  for a suitable $\rh,$
with respect to the metric $d (x, y) = \sum_{j = 1}^n |y_j - x_j|.$
Using induction and Lemma 5.6,
we find that
$B$ is $(Y_k, 10^k \rh)$-embeddable in $\OA{2}$  for all $k.$ In
particular, taking $k = n,$ we see that
$B$ has an injective unital representation in $\OA{2}.$

Restriction to $X$ gives an injective unital representation $\ph$ of
$(\Kt A)\unit$ in $\OA{2}.$ Let $e \in K$ be a rank one \pj, and
define a continuous section $p$ of $(\Kt A)\unit$ by
$p (x) = e \otimes 1_{A (x)}.$ Then the function
$x \mapsto \ph_x (p(x))$ is a \pj\  in $C(X, \OA{2}).$ It follows
from \cite{Zh1} that the set $P$ of projections in $\OA{2}$ other
than $0$ and $1$ is weakly contractible (all homotopy groups trivial).
Since it is homotopy equivalent to an open subset of a Banach space,
it is contractible. Therefore the projection $x \mapsto \ph_x (p(x))$
is homotopic to
a constant \pj\  $x \mapsto p_0$ for some $p_0 \in \OA{2}.$
It follows that there is a unitary $v \in C(X, \OA{2})$ such that
$v (x) \ph_x (p(x)) v(x)^* = p_0$ for all $x \in X.$
Then $a \mapsto v (x) \ph_x (e \otimes a) v(x)^*$ is an
injective unital \hm\  from $A (x)$ to $p_0 \OA{2} p_0,$
and the family of all these \hm s is
an injective unital representation of $A$ in $p_0 \OA{2} p_0.$
Since $p_0 \OA{2} p_0 \cong \OA{2},$ the existence part is proved.

Now we do the approximate
uniqueness part. Let the notation be as in the existence
part. Choose a \hm\  $\mu_0 : \Kt \OA{2} \to \OA{2}$ which is an
isomorphism onto a (nonunital) hereditary subalgebra of $\OA{2},$
and let $\mu : (\Kt \OA{2})\unit \to \OA{2}$ be the unital extension.
Then the definition
$\tilde{\ph}^{(i)} = \mu \circ (\id_K \otimes \ph^{(i)}_x)$ gives
two injective unital representations of $(\Kt A)\unit.$
Extend these to unital
injective representations of $B$ by setting
$\ps^{(i)}_x = \tilde{\ph}^{(i)}_{f (x)}$ for $x \in U$ and
$\ps^{(i)}_x (\ld) = \ld \cdot 1$ for $x \not\in U.$ In the existence
part of the proof, we saw that $B$ is $(Y_n, 10^n \rh)$-embeddable in
$\OA{2}.$ It follows that $\ps^{(1)}$ is \ayue\  to $\ps^{(2)}.$
Restricting to $X,$ we see that $\tilde{\ph}^{(1)}$ is \ayue\  to
$\tilde{\ph}^{(2)}.$

Let $F \in \Gm (A)$ be finite, and assume
all elements of $F$ have norm at most $1.$ Regard $F$
as a subset of $\Gm ((\Kt A))\unit$ via $a \to p \otimes a$
(with $p$ as above).  Choose $\dt > 0$ such
that $\dt < \frac{\ep}{3},$ and also small enough that if projections
$q_1$ and $q_2$ satisfy $\| q_1 - q_2 \| < \dt,$ then there is a unitary
$z$ with $z q_1 z^* = q_2$ and $\| z - 1 \| < \frac{\ep}{3}.$
Apply the definition of \aue\ to $G = F \cup \{p\},$ using $\dt$
for $\ep.$ Call the resulting unitary $w.$ Find $z$ as above with
$q_1 = w p w^*$ and $q_2 = p.$  Then $u = p z w p$ is a unitary function
with values in $\OA{2} = \C e \otimes \OA{2} \subset (\Kt \OA{2})\unit$
which satisfies
$\| u (x) \ph^{(1)}_x (a (x)) u (x)^* - \ph^{(2)}_x (a (x)) \| < \ep$
for all $a \in F$ and $x \in X.$
\end{pff}

We now want to extend this theorem to finite CW complexes. For this,
we need a relative version. The following lemma is easy to prove and
suffices.

\begin{lem}
Let $D^n$ be the closed unit ball in $\R^n,$ with boundary $S^{n - 1}.$
Let $A$ be a continuous field of \ca s over $D^n,$ such that the
section algebra $\Gm (A)$ is separable and exact.

(1) Let $\ph$ be an injective representation of $A |_{S^{n - 1}}$
in $\OA{2}.$ Then there exists an injective representation $\ps$ of $A$
in $\OA{2}$ such that $\ps |_{S^{n - 1}} = \ph.$
If $A$ and $\ph$ are unital, so is $\ps.$

(2) Assume $A$ is unital, and let $\ph^{(1)}$ and $\ph^{(2)}$
be injective unital representations of $A$
in $\OA{2}.$ Let $u_1, \dots, u_n$ be unitary sections of
$A,$ and let $z_0 : S^{n - 1} \to U (\OA{2})$ be continuous and
satisfy
$\ds (z_0 (\ph^{(1)} |_{S^{n - 1}}) z_0^*, \ph^{(2)} |_{S^{n - 1}})
           < \dt.$
Then there is a continuous unitary $z : D^n \to U (\OA{2})$ such that
$z |_{S^{n - 1}} = z_0$ and $\ds (z \ph^{(1)} z^*, \ph^{(2)}) < 10\dt.$
\end{lem}

\begin{pff}
It is convenient to define
\[
S (r) = \{ x \in D^n : \| x \| = r \} \andeqn
S (r_1, r_2) = \{ x \in D^n : r_1 \leq \| x \| \leq r_2 \}.
\]
We identify $S (r)$ with $S^{n - 1}$ in the obvious way. We further
identify $S (r_1, r_2)$ with $S^{n - 1} \times [0, 1]$ by starting
at the outside edge: $S (r_2)$ goes to $S^{n - 1} \times \{0\}.$

We first prove (1). Unitizing $A$ (and correspondingly $\ph$) as in
the proof of the previous theorem, we may assume $A$ and $\ph$ are
unital.  Let $\ps^{(0)}$ be an arbitrary
injective unital representation of $A$ in $\OA{2}$ (from the previous
theorem).
Choose a sequence $u_1, u_2, \dots$ of unitary sections of $A$ such
that $u_1 (x), u_2 (x), \dots$ generate $A (x)$ for every $x \in D^n.$
Write $\ds^{(n)}$ for the sectional distance with respect to
$u_1, \dots, u_n.$ Choose $0 < r_1 < r_2 < \cdots < 1,$ with
$r_n \to 1,$ such that
$\ds^{(n)} (\ps^{(0)} |_{S (r)}, \ps^{(0)} |_{S (1)}) \leq 2^{-n}$
for $r_n \leq r \leq 1.$ By the previous theorem,
$\ps^{(0)} |_{S (1)}$ is \ayue\  to $\ph,$ so
there are continuous functions $w_{n, r} : S (r) \to U (\OA{2}),$
for $r_n \leq r \leq r_{n + 1},$ such that
$\ds^{(n)} (w_{n, r} (\ps^{(0)} |_{S (r)}) w_{n, r}^*, \ph)
                     \leq 2 \cdot 2^{-n}.$
We may take $w_{1, r_1} = 1.$ By Lemma 5.4, there exist
$z_n : S (r_n, r_{n + 1})  \to U (\OA{2})$ such that
\[
z_n |_{S (r_n)} = w_{n, r_n}, \,\,\,\,\,\,\,
z_n |_{S (r_{n + 1})} = w_{{n + 1}, r_{n + 1}}, \andeqn
\ds^{(n)} (z_n (x) \ps^{(0)}_x z_n (x)^*, \ph_{x / \| x \|} )
                       \leq 20 \cdot 2^{-n}.
\]
Define a continuous function $z$ from the interior of $D^n$ to
$U (\OA{2})$ by setting $z = z_n$ on $S (r_n, r_{n + 1})$ and
$z (x) = 1$  for $\|x\| < r_1.$ Then set $\ps_x = \ph_x$ for
$\| x \| = 1$ and
$\ps_x = z (x) \ps^{(0)}_x z (x)^*$ for  $\| x \| < 1.$ This defines the
required $\ps.$

Now we prove (2). The previous theorem provides a continuous
$w : D^n \to U (\OA{2})$ such that
$\ds (w \ph^{(1)} w^*, \ph^{(2)}) < \dt.$
Use Lemma 5.4 to find $z : S (\frac{1}{2}, 1) \to U (\OA{2})$ such that
\[
z |_{S^{n - 1}} = z_0, \,\,\,\,\,\,\, z |_{S (1/2)} = w |_{S (1/2)},
   \andeqn
\ds (z (\ph^{(1)} |_{S (1/2, 1)}) z^*,
           \ph^{(2)} |_{S (1/2, 1)})) < 10 \dt.
\]
Then define $z (x) = w (x)$ for $\| x \| < \frac{1}{2}.$
\end{pff}

\begin{thm}
Let $X$ be a finite CW complex. Let $A$
be a continuous field of \ca s over $X,$ such that
$\Gm (A)$ is separable and exact. Then
$A$ has an injective representation in $\OA{2},$ which can be taken
unital if $A$ is unital. Moreover, if $A$ is unital then any
two injective unital representations $\ph^{(1)}$ and $\ph^{(2)}$ of $A$
are \ayue\  in the sense of Theorem 5.7.
\end{thm}

\begin{pff}
Both parts are proved by induction over the cells, and both are
immediate from previous theorems for a zero dimensional finite
CW complex.

For the existence induction step, assume the theorem is known for
$X,$ and let $Y = X \cup_f D_n,$ where $f : S^{n - 1} \to X$ is
the attaching map. Let $g : D_n \to Y$ be the map
extending $f.$ Let $\ph$ be a (unital)
injective representation of $A |_X$ in $\OA{2}.$ Then
$x \mapsto \ph_{f (x)}$ is a (unital) injective representation of
$g^* (A) |_{S^{n - 1}}$ in $\OA{2}.$ Now $g^* (A)$ is exact
(one checks condition (2) of Theorem 4.7; it is easier than in
the proof of Theorem 5.7) and has separable section algebra. The
previous lemma therefore provides a
(unital) injective representation of $g^* (A)$ in $\OA{2}$ which
extends $x \mapsto \ph_{f (x)}.$ Use it to extend $\ph$ to a
representation of $A.$

For the approximate
uniqueness result, let a finite set $F$ of sections and
$\ep > 0$ be given. \Wolog\  the sections in $F$ are all unitary.
Let $N$ be the number of cells of strictly positive dimension, and let
$X_0$ be the zero skeleton. Given $\ph^{(1)}$ and $\ph^{(2)},$
choose a unitary $v_0 : X_0 \to \OA{2}$ such that
$\ds (v_0 (\ph^{(1)} |_{X_0}) v_0^*, \ph^{(2)} |_{X_0}) < 10^{-N} \ep.$
Use the uniqueness part of the previous lemma, in the same way the
existence part was used in the previous paragraph, to extend $v_0$
cell by cell. If $X_n$ is the subcomplex obtained by adding $n$
cells and $v_n$ is the unitary defined on it, we will have
$\ds (v_n (\ph^{(1)} |_{X_n}) v_n^*, \ph^{(2)} |_{X_n})
                            < 10^{-N + n} \ep.$
\end{pff}

The advantage of the methods of this section is that they give control
over the ``smoothness'' of the images under the representation
of the generating sections. Here is the one dimensional version,
which is easy to prove.

\begin{thm}   
Let $A$ be a unital  continuous field over $[0, 1]$ such
that $\Gm (A)$ is exact, and let $u_1, \dots, u_m \in \Gm (A)$ be
unitary sections of $A$ such that, for each $x \in X,$ the
elements $u_1 (x), \dots, u_m (x)$ generate $A (x).$ Suppose the
function $\rh_0$ of Proposition 4.11
(using
$u_1, \dots, u_m$ in place of $a_1, \dots, a_m$)
is $\Lip^{\af}$ for some
$\af \in (0, 1],$ that is, there is a constant $C_0$ such that
$\rh_0 (x, y) \leq C_0 |x - y |^{\af}$ for all $x,$ $y \in [0, 1].$
Then there is a unital injective representation $\ph$ of
$A$ in $\OA{2}$ such that the the functions $x \mapsto \ph_x (u_l (x))$
are $\Lip^{\af / 2},$ that is, there is a constant $C$ such that
$\| \ph_x (u_l (x)) - \ph_y (u_l (y)) \| \leq C |x - y |^{\af / 2}$
for all $x,$ $y \in [0, 1].$ Moreover, $C$ depends only on $C_0$ and
$\af.$
\end{thm}

\begin{pff}
By Proposition 4.13,
the function $\rh$ defined there
satisfies
\[
\rh (x, y) \leq 11 C_0^{1/2} |x - y |^{\af/2}.
\]
That is,
taking $X_0$ to be a one point space, $A$ is $(X_0, \rh)$-embeddable in
$\OA{2}$ with $\rh (t) = 11 C_0^{1/2} t^{\af/2}.$ (See Definition 4.15
and Remark 4.16.)

We now follow the proof of Lemma 5.6,
with the spaces
$X$ and $Y$ there both one point spaces, and making suitable
minor modifications. Set $C_1 = 11 C_0^{1/2},$ and set
\[
\bt = \frac{1}{1 - \af / 2} \in [1, 2].
\]
Choose some integer $n$ with $181^{\bt}\leq n \leq 181^{\bt} + 1.$
Choose the numbers of the proof of Lemma 5.6
to be
$n_1 = n_2 = \cdots = n' = n,$ so that $N_k = n^k.$ Take
$d_0 = 181 C_1,$ and define
$d_{k + 1} = 91 \rh (n^{- (k + 1)}) + 90 n^{-1} d_k$ inductively,
as in the proof of Lemma 5.6.
We prove by induction that
$d_k \leq 181 C_1 n^{- k \af / 2}.$ The one slightly nontrivial
step is the observation that
\[
181 n^{ - (1 - \af / 2)}  \leq 181^{1 - \bt (1 - \af / 2)} = 1.
\]

Following the procedure of the proof of Lemma 5.6,
we now
obtain, for $n^{ - (k + 1)} < t_2 - t_1 \leq n^{-k},$ the estimate
\beqr
\ds (\ph^{(t_1)}, \ph^{(t_2)} ) & < &
   \sum_{s = k}^{\infty} ( 3 d_s + 2 n d_{s + 1})
   \leq  3 \cdot 181 C_1 \cdot
                       \frac{n^{- k \af / 2}}{1 - n^{- \af / 2}} +
      2 n \cdot 181 C_1 \cdot
                       \frac{n^{- (k + 1) \af / 2}}{1 - n^{- \af / 2}}
                        \\
  & = & \frac{181 C_1 (3 n^{\af / 2} + 2 n)}{1 - n^{- \af / 2}}
                   \cdot n^{- (k + 1) \af / 2}
  \leq  \frac{181 C_1 \left(5 (181^{\bt} + 1) \right)}
            {1 - 181^{ - \af \bt / 2}}
                   \cdot (t_2 - t_1)^{\af / 2}.
\eeqr
The rest of the proof goes through as it stands, and we obtain in the
end
\[
\| \ps_{t_1} (u_l (t_1)) - \ps_{t_2} (u_l (t_2)) \| \leq
     M (\af) C_0^{1/2} | t_1 - t_2 |^{\af / 2},
\]
with
\[
M (\af) = 11 \cdot \frac{181 \cdot 5 (181^{\bt} + 1)}{1 - 181^{ - \af \bt / 2}}.
\]
This proves the theorem with $C = M (\af) C_0^{1/2}.$
\end{pff}

\begin{rmk}
With more care, the choice of $M (\af)$ in the proof of Theorem 5.10
can be improved considerably.
To keep down the sizes of some of the numbers in the next section, we
describe how to get $C \leq 330,000 C_0^{1/2}$ in the case $\af = 1.$
First, taking $n' = n$ in the proof of Lemma 5.5 allows considerable
simplification (this is the case done in \cite{HR}) and improvement of
the conclusion to
\[
\ds (\gm^{(j - 1)}, \gm^{(j)} ) <
    46 \cdot \rh \left(\frac{t_1 - t_0}{n} \right) + \frac{45}{n} \cdot d_0.
\]
Let $C_1$ be as in the proof of Theorem 5.10, but take $n = n' = 90^2$
(as in \cite{HR}).
Take $d_0 = 46 C_1$ and
$d_{k + 1} = 46 \rh (n^{- (k + 1)} ) + 45 n^{-1} d_k.$
Then $d_k \leq 92 C_1 n^{- k/2}.$
For $t_1$ and $t_2$ of the form $j_1 / n^{k_1}$
and $j_2 / n^{k_2},$ estimate $\ds (\ph^{(t_1)}, \ph^{(t_2)} )$ by
the more careful method in the proof of Theorem 5.4 of \cite{HR}.
One obtains
\[
\ds (\ph^{(t_1)}, \ph^{(t_2)} )
    \leq 320 \cdot 92 \cdot C_1 \cdot |t_2 - t_1|^{1/2}
    \leq 30,000 C_1 |t_2 - t_1|^{1/2}.
\]
\end{rmk}

\section{The field of rotation algebras}

In this section, we apply the results of the previous section
specifically to the continuous field of rotation algebras. In
Theorem 5.4 and Corollary 5.5 of \cite{HR}, Haagerup and R\o rdam
construct a $\Lip^{1/2}$ (with respect to the sections given by the
standard generators) representation of the field of rotation
algebras in $L (H)$ for a separable infinite dimensional Hilbert space
$H.$
In this section, we produce a $\Lip^{1/2}$ representation in $\OA{2}.$

By combining the estimates of  Haagerup and R\o rdam (Theorem 4.9 (1)
of \cite{HR}), the explicit estimate in Lemma 1.8,
an easy
computation to show that the constant $M$ there is $1,$ and
Lemma 1.10,
one finds that the function $\rh_0$ of Proposition 4.11
satisfies
$\rh_0 (\te_1, \te_2) \leq 480 \cdot | \te_1 - \te_2 |^{1/2}.$
A slight modification of Theorem 5.10
(to use the circle
instead of $[0, 1]$) then gives a $\Lip^{1/4}$ representation in
$\OA{2}.$

We improve this procedure by estimating $\rh_0 (\te_1, \te_2)$
directly. By explicitly estimating the completely bounded norms of
certain linear maps between finite dimensional operator spaces in the
rotation algebras, we prove an inequality of the form
$\rh_0 (\te_1, \te_2) \leq C_0 | \te_1 - \te_2 |$ (not just
$C_0 | \te_1 - \te_2 |^{1/2}$) for some constant $C_0.$ This implies
the existence of a $\Lip^{1/2}$ representation in $\OA{2}.$
We therefore have an alternate proof (with different constants) of
Theorems 4.9 and 5.4 of \cite{HR}. This proof makes no use of
unbounded operators or canonical commutation relations.

We begin by establishing notation.

\begin{ntn}   
For $\te \in \R$ let $A (\te)$ be the rotation \ca, the universal
\ca\   generated by unitaries $u (\te)$ and $v (\te)$ satisfying
$u (\te) v (\te) = \exp (2 \pi i \te) v (\te) u (\te).$
By Corollary 3.6 of \cite{Rf}, the rotation algebras are the fibers of
a continuous field over the circle $S^1,$ which we think of as
$[0, 1]$ with the endpoints identified, or as $\R / \Z.$
Moreover, $u$ and $v$ are continuous sections.
(The fact that the rotation algebras form a continuous field in this
manner was known to Elliott and others before \cite{Rf}.)
Let $E (\te) \subset A (\te)$ be the operator space
$E (\te) = \spn (1, u (\te), u (\te)^*, v (\te), v (\te)^*).$

Further, for each $\te \in \R,$ fix a unital embedding
$\io_{\te} : A (\te) \to \OA{2},$ and use it to regard $A (\te)$ as a
unital subalgebra of $\OA{2}.$ Define
\[
\rh_0 (\te_1, \te_2) =
    \inf_T ( \max ( \| T (u (\te_1)) - \io_{\te_2} (u (\te_2)) \|,
          \| T (v (\te_1)) - \io_{\te_2} (v (\te_2)) \| ) ),
\]
where the infemum is taken over all \ucp\  maps
$T : A (\te_1) \to \OA{2}.$ By Proposition 4.11 (1),
this function does not depend on
the embeddings used. (In the case at hand, since all algebras involved
are nuclear, Remark 4.12
implies that one gets the same function by taking the
infemum over all \ucp\  maps $T : A (\te_1) \to A (\te_2).$)
\end{ntn}

The following two lemmas constitute an analog of Proposition 4.5 of
\cite{HR}.

\begin{lem}   
The function $\rh_0$ is continuous and translation invariant, that is,
$\rh_0 (\te_1 + \te, \te_2 + \te) = \rh_0 (\te_1, \te_2)$
for all $\te_1,$ $\te_2,$ $\te \in \R.$
\end{lem}

\begin{pff}
The function $\rh_0$ is continuous by Proposition 4.11 (2).
For translation invariance, it suffices to prove that
$\rh_0 (\te_1 + \te, \te_2 + \te) \leq \rh_0 (\te_1, \te_2).$
By continuity, we may restrict
to $\te_1,$ $\te_2$ irrational and $\te$ rational.

Let $\ep > 0,$ and let $T : A (\te_1) \to \OA{2}$ satisfy
\[
\max (\| T (u (\te_1)) - \io_{\te_2} (u (\te_2)) \|,
                        \| T (v (\te_1)) - \io_{\te_2} (v (\te_2)) \|)
       < \rh_0 (\te_1, \te_2) + \frac{\ep}{2}.
\]
Set $u_1 = u (\te_1) \otimes u (\te)$ and
$v_1 = v (\te_1) \otimes v (\te),$ which are unitaries in
$A (\te_1) \otimes A (\te)$ satisfying
$u_1 v_1 = \exp ( 2 \pi i (\te_1 + \te)) v_1 u_1.$ Since
$\te_1 + \te$ is irrational, this gives a unital embedding $\ld$ of
$A (\te_1 + \te)$ in $A (\te_1) \otimes A (\te).$ Similarly,
the unitaries $u_2 = u (\te_2) \otimes u (\te)$ and
$v_2 = v (\te_2) \otimes v (\te)$ give a unital embedding $\ps$ of 
$A (\te_2 + \te)$ in $A (\te_2) \otimes A (\te).$
We will estimate $\rh_0 (\te_1 + \te, \te_2 + \te)$ by considering
the \ucp\  map $T \otimes \id_{A (\te)}$ followed by a suitable
embedding in $\OA{2}.$

Let $\mu : \OT{\OA{2}} \to \OA{2}$ be an isomorphism
(from Theorem   0.8),
and let
\[
\ph = \mu \circ (\io_{\te_2} \otimes \io_{\te}) \circ \ps :
             A (\te_2 + \te) \to \OA{2}.
\]
Then $\ph$ is \ayue\  to $\io_{\te_2 + \te}$ by Theorem 1.15,
so there is a unitary $w \in \OA{2}$ such that
$\| w \ph (a) w^* - \io_{\te_2 + \te} (a) \| \leq \half \ep \| a \|$
for $a \in E (\te_2 + \te).$

Now define $S : A (\te_1 + \te) \to \OA{2}$ by
\[
S_0 = \mu \circ (T \otimes \io_{\te}) \circ \ld  \andeqn
S (a) = w S_0 (a) w^*.
\]
This map is \uacp\  because $T \otimes \id_{A (\te)}$ is
(by Proposition IV.4.23 (i) of \cite{Tk}).
Furthermore, it is easy to check that
\beqr
\lefteqn{
\| S (u (\te_1 + \te)) - \io_{\te_2 + \te} (u (\te_2 + \te)) \|
     }       \\
 & \leq &
\| S_0 (u (\te_1 + \te)) - \ph (u (\te_2 + \te)) \| +
 \| w \ph (u (\te_2 + \te)) w^* - \io_{\te_2 + \te} (u (\te_2 + \te)) \|
                  \\
  &  < & \| T (u (\te_1)) \otimes u (\te) -
                   \io_{\te_2} (u (\te_2)) \otimes u (\te) \|
             + \frac{\ep}{2}
    <  \rh_0 (\te_1, \te_2) + \frac{\ep}{2} + \frac{\ep}{2}
        = \rh_0 (\te_1, \te_2) + \ep,
\eeqr
and similarly for $v$ in place of $u.$ This shows that
$\rh_0 (\te_1 + \te, \te_2 + \te) < \rh_0 (\te_1, \te_2) + \ep.$
\end{pff}

\begin{lem}   
The function $\rh_0$ is a continuous pseudometric on $\R.$
\end{lem}

\begin{pff}
Continuity was shown in the previous lemma.
The triangle inequality is Proposition 4.11 (3).
It remains to prove symmetry. It suffices to prove that
$\rh_0 (0, - \te) = \rh_0 (0, \te).$ Indeed, the previous lemma then
gives
\[
\rh_0 (\te_2, \te_1) = \rh_0 (0, \te_1 - \te_2) =
         \rh_0 (0, \te_2 - \te_1) = \rh_0 (\te_1, \te_2)
\]
for any $\te_1$ and $\te_2.$

We may clearly assume
$\te \neq 0.$ Let $\ph_{\te} : A (\te) \to A (- \te)$ be the
isomorphism given by
\[
\ph_{\te} (u (\te)) = v (\te) \andeqn \ph_{\te} (v (\te)) = u (\te).
\]
Since $\rh (0, \te)$ does not depend on the choice of embedding
$\io_{\te}$ of $A (\te)$ in $\OA{2},$ we may assume that
$\io_{- \te} = \io_{\te} \circ \ph_{\te}^{-1},$ that is, that
in $\OA{2}$ we have
\[
\io_{ - \te} (u ( - \te)) = \io_{\te} (v (\te)) \andeqn
\io_{ - \te} (v ( - \te)) = \io_{\te} (u (\te)).
\]
Given a \ucp\  map $T : A (0) \to \OA{2},$
the map $T \circ \ph_0 : A (0) \to \OA{2}$ satisfies
\[
\| (T \circ \ph_0) (u (0)) - \io_{ - \te} (u ( - \te)) \| =
    \| T (v (0)) -  \io_{\te} (v (\te)) \|
\]
and
\[
\| (T \circ \ph_0) (v (0)) - \io_{ - \te} (v ( - \te)) \| =
         \| T (u (0)) -  \io_{\te} (u (\te)) \|.
\]
Taking the infemum over all $T,$ we find that
$\rh_0 (0, - \te) \leq \rh_0 (0, \te).$ Replacing $\te$ by $- \te,$
we get equality.
\end{pff}

We now prepare to estimate $\rh_0 (0, \te).$

\begin{lem}   
Let $H$ be a Hilbert space, and let $\xi,$ $\et \in H$ with
$\| \xi \| = \| \et \| = 1.$ Then
$\Real ( \langle \xi, \et \rangle ) = 1 - \half \| \xi - \et \|^2.$
\end{lem}

\begin{pff}
Write $\| \xi - \et \|^2 = \langle \xi - \et, \xi - \et \rangle$ and
calculate.
\end{pff}

\begin{lem}   
Let $E_0 (\te) = \spn (1, u (\te), v (\te)) \subset A (\te).$
Define $T_{\te} : E_0 (\te) \to E_0 (0)$ by
\[
T_{\te} (1) = 1,  \,\,\,\,\,\,
T_{\te} (u (\te)) = u (0), \andeqn T_{\te} (v (\te)) = v (0).
\]
If $\te$ is a rational number of the form
$m^2 / (2 n + 1)^2$ with $m,$ $n \in \N$ (not necessarily in
lowest terms), and if $\te < 2/25,$ then
$\| T_{\te} \|_{\cb} \leq (1 - \xfrac{25}{2} \te)^{- 1/2}.$
\end{lem}

\begin{pff}
We start by making several reductions. First, note that if $E$ is an
operator space, then a bounded
linear functional $\om : E \to \C$ is completely bounded and
satisfies $\|\om \|_{\cb} = \|\om \|.$
(See, for example, Proposition 3.7 of \cite{Pl}.)
Therefore, using $A (0) = C (S^1 \times S^1),$ we have
\[
\| T_{\te} \|_{\cb} =
   \sup_{x \in S^1 \times S^1} \| \ev_x \circ T_{\te} \|_{\cb} =
  \sup_{x \in S^1 \times S^1} \| \ev_x \circ T_{\te} \| = \| T_{\te} \|.
\]
It therefore suffices to show that
$\| T_{\te} \| \leq (1 - \xfrac{25}{2} \te)^{- 1/2}.$

For the next step, it is convenient to use the (nonstandard)
convention $\sgn (0) = 1.$ Our problem is equivalent to showing that
for $\af,$ $\bt,$ $\gm \in \C$ we have:
\[
\| \af \cdot 1 + \bt u (0) + \gm v (0) \| = 1
  \,\,\,\,\,\, {\mathrm{implies}} \,\,\,\,\,\,
    \| \af \cdot 1 + \bt u (\te) + \gm v (\te) \| 
                   \geq \left(1 - \xfrac{25}{2} \te \right)^{1/2}.
\]
Multiplying the two expressions inside the norm signs by
$\overline{\sgn (\af)},$ we see that it suffices to prove this with
$\af \geq 0.$ Next, note that clearly
$\| \af \cdot 1 + \bt u (0) + \gm v (0) \| \leq \af + | \bt| + |\gm|,$
while the reverse inequality follows by considering the point
$x = \left(\overline{\sgn (\bt)}, \overline{\sgn (\gm)} \right)
                     \in S^1 \times S^1.$
So it suffices to show that
for $\af,$ $\bt,$ $\gm \in \C$ with $\af \geq 0$ we have
\[
\af + | \bt| + |\gm| = 1
  \,\,\,\,\,\, {\mathrm{implies}} \,\,\,\,\,\,
    \| \af \cdot 1 + \bt u (\te) + \gm v (\te) \|
                   \geq \left(1 - \xfrac{25}{2} \te\right)^{1/2}.
\]

Now let $q = (2 n + 1)^2,$ and define $\zt = \exp (2 \pi i /q),$ a
primitive $q$-th root of $1.$ Define unitaries $y_0,$ $z_0 \in M_q$ by
\[
y_0 = \diag (1, \zt, \dots, \zt^{q - 1})
    \andeqn
z_0 = \left( \begin{array}{cccccc}
            0     &    0   &    0   &   \cdots &    0   &    1     \\
            1     &    0   &    0   &   \cdots &    0   &    0     \\
            0     &    1   &    0   &   \cdots &    0   &    0     \\
           \vdots & \vdots & \vdots &   \ddots & \vdots & \vdots   \\
            0     &    0   &    0   &   \cdots &    0   &    0     \\
            0     &    0   &    0   &   \cdots &    1   &    0     \\
    \end{array} \right).
\]
Then $y_0 z_0 = \zt z_0 y_0.$ Set $y = y_0^m$ and $z = z_0^m,$
so that $y z = \zt^{m^2} z y = \exp (2 \pi i \te ) z y.$
Thus there is a unital \hm\  $\ph : A (\te) \to M_q$ given by
$\ph ( u (\te)) = \overline{\sgn (\bt)} y$ and
$\ph ( v (\te)) = \overline{\sgn (\gm)} z.$
We have
\[
\ph (\af \cdot 1 + \bt u (\te) + \gm v (\te) ) =
      \af \cdot 1 + |\bt| y + |\gm| z.
\]
Note that we can write $\te = m^2 / (2 n + 1)^2$ with $n$ arbitrarily
large. Since $\| \ph \| \leq 1,$ it therefore suffices to show that
for $\af,$ $\bt,$ $\gm \in \C$ with $\af \geq 0$ and for all large
enough $n,$ we have
\[
\af + | \bt| + |\gm| = 1
  \,\,\,\,\,\, {\mathrm{implies}} \,\,\,\,\,\,
    \| \af \cdot 1 + |\bt| y + |\gm| z \|
                   \geq \left(1 - \xfrac{25}{2} \te\right)^{1/2}.
\]
Equivalently, we assume that $\af,$ $\bt,$ $\gm \geq 0,$ and show
that
\[
\af +  \bt + \gm = 1
  \,\,\,\,\,\, {\mathrm{implies}} \,\,\,\,\,\,
    \| \af \cdot 1 + \bt y + \gm z \| \geq
                          \left(1 - \xfrac{25}{2} \te\right)^{1/2}
\]
for large $n.$

Define
\[
\xi_0 = \left(1, \frac{n - 1}{n}, \dots, \frac{2}{n}, \frac{1}{n}, 0, 0,
      \dots, 0, \frac{1}{n}, \frac{2}{n}, \dots, \frac{n - 1}{n} \right)
  \in \C^q \andeqn
\xi = \frac{1}{ \| \xi_0 \|} \xi_0.
\]
Using the formula
\[
1^2 + 2^2 + \cdots + n^2 = \xfrac{1}{6} n (n + 1) (2 n + 1),
\]
we can calculate
\[
\| \xi_0 \|^2 = 1 + 2 \sum_{k = 1}^{n - 1} \left( \frac{k}{n} \right)^2
     = 1 + \xfrac{1}{3} n^{-1} (n + 1) (2 n + 1) \geq \xfrac{2}{3} n.
\]
Further, using
\[
| 1 - \zt^{k m} | \leq k m | 1 - \zt | \leq
                  \frac{2 \pi k m}{(2 n + 1)^2}
\]
and $(n - k) k \leq n^2 / 4$ for $0 \leq k \leq n,$ we obtain
\beqr
\| \xi_0 - y\xi_0 \|^2 &  = &
   2 \sum_{k = 1}^{n - 1} \left( \frac{n - k}{n} \right)^2
                                 | 1 - \zt^{k m} |^2    \\
  & \leq & \left( \frac{2}{n^2} \right)
        \left( \frac{4 \pi^2 m^2}{(2 n + 1)^4} \right)
                      \sum_{k = 1}^{n - 1} (n - k)^2 k^2
   \leq \left( \frac{\pi^2}{2} \right)
          \left( \frac{(n - 1) n^2 m^2}{(2 n + 1)^4} \right).
\eeqr
So
\[
\| \xi_0 - y\xi_0 \| \leq
      \frac{\pi n m}{(2 n + 1)^2} \sqrt{ \frac{n - 1}{2}}.
\]

Next, we estimate $\| \xi_0 - z \xi_0 \|.$ The components of
$\xi_0 - z \xi_0$ for which one of $(\xi_0)_j$ and $(z \xi_0)_j$
is zero and the other is not are
$\pm \frac{1}{n}, \pm \frac{2}{n}, \dots, \pm  \frac{m}{n},$
each occurring once. The sum of their squares is, using
$k^2 + (m - k)^2 \leq m^2$ for $1 \leq k \leq m,$ equal to
\[
2 \sum_{k = 1}^{m} \frac{k^2}{n^2}
  = \frac{2 m^2}{n^2} + \sum_{k = 1}^{m - 1} \frac{k^2 + (m - k)^2}{n^2}
 \leq (m + 1) \frac{m^2}{n^2}.
\]
There are $2 n - 1 - m$ other nonzero components, all of absolute
value at most $m/n,$ so
\[
\| \xi_0 - z \xi_0 \| \leq
    \sqrt{ (2 n - 1 - m) \frac{m^2}{n^2} + (m + 1) \frac{m^2}{n^2} }
    = m \sqrt{ \frac{2}{n} }.
\]

Using $\| \xi_0 \| \geq \sqrt{2n / 3},$ we therefore obtain
\[
\| \xi - y \xi \| \leq
      \frac{\pi n m}{(2 n + 1)^2} \sqrt{ \frac{n - 1}{2}}
                             \sqrt{ \frac{3}{2 n} }
             \leq \left( \frac{ \pi \sqrt{3} }{4} \right)
                                \left( \frac{m}{2 n + 1} \right)
    = \frac{ \pi \sqrt{3} }{4} \te^{1/2}
\]
and
\[
\| \xi - z \xi \| \leq
         m \sqrt{ \frac{2}{n} } \sqrt{ \frac{3}{2 n} }
       = 2 \sqrt{3} \left(1 + \frac{1}{2n} \right) \te^{1/2}.
\]
If $n$ is sufficiently large, it follows that
\[
\| \xi - y \xi \| \leq \xfrac{3}{2} \te^{1/2}, \,\,\,\,\,\,
\| \xi - z \xi \| \leq \xfrac{7}{2} \te^{1/2}, \andeqn
\| y \xi - z \xi \| \leq 5 \te^{1/2}.
\]

Using Lemma 6.4
on the real parts of the scalar products,
and taking $\af,$ $\bt,$ $\gm \geq 0,$ we calculate
\beqr
\| (\af \cdot 1 + \bt y + \gm z) \xi \|^2
  & = & \af^2 + \bt^2 + \gm^2 +
        2 \af \bt \left( 1 - \xfrac{1}{2} \| \xi - y \xi \|^2 \right) \\
 &  & \mbox{} +
      2 \af \gm \left( 1 - \xfrac{1}{2} \| \xi - z \xi \|^2 \right) +
     2 \bt \gm \left( 1 - \xfrac{1}{2} \| y \xi - z \xi \|^2 \right)  \\
 & \geq & (\af + \bt + \gm)^2 \left( 1  - \xfrac{25}{2} \te \right)
\eeqr
for large $n.$ Since $\| \xi \| = 1,$ this shows that, for
$\af,$ $\bt,$ $\gm \geq 0$ and $n$ large,
\[
\af +  \bt + \gm = 1
  \,\,\,\,\,\, {\mathrm{implies}} \,\,\,\,\,\,
   \| \af \cdot 1 + \bt y + \gm z \| \geq
       \left(1 - \xfrac{25}{2} \te \right)^{1/2}.
\]
\end{pff}

\begin{lem}   
Let $A$ and $B$ be unital \ca s, with $B$ nuclear. Let $E \subset A$
be a \fd\  operator system, and let $E_0 \subset E$ be a subspace
such that $1 \in E_0$ and $E_0 + E_0^* = E.$ If $\ep > 0$ and
$T : E \to B$ is a unital selfadjoint
linear map with $\| T |_{E_0} \|_{\cb} < 1 + \ep,$ then there exists a
\ucp\  map $S : A \to B$ such that $\| S |_{E_0} - T |_{E_0} \| < \ep.$
\end{lem}

\begin{pff}
Choose $\dt > 0$ such that $\| T |_{E_0} \|_{\cb} < 1 + \ep - \dt.$
Since $B$ is nuclear, there are $n$ and \ucp\  maps
$Q_1 : B \to M_n$ and $Q_2 : M_n \to B$ such that
$\| Q_2 \circ Q_1 |_{T (M)} - \id_B |_{T (M)} \| < \dt.$ We have
$\| Q_1 \circ T |_{E_0} \|_{\cb} < 1 + \ep - \dt,$ so
Wittstock's generalization of the Arveson extension theorem
(Theorem 7.2 of \cite{Pl}) yields a linear map $Q_0 : E \to M_n$ such
that $\| Q_0 \|_{\cb} < 1 + \ep - \dt$ and
$Q_0 |_{E_0} = Q_1 \circ T |_{E_0}.$
Set $Q (a) = \half (Q_0 (a) + Q_0 (a^*)^*).$ Then $Q : E \to M_n$ is
linear, selfadjoint, and satisfies $\| Q\|_{\cb} < 1 + \ep - \dt$ and
$Q |_{E_0} = Q_1 \circ T |_{E_0}.$ Lemma 1.9
provides a \ucp\  map $S_0 : A \to M_n$ such that
$\| S_0 |_E - Q \| < \ep - \dt.$ Then $S = Q_2 \circ S_0 : A \to B$
is  a \ucp\  map such that $\| S |_{E_0} - T |_{E_0} \| < \ep.$
\end{pff}

\begin{thm}   
The function $\rh_0$ of Notation 6.1
satisfies $\rh_0 (\te_1, \te_2) \leq \xfrac{25}{4} | \te_1 - \te_2|$
for all $\te_1,$ $\te_2 \in \R.$
\end{thm}

\begin{pff}
By Lemmas 6.2
and 6.3,
$\rh_0$ is a continuous
translation invariant pseudometric. We therefore easily see that it
is enough to show that for every $\ep > 0$ there is $\dt > 0$ and
a dense subset $S$ of $(0, \dt)$ such that
$\rh_0 (\te, 0) \leq (\frac{25}{4} + \ep ) \te$ for $\te \in S.$

Take $\ep_0 = \frac{2}{25} \ep,$ and choose $\dt_0 > 0$ such that
$( 1 - r)^{- 1/2} < 1 + \half (1 + \ep_0) r$ for $r \in (0, \dt_0).$
Take $\dt = \frac{2}{25} \dt_0$ and take
\[
S = (0, \dt) \cap
        \left\{ \frac{m^2}{(2 n + 1)^2} : m, \, n \in \N \right\}.
\]
Then for $\te \in S,$ Lemma 6.5
shows that the map
$T_{\te} : E_0 (\te) \to E_0 (0)$ of that lemma satisfies
\[
\| T_{\te} \|_{\cb} \leq \left(1 - \xfrac{25}{2} \te \right)^{- 1/2}
     < 1 + (\xfrac{25}{4} + \ep ) \te.
\]
Now apply Lemma 6.6,
with $E_0 = E_0 (\te),$ $E = E (\te),$
and $T$ defined by
\[
T (1) = 1,\,\,\,\,\,\,
T (u (\te)) = u (0), \,\,\,\,\,\, T (u (\te)^*) = u (0)^*, \,\,\,\,\,\,
T (v (\te)) = v (0), \andeqn T (v (\te)^*) = v (0)^*,
\]
so that $T |_{E_0} = T_{\te}.$ This gives a \ucp\  map
$R_0 : A (\te) \to A (0)$ such that
$\| R_0 |_{E_0 (\te)} - T_{\te} \| < (\xfrac{25}{4} + \ep ) \te.$
In particular, $R = \io_0 \circ R_0 : A (\te) \to \OA{2}$ is a
\ucp\  map satisfying
\[
\| R (u (\te)) - \io_0 (u (0)) \| < (\xfrac{25}{4} + \ep ) \te
\andeqn
\| R (v (\te)) - \io_0 (v (0)) \| < (\xfrac{25}{4} + \ep ) \te.
\]
So $\rh_0 (\te, 0) \leq (\frac{25}{4} + \ep ) \te.$
\end{pff}

\begin{cor}   
The field of rotation algebras, with the unitary sections defined by
the standard generators, is $(X_0, \rh)$-embeddable in $\OA{2}$ for
a one point space $X_0$ and with $\rh (r) = 28 r^{1/2}.$
\end{cor}

\begin{pff}
This follows from Remark 4.16
and the inequality
$11 \sqrt{25/4} < 28.$
\end{pff}

\begin{thm}   
There exists a unital injective representation $\ph$ of the
continuous field of rotation algebras, regarded as defined over $\R,$
which is periodic in the sense that $\ph_{\te + 1} = \ph_{\te}$ for
all $\te \in \R,$ and for which there is a constant $C$ such that
for all $\te_1,$ $\te_2 \in \R,$
\[
\max ( \| \ph_{\te_1} (u (\te_1)) - \ph_{\te_2} (u (\te_2)) \|, \,
       \| \ph_{\te_1} (v (\te_1)) - \ph_{\te_2} (v (\te_2)) \| )
  < C | \te_1 - \te_2 |^{1/2}.
\]
Moreover, $C$ can be chosen smaller than $840,000 = 8.4 \cdot 10^5.$
\end{thm}

\begin{pff}
We describe the changes that must be made to
the proofs of Lemma 5.6
and Theorem 5.10.
Take the spaces $X$ and $Y$ of Lemma 5.6
to be one point spaces
(as in the proof of Theorem 5.10),
and use $\R$ in place of
$[0, 1].$ Choose a single injective unital
\hm\  $\af : A (0) \to \OA{2},$ and take $\ph_n = \af$ for every
integer $n.$ We carry out the rest of the construction on $[0, 1],$
repeating each step using
periodicity in each $[n, n + 1].$ (The only reason for not
restricting to $[0, 1]$ at the beginning is to ensure that we
obtain the distance estimate of the theorem for, say,
$\te_1 < 1 < \te_2.$)

We will take all sectional distances $\ds (\ph_{\te_1}, \ph_{\te_2})$
with respect to the continuous sections $u$ and $v.$
By Theorem 6.7,
we have
$\rh_0 (\te_1, \te_2) < \frac{25}{4} | \te_1 - \te_2 |.$
So, in the proof of Theorem 5.10,
we take $\af = 1.$
We follow the modification described for this value of $\af$ in Remark
5.11.

Construct injective unital \hm s $\ph_{\te} : A(\te) \to \OA{2}$
for $\te$ in the set
\[
S_0 = \{ j/n^k : k \geq 0, \,\, 0 \leq j \leq n^k\}
\]
as in the proof of Lemma 5.6,
with $n_1,$ $n_2, \dots,$ and $n'$ all equal to $90^2$
(see Remark 5.11),
and extend over $S = \{ j/n^k : k \geq 0, \,\, j \in \Z \}$
by periodicity. As in Remark 5.11, if $\te_1,$ $\te_2 \in S$ then
\[
\ds (\ph_{\te_1}, \ph_{\te_2} ) \leq
    330,000 \left( \xfrac{25}{4} \right)^{1/2} | \te_1 - \te_2 |^{1/2}
  \leq 840,000 | \te_1 - \te_2 |^{1/2}.
\]
This is true for all $\te_1,$ $\te_2 \in S,$ but clearly extends by
continuity to all $\te_1,$ $\te_2 \in \R.$
\end{pff}

\begin{cor}   
There exist continuous functions $u,$ $v : S^1 \to U (\OA{2})$
such that:

(1) $u (\zt) v (\zt) = \zt v (\zt) u (\zt)$ for all $\zt \in S^1.$

(2) $C^* (u (\zt), v (\zt))$ is isomorphic to the universal
\ca\  on unitaries $u$ and $v$ satisfying $u v = \zt v u.$

(3) There is a constant $C$ such that for all $\zt_1,$ $\zt_2 \in S^1,$
we have
\[
\| u (\zt_1) - u (\zt_2) \| \leq C | \zt_1 - \zt_2 |^{1/2}  \andeqn
\| v (\zt_1) - v (\zt_2) \| \leq C | \zt_1 - \zt_2 |^{1/2}.
\]

Moreover, $C$ can be chosen less than $420,000.$
\end{cor}

\begin{pff}
This follows from the theorem, because if $\zt_1,$ $\zt_2 \in S^1,$
then there exist $\te_1,$ $\te_2 \in \R$ such that
\[
\exp (2 \pi i \te_1) = \zt_1, \,\,\,\,\,\,
\exp (2 \pi i \te_2) = \zt_2, \andeqn
| \te_1 - \te_2 | \leq \xfrac{1}{4} | \zt_1 - \zt_2 |.
\]
It follows that
$| \te_1 - \te_2 |^{1/2} \leq \xfrac{1}{2} | \zt_1 - \zt_2 |^{1/2}.$
\end{pff}

\begin{rmk}    
No better exponent is possible in Theorem 6.9
or in
Corollary 6.10,
because Proposition 4.6 of \cite{HR} shows
no better exponent is possible even for representations on a
Hilbert space. This implies that the square root in
Proposition 4.13
can't be removed, and also provides a (rather indirect) proof that the
exponent $\half$ in Lemma 1.12
can't be improved. Moreover,
in the construction of embeddings in $\OA{2},$ the exponent $\half$
in Lemma 1.12
can't be evaded by using some other proof.
\end{rmk}


\begin{thebibliography}{BEEK}

\bibitem[AAP]{AAP} C. A. Akemann, J. Anderson, and G. K. Pedersen,
{\em Excising states of \ca s,}
Canad. J. Math. {\bf 38}(1986), 1239-1260.

\bibitem[AP1]{AP1} C. A. Akemann and G. K. Pedersen,
{\em Ideal perturbations of elements in \ca s},
Math. Scand. {\bf 41}(1977), 117-139.

\bibitem[AP2]{AP2} C. A. Akemann and G. K. Pedersen,
{\em Central sequences and inner derivations of \sep\  \ca s,}
Amer. J. Math. {\bf 101}(1979), 1047-1061.

\bibitem[Ar]{Ar}  W. Arveson, {\em Notes on extensions of \ca s},
Duke Math. J. {\bf 44}(1977), 329-355.

\bibitem[BK]{BK} B. Blackadar and E. Kirchberg, {\em Generalized
inductive limits of finite dimensional \ca s}, preprint.

\bibitem[BKR]{BKR} B. Blackadar, A. Kumjian, and M. R\o rdam,
{\em Approximately central matrix units and the structure of
non-commutative tori}, $K$-Theory {\bf 6}(1992), 267-284.

\bibitem[Bl]{Bl} E. Blanchard, {\em Subtriviality of
continuous fields of nuclear C*-algebras}, preprint.

\bibitem[BEEK]{BEEK} O. Bratteli, G. A. Elliott, D. E. Evans, and
A. Kishimoto, {\em On the classification of C*-algebras of real
rank zero III: The infinite case}, in preparation.

\bibitem[BKRS]{BKRS} O. Bratteli, A. Kishimoto, M. R\o rdam, and E.
St\o rmer, {\em The crossed product of a UHF algebra by a shift},
Ergod. Th. Dynam. Sys. {\bf 13}(1993), 615-626.

\bibitem[CE1]{CE1} M.-D. Choi and E. G. Effros, {\em Nuclear \ca s and
the approximation property}, Amer. J. Math. {\bf 100}(1978), 61-79.

\bibitem[CE2]{CE2} M.-D. Choi and E. G. Effros, {\em Injectivity and
operator spaces}, J. Functional Analysis {\bf 24}(1977), 156-209.

\bibitem[Cn1]{Cu0} J. Cuntz, {\em Simple C*-algebras generated by
isometries}, Comm. Math. Phys. {\bf 57}(1977), 173-185.

\bibitem[Cn2]{Cu1} J. Cuntz, {\em K-theory for certain
\ca s}, Ann. Math. {\bf 113}(1981), 181-197.

\bibitem[Dx]{Dx} J. Dixmier, {\em C*-Algebras}, North-Holland,
Amsterdam, New York, Oxford, 1977.

\bibitem[EH]{EH} E. G. Effros and U. Haagerup, {\em Lifting
problems and local reflexivity for \ca s}, Duke Math. J.
{\bf 52}(1985), 103-128.

\bibitem[EfR]{EfR} E. G. Effros and J. Rosenberg,
{\em C*-algebras with approximately inner flip}, Pacific J. Math.
{\bf 77}(1978), 417-443.

\bibitem[El]{El1} G. A. Elliott, {\em The classification
of \ca s of real rank zero}, J. reine angew. Math., {\bf 443}(1993),
179-219.

\bibitem[ElR]{ElR} G. A. Elliott and M. R\o rdam,
{\em Classification of certain infinite simple C*-algebras II},
Comment. Math. Helvetici {\bf 70}(1995), 615-638.

\bibitem[HR]{HR} U. Haagerup and M. R\o rdam, {\em Perturbations of
the rotation \ca s and of the Heisenberg commutation relation},
Duke Math. J. {\bf 77}(1995), 627-656.

\bibitem[HK]{HK} K. H. Hofmann and K. Keimel,
{\em Sheaf theoretic concepts in analysis: bundles and sheaves of
Banach spaces, Banach $C (X)$-modules}, pages 415-441 in:
{\em Applications of Sheaves}, Lecture Notes in Math.\ no.\  753,
Springer-Verlag, New York, Heidelberg, Berlin, 1979.

\bibitem[JP]{JP} M. Junge and G. Pisier, {\em Bilinear forms on exact
operator spaces and $B (H) \otimes B (H)$}, Geometric and Functional
Analysis {\bf 5}(1995), 329-363.

\bibitem[KR]{KR} R. V. Kadison and J. R. Ringrose,
{\em Fundamentals of the Theory of Operator Algebras, Volume II},
Academic Press, New York, London, Paris, San Diego, San Francisco,
S\~{a}o Paulo, Sydney, Tokyo, Toronto, 1983.

\bibitem[Ks]{Ks1} G. G. Kasparov, {\em Hilbert C*-modules: theorems
of Stinespring and Voiculescu},
J. Operator Theory {\bf 4}(1980), 133-150.

\bibitem[Kr1]{Kr0} E. Kirchberg, {\em C*-nuclearity implies CPAP},
Math. Nachr. {\bf 76}(1977), 203-212.

\bibitem[Kr2]{Kr0.5} E. Kirchberg, {\em On non-semisplit extensions,
tensor products and exactness of group \ca s}, Invent. Math.
{\bf 112}(1993), 449-489.

\bibitem[Kr3]{Kr1} E. Kirchberg, {\em On subalgebras of the CAR
algebra}, J. Functional Analysis {\bf 129}(1995), 35-63.

\bibitem[Kr4]{Kr2} E. Kirchberg, {\em Commutants of unitaries in
UHF algebras and functorial properties of exactness}, J. reine angew.
Math. {\bf 452}(1994), 39-77.

\bibitem[Kr5]{Kr3} E. Kirchberg, {\em  The classification of purely
infinite C*-algebras using Kasparov's theory},
preliminary preprint (3rd draft).

\bibitem[KW]{KW} E. Kirchberg and S. Wassermann, {\em Operations on
continuous bundles of \ca s}, Math. Annalen {\bf 303}(1995), 677-697.

\bibitem[Ln1]{Ln} H. Lin, {\em Exponential rank of C*-algebras with
      real rank zero and the Brown-Pedersen conjectures},
       J. Functional Analysis {\bf 114}(1993), 1-11.

\bibitem[Ln2]{Ln2} H. Lin, {\em Almost commuting unitaries
  in purely infinite simple C*-algebras}, Math. Annalen
{\bf 303}(1995), 599-616.

\bibitem[Ln3]{Ln3} H. Lin, {\em Almost commuting unitaries
and classification of purely infinite simple C*-algebras}, preprint.

\bibitem[LP1]{LP0} H. Lin and N. C.  Phillips, {\em Classification of
direct limits of even Cuntz-circle algebras}, Memoirs Amer. Math. Soc.
no. 565 (1995).

\bibitem[LP2]{LP} H. Lin and N. C.  Phillips, {\em Approximate
unitary equivalence of \hm s from $\OI$}, J. reine angew. Math.
{\bf 464}(1995), 173-186.

\bibitem[Lr]{Lr} T. A. Loring,
{\em C*-algebras generated by stable relations},
J. Functional Analysis {\bf 112}(1993), 159-203.

\bibitem[OZ]{OZ} C. L. Olsen and W. R. Zame, {\em Some C*-algebras with
a single generator}, Trans. Amer. Math. Soc. {\bf 215}(1976), 205-217.

\bibitem[Pl]{Pl} V. I Paulsen, {\em Completely Bounded Maps and
Dilations,} Pitman Research Notes in Math. no. 146, Longman Scientific
and Technical, Harlow, Britain, 1986.

\bibitem[Pd]{Pd} G. K. Pedersen, {\em C*-Algebras and their
Automorphism Groups,} Academic Press, London, New York, San Francisco,
1979.

\bibitem[Ph1]{Ph} N. C.  Phillips, {\em Approximation by unitaries with
       finite spectrum in purely infinite simple C*-algebras},
    J. Funct. Anal. {\bf 120}(1994), 98-106.

\bibitem[Ph2]{Ph2} N. C.  Phillips, {\em A classification theorem for
nuclear purely infinite simple \ca s}, preprint.

\bibitem[Rf]{Rf} M. Rieffel, {\em Continuous fields of \ca s coming
from group cocycles and actions}, Math. Ann. {\bf 283}(1989),
631-643.

\bibitem[Rn]{Rn} J. Ringrose, {\em Exponential length and
exponential rank of \ca s}, Proc. Royal Soc. Edinburgh (Sect. A)
{\bf 121}(1992), 55-71.

\bibitem[Rr1]{Rr1} M. R\o rdam, {\em Classification of inductive limits
of Cuntz algebras}, J. reine angew. Math. {\bf 440}(1993), 175-200.

\bibitem[Rr2]{Rr2} M. R\o rdam, {\em Classification of
Cuntz-Krieger algebras}, K-Theory {\bf 9}(1995), 31-58.

\bibitem[Rr3]{Rr3} M. R\o rdam, {\em A short proof of Elliott's
Theorem: $\OA{2}\otimes \OA{2}\cong \OA{2}$},
C. R. Math. Rep. Acad. Sci. Canada {\bf 16}(1994), 31-36.

\bibitem[Rr4]{Rr4} M. R\o rdam,
{\em Classification of certain infinite simple C*-algebras},
J. Funct. Anal. {\bf 131}(1995), 415-458.

\bibitem[Sc]{Sch} C. Schochet, {\em Topological methods for
     \ca s II: geometric resolutions and the  K\"{u}nneth formula},
Pacific J. Math. {\bf 98}(1982), 443-458.

\bibitem[Tk]{Tk}  M. Takesaki, {\em Theory of Operator Algebras I},
Springer-Verlag, New York, Heidelberg, Berlin, 1979.

\bibitem[Vc]{Vc} D. Voiculescu, {\em A note on quasidiagonal \ca s
and homotopy}, Duke Math. J. {\bf 62}(1991), 267-271.

\bibitem[Ws]{Ws} S. Wassermann, {\em Exact C*-Algebras and Related
Topics}, Lecture Notes Series no.~19, GARC, Seoul National University,
1994.

\bibitem[Zh1]{Zh} S. Zhang, {\em A property of purely infinite simple
\ca s}, Proc. Amer. Math. Soc. {\bf 109}(1990), 717-720.

\bibitem[Zh2]{Zh1} S. Zhang, {\em On the homotopy type of the unitary
group and the Grassmann space of \pisca s}, K-Theory, to appear.

\end{thebibliography}
\end{document}